\documentclass[a4paper,12pt,twoside]{book}

\usepackage[T1]{fontenc}
\usepackage[latin1]{inputenc}

\usepackage{lmodern}
\usepackage{microtype}

\usepackage{fancyvrb}
\usepackage{pdfpages}
\usepackage{tabularx}
\usepackage{array}

\usepackage{rotating, multirow}
\usepackage{yhmath}

\usepackage{graphicx}
\usepackage[ED=MITT-STICIA, Ets=INP]{tlsflyleaf}

\usepackage{float}
\usepackage[pdfpagelabels, pagebackref=true, hyperfootnotes=false, bookmarksnumbered=true]{hyperref}

\definecolor{amblu}{RGB}{51, 102, 153}
\definecolor{amred}{RGB}{230,0,0}
\definecolor{amgreen}{RGB}{96,107,47}

\usepackage{amsmath,amsfonts,amssymb}
\usepackage{fixltx2e}
\usepackage{algorithm}
\usepackage{algpseudocode}
\MakeRobust{\Call}
\MakeRobust{\If}

\usepackage{enumitem}
\usepackage{bm}
\usepackage{color}
\usepackage{cancel}
\usepackage{pdflscape}
\usepackage{subfig}

\usepackage{pgfplots}
\usetikzlibrary{shapes,arrows}

\usepackage[toc, nonumberlist]{glossaries}
\usepackage{minitoc}
\dominitoc
\doparttoc

\hypersetup{
	pdftitle = {Hybridization of interval methods and evolutionary algorithms for solving difficult optimization problems},
	pdfkeywords = {global optimization, interval analysis, evolutionary algorithms},
	pdfborderstyle = {/S/U/W 1},
	pdfauthor = {Charlie Vanaret},
	pdfdisplaydoctitle = true
}

\newtheorem{theorem}{Theorem}
\newtheorem{example}{Example}

\newtheorem{definition}{Definition}
\newtheorem{remark}{Remark}

\newcommand{\eqdef}{:=}
\newcommand{\resp}{respectively}
\newcommand{\inout}{\textbf{in-out}}

\newcommand{\hull}{\mathop{\square}}

\graphicspath{{figures/}}

\makeglossaries

\title{Hybridization of interval methods and evolutionary algorithms for solving difficult optimization problems}
\author{Charlie Vanaret}
\defencedate{January 27, 2015}
\lab{IRIT-APO (UMR 5505)}
\nboss{2}
\makesomeone{boss}{1}{Nicolas Durand}{}{}
\makesomeone{boss}{2}{Jean-Baptiste Gotteland}{}{}
\nreferee{0}
\njudge{8}
\makesomeone{judge}{1}{Nicolas Durand}{ENAC}{Directeur de thèse}
\makesomeone{judge}{2}{Jean-Baptiste Gotteland}{ENAC}{Co-encadrant de thèse}
\makesomeone{judge}{3}{El-Ghazali Talbi}{Université de Lille}{Rapporteur}
\makesomeone{judge}{4}{Gilles Trombettoni}{Université de Montpellier}{Rapporteur}
\makesomeone{judge}{5}{Jean-Marc Alliot}{IRIT}{Examinateur}
\makesomeone{judge}{6}{Jin-Kao Hao}{Université d'Angers}{Examinateur}
\makesomeone{judge}{7}{Thomas Schiex}{INRA Toulouse}{Examinateur}
\makesomeone{judge}{8}{Marc Schoenauer}{INRIA Saclay}{Examinateur}

\usepackage{a4wide}

\begin{document}
\makeflyleaf

\pagenumbering{roman}
\setcounter{page}{1}

\section*{Abstract}

Reliable global optimization is dedicated to finding a global minimum in the presence of rounding
errors. The only approaches for achieving a numerical proof of optimality in global optimization are
interval-based methods that interleave branching of the search-space and
pruning of the subdomains that cannot contain an optimal solution. The exhaustive interval branch
and bound methods have been widely studied since the 1960s and have benefitted from the development
of refutation methods and filtering algorithms, stemming from the interval analysis and interval
constraint programming communities. It is of the utmost importance: i) to compute sharp enclosures
of the objective function and the constraints on a given subdomain; ii) to find a good approximation
(an upper bound) of the global minimum.

State-of-the-art solvers are generally \textit{integrative} methods, that is they embed local
optimization algorithms to compute a good upper bound of the global minimum over each subspace. In
this document, we propose a \textit{cooperative} framework in which interval methods
cooperate with evolutionary algorithms. The latter are stochastic algorithms in which a population
of individuals (candidate solutions) iteratively evolves in the search-space to reach satisfactory
solutions. Evolutionary algorithms, endowed with operators that help individuals escape from local
minima, are particularly suited for difficult problems on which traditional methods struggle to
converge.

Within our cooperative solver Charibde, the evolutionary algorithm and the interval-based
algorithm run in parallel and exchange bounds, solutions and search-space via message passing. A
strategy combining a geometric exploration heuristic and a domain reduction operator prevents
premature convergence toward local minima and prevents the evolutionary algorithm from exploring
suboptimal or unfeasible subspaces. A comparison of Charibde with state-of-the-art solvers based on
interval analysis (GlobSol, IBBA, Ibex) on a benchmark of difficult problems shows that
Charibde converges faster by an order of magnitude. New optimality results are provided for five
multimodal problems, for which few solutions were available in the literature. Finally,
we certify the optimality of the putative solution to the Lennard-Jones cluster problem for five
atoms, an open problem in molecular dynamics.

\cleardoublepage
\section*{Acknowledgement}

Je dois mon entrée -- somme toute assez imprévue -- dans le monde académique à Jean-Marc Alliot.
Mes premières semaines à la DTI, à la découverte de la programmation par contraintes, m'ont
laissé entrevoir certains aspects de la recherche que je ne soupçonnais pas ; grand bien lui en a
pris. Je remercie Jean-Baptiste Gotteland pour sa confiance, sa disponibilité et la grande latitude
qu'il m'a laissée. Merci à Nicolas Durand pour sa générosité et ses encouragements. Je garde
un souvenir particulièrement émouvant de notre tentative de record de la traversée du Massif central
en TB-20.

Je tiens à remercier Gilles Trombettoni d'avoir répondu patiemment à toutes mes questions
relatives à l'implémentation des contracteurs. Je suis reconnaissant à Gilles et à El-Ghazali Talbi
d'avoir accepté de rapporter ma thèse. Je remercie mes examinateurs Thomas Schiex et Marc
Schoenauer d'avoir fait le déplacement à Toulouse, et Jin-Kao Hao de m'avoir donné l'opportunité de
faire ma première télé.

Que dire de mes "compagnons de galère" Richard Alligier et Mohammad Ghasemi Hamed...
Partager un bureau à trois n'est pas toujours chose aisée, surtout lorsque mes jeux de mots du lundi
matin valent ceux d'un vendredi après-midi. Z05 a été le théâtre de discours parfois animés,
toujours passionnés, d'échanges constructifs et de synchronisation pour la pause thé. 

Mon camarade de marave Cyril Allignol, le baryton Nicolas Barnier, David Gianazza maître
ès foncteurs, Alexandre "Jean-Michel" Gondran, Sonia Cafieri, Loïc Cellier, Brunilde
Girardet, Laureline Guys, Olga Rodionova, Nicolas Saporito et Estelle Malavolti -- pour un temps mes
voisins, entre deux valses des bureaux -- tous ont contribué à la bonne humeur et à l'ambiance
chaleureuse du bâtiment Z.

Je salue Daniel Ruiz pour sa gentillesse, et l'équipe APO pour leur accueil. Frédéric Messine et
Jordan Ninin ont aimablement répondu à mes interrogations affines. Ma gratitude va à Christine
Surly, garante de ma logistique pendant ces trois années chez Midival, et à Jean-Pierre Baritaud
pour le soutien technique lors de ma soutenance.

Merci à Fabien Bourrel, mon relecteur officiel qui n'a jamais trouvé une seule typo, alors que...
J'en profite pour passer un petit coucou à ma famille, aux Barousse, aux boys d'Hydra, à mes
camarades de promo n7, aux expatriés viennois, au TUC Escrime et au Péry.

Merci à mes parents et mon frangin d'avoir fait le déplacement lors de ma soutenance.

Désolé pour le slide numéro 2.

\cleardoublepage

\newacronym[\glslongpluralkey={floating-point units},\glsshortpluralkey={FPUs}]{FPU}{FPU}{floating-point unit}
\newacronym[\glslongpluralkey={evolutionary algorithms},\glsshortpluralkey={EAs}]{EA}{EA}{evolutionary algorithm}
\newacronym[\glslongpluralkey={genetic algorithms},\glsshortpluralkey={GAs}]{GA}{GA}{genetic algorithm}
\newacronym[]{FPA}{FPA}{floating-point arithmetic}
\newacronym[]{IA}{IA}{interval arithmetic}
\newacronym[]{AA}{AA}{affine arithmetic}
\newacronym[]{BB}{BB}{branch and bound}
\newacronym[]{IBB}{IBB}{interval branch and bound}
\newacronym[]{IBC}{IBC}{interval branch and contract}
\newacronym[]{AD}{AD}{automatic differentiation}
\newacronym[]{DE}{DE}{differential evolution algorithm}
\newacronym[]{CSP}{CSP}{constraint satisfaction problem}
\newacronym[]{NCSP}{NCSP}{numerical constraint satisfaction problem}
\newacronym[]{CID}{CID}{constructive interval disjunction}
\newacronym[]{ICP}{ICP}{interval constraint programming}
\newacronym[plural={CAS}, \glsshortpluralkey={CAS}]{CAS}{CAS}{computer algebra system}

\tableofcontents
\listoffigures
\listoftables
\listofalgorithms

\printglossary[title=Glossary,toctitle=Glossary]
\markboth{GLOSSARY}{}

\setcounter{chapter}{0}
\chapter*{Introduction}
\markboth{INTRODUCTION}{} 
\pagenumbering{arabic}
\setcounter{page}{1}

\section*{Motivation}
Numerical computations based on floating-point arithmetic may be subject to roundoff errors ; roundoff accumulation sometimes produces irrelevant results that are disastrous for critical systems (for instance, in aerospace). The only methods capable of rigorously bounding the intermediary steps of a numerical computation are based on interval analysis, a branch of numerical analysis that extends floating-point arithmetic to intervals. The ability of interval analysis to compute with sets unraveled new horizons for the global optimization community.

Reliable global optimization methods based on interval analysis, called interval branch and bound, partition the search space and discard subspaces that cannot contain an optimal solution using refutation arguments: whenever a lower bound of the range of the objective function on a subspace is larger than an upper bound of the global minimum (the objective value of any feasible point), it is numerically guaranteed that the subspace cannot contain an optimal solution.
Nowadays, cutting-edge solvers embed filtering (or contraction) operators that stem from the numerical analysis and the discrete optimization communities ; they aim at reducing the bounds of the variables without losing the optimal solution. Bisection however remains sometimes unavoidable. On account of its exponential complexity in the number of variables, one cannot hope to solve instances larger than a few dozen variables. 

Invoking exhaustive methods to solve nonconvex and highly multimodal optimization problems may seem hopeless. In this case, metaheuristics usually provide satisfactory solutions within a reasonable time, albeit with no guarantee of optimality. Among the population-based metaheuristics that maintain a population of individuals (a set of candidate solutions), evolutionary algorithms mimic mechanisms inspired by nature, in order to guide a random walk towards good solutions. Because they embed operators that help escape from local minima, metaheuristics are widely used in the optimization community when other methods fail to converge.

In general, exhaustive solvers integrate a local method to the branch bound scheme in order to compute approximations (upper bounds) of the global minimum. Very few combine an exact global method (branch and bound) and a stochastic method (such as an evolutionary algorithm) ; the existing approaches are essentially sequential (one method runs after the other) or integrative (one is embedded within the other).

\section*{Adopted approach and contributions}
\begin{figure}[htbp]
\centering
\def\svgwidth{0.85\columnwidth}
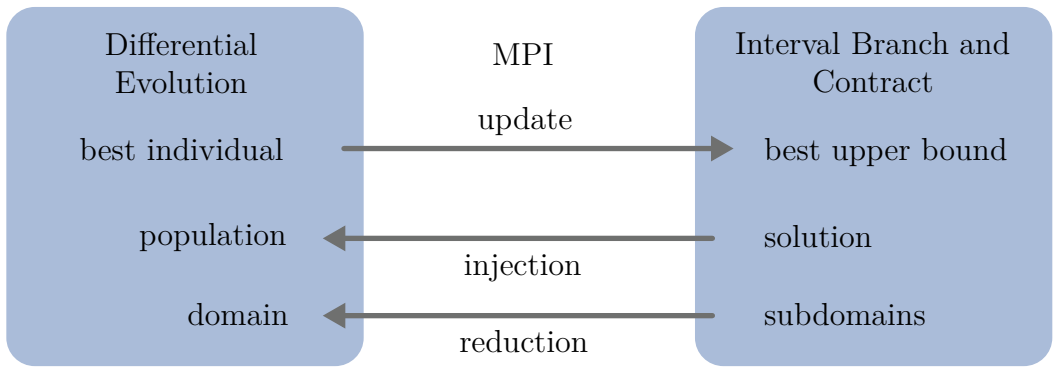
\caption{Cooperation scheme of Charibde}
\label{fig:intro-charibde}
\end{figure}

Branch and bound methods require good feasible solutions (whose objective values are upper bounds of the global minimum), and accurate enclosures of the objective function and the constraints on a given subspace.
In this document, we introduce a cooperative framework that combines state-of-the-art interval methods and evolutionary algorithms. Within our hybrid solver Charibde, an interval branch and contract method and a differential evolution algorithm run in parallel and exchange bounds, solutions and domain using message passing (Figure \ref{fig:intro-charibde}).

The evolutionary algorithm quickly explores the search space in the search of a satisfactory feasible solution. Its evaluation is sent to the interval method in order to intensify the pruning of the infeasible and suboptimal subspaces of the search space.
Whenever the interval method finds a points improves the best known solution, it is injected into the population of the evolutionary algorithm in order to avoid premature convergence towards local minima.
A combination of a novel exploration strategy and a periodic domain reduction of the differential evolution algorithm avoids the generation of individuals within infeasible or suboptimal subspaces.
Cutting-edge filtering operators (contractors) reduce the bounds of the variables on each subspace by discarding values that are inconsistent with respect to the constraints. Some exploit the syntax tree of a single constraint at the time, others convexify the problems and consider all the constraints simultaneously. Combining contraction and automatic differentation produces tighter enclosures of the partial derivatives of the functions, which in turns makes the first-order refutation test more efficient.

Charibde has proven competitive with cutting-edge reliable interval-based solvers and unreliable NLP solvers: 
it outperforms GlobSol, IBBA and Ibex by an order of magnitude on a subset of difficult COCONUT\footnote{The COCONUT benchmark is available at \url{http://www.mat.univie.ac.at/~neum/glopt/coconut/Benchmark/Benchmark.html}} problems.
We provide new optimality results for five multimodal problems (Michalewicz, Sine Wave Sine Envelope, Eggholder, Keane, Rana) for which few solutions, even approximate, are known.
Finally, we present the first numerical proof of optimality for the open Lennard-Jones cluster problem with five atoms. We show that interval-based solvers do not converge within reasonable time, and that NLP solvers BARON and Couenne provide numerically erroneous results that cannot be trusted. 

\section*{Organization of the document}

This document is composed of seven chapters.
Chapter \ref{chap:opti} exposes the mathematical context of the study and introduces the basic theory of nonlinear optimization, the first-order optimality conditions and resolution methods for convex and nonconvex problems.
Evolutionary algorithms, including genetic algorithms and differential evolution algorithms, are presented in Chapter \ref{chap:evolutionary-algorithms}.
Chapter \ref{chap:ia} introduces interval methods, their application to global optimization and automatic differentiation techniques.
Chapter \ref{chap:contractors} extends the previous chapter and compares filtering algorithms (also known as contractors) for interval domains. They stem from the numerical analysis and the constraint programming communities.
Our reliable solver Charibde is described in Chapter \ref{chap:charibde}. We explain its architecture in detail and the advanced techniques devised to exploit the combination between interval methods and metaheuristics.
In Chapter \ref{chap:lj}, we close the open Lennard-Jones cluster problem with five atoms by providing the first numerical proof of optimality of the solution.

\adjustmtc[+1]

\chapter{Nonlinear optimization}
\label{chap:opti}

\minitoc

Optimization is the discipline that determines in an analytical or numerical fashion the
best solution to a problem, with respect to a certain criterion. It is fundamental for
solving countless problems in industry, economics and physics in order to reduce costs or
computing time. The quality of the solution computed by an optimization process generally
depends upon the model used to approximate real data, and the resolution method.
Section \ref{sec:optimization} introduces unconstrained and constrained optimization, and
necessary conditions of optimality. Optimization techniques are
mentioned in Section \ref{sec:optimization-techniques}.

\section{Optimization theory}\label{sec:optimization}

A continuous optimization problem can be written in standard form:
\begin{equation}
\min_{\bm{x} \in \mathbb{R}^n} f(\bm{x}) \text{ subject to } \bm{x} \in D
\label{eq:optimization-problem}
\end{equation}
$\bm{x} = (x_1, \ldots, x_n)$ are \emph{decision variables}. $f:
D \subset \mathbb{R}^n \rightarrow \mathbb{R}$ is the \emph{objective function} and $D$ is
the \emph{feasible set}. Any point $\bm{x} \in \mathbb{R}^n$ that belongs to $D$ is called a
\emph{feasible point}. Note that maximizing a function $f$ is equivalent to minimizing
the function $-f$:
\begin{equation}
\max_{\bm{x} \in D} f(\bm{x}) = -\min_{\bm{x} \in D} (-f(\bm{x}))
\end{equation}

\subsection{Local and global minima}
Solving an optimization problem boils down to seeking a local or global minimum (Definition
\ref{def:min-local-global}) of a function, and (or) the set of corresponding minimizers.

\begin{definition}[Minima and minimizers]
Let $\bm{x} \in \mathbb{R}^n$.
\begin{itemize}
\item $\bm{x}$ is a \emph{local minimizer} of $f$ in $D \subset \mathbb{R}^n$
if $\bm{x} \in D$ and there exists an open neighborhood $N$ of $\bm{x}$ such that:
\begin{equation}
\forall \bm{y} \in D \cap N, \quad f(\bm{x}) \le f(\bm{y})
\end{equation}
$f(\bm{x})$ is a \emph{local minimum} of $f$ in $D$ ;
\item $\bm{x}$ is a \emph{global minimizer} of $f$ in $D \subset \mathbb{R}^n$
if $\bm{x} \in D$ and:
\begin{equation}
\forall \bm{y} \in D, \quad f(\bm{x}) \le f(\bm{y})
\end{equation}
$f(\bm{x})$ is a \emph{global minimum} of $f$ in $D$.
\end{itemize}
\label{def:min-local-global}
\end{definition}
Local and global maximizers and maxima are defined likewise.

\begin{figure}[htbp]
\centering
\def\svgwidth{0.6\columnwidth}
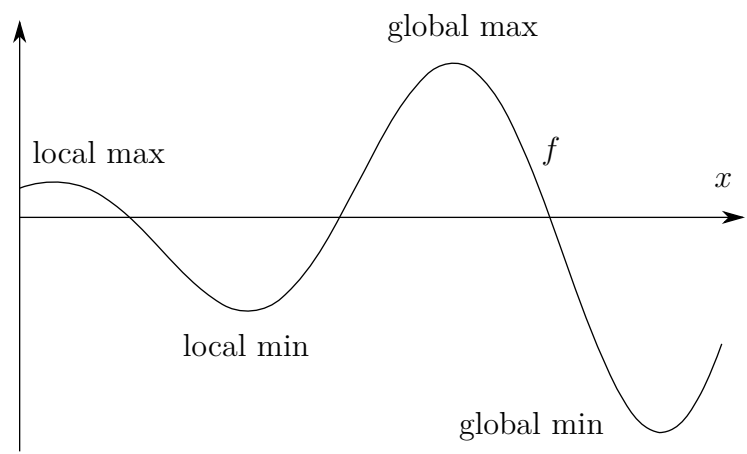
\caption{Local and global extrema}
\label{fig:extrema}
\end{figure}

Figure \ref{fig:extrema} illustrates local and global extrema (minima and maxima)
of a continuous univariate function.

\subsection{Existence of a minimum}
The extreme value theorem (Theorem \ref{th:weierstrass}) states that Problem
\ref{eq:optimization-problem} has a minimum when $f$ is continuous and $D$ is a non-empty compact set.

\begin{theorem}[Extreme value theorem]
\label{th:weierstrass}
A continuous function $f: D \rightarrow \mathbb{R}$, where $D \subset \mathbb{R}^n$ is a
non-empty compact set, attains a maximum and a minimum. In particular, there exists
$\bm{x} \in D$ such that
\begin{equation}
\forall \bm{y} \in D, \quad f(\bm{x}) \le f(\bm{y})
\end{equation}
\end{theorem}

\subsection{Unconstrained optimization}
In this section, we characterize the points that minimize a function $f:
\mathbb{R}^n \rightarrow \mathbb{R}$:
\begin{equation}
\min_{\bm{x} \in \mathbb{R}^n} f(\bm{x})
\end{equation}
where $f$ is assumed at least differentiable.

$\nabla f$ is the gradient of a differentiable function $f$. $\nabla^2 f$ is the Hessian matrix
of a twice-differentiable function $f$, whose element at row $i$ and column $j$ is $\frac{\partial^2
f}{\partial x_i\partial x_j}$.

Theorem \ref{th:necessary-conditions} introduces necessary conditions of optimality that
characterize the local minima of $f$.

\begin{definition}[Stationary point]
\label{def:stationary-point}
A point $\bm{x} \in \mathbb{R}^n$ is a \emph{stationary point} of a differentiable function
$f: \mathbb{R}^n \rightarrow \mathbb{R}$ if:
\begin{equation}
\nabla f(\bm{x}) = 0
\end{equation}
\end{definition}

\begin{theorem}[Necessary conditions of optimality]
\label{th:necessary-conditions}
Let $\bm{x}^* \in \mathbb{R}^n$ be a local minimum of a differentiable function $f:
\mathbb{R}^n \rightarrow \mathbb{R}$. Then:
\begin{enumerate}
\item $\bm{x}^*$ is a stationary point (Definition \ref{def:stationary-point}) of $f$
(first-order condition) ;
\item if $f$ is twice differentiable in an open neighborhood of $\bm{x}^*$, then
$\nabla^2 f(\bm{x}^*)$ is positive semi-definite (second-order condition).
\end{enumerate}
\end{theorem}

\begin{remark}
The stationarity of a local minimum is a \emph{necessary} but not \emph{sufficient} condition:
the function $f(x) = x^3$ has a stationary point $x = 0$ that verifies the second-order condition,
however $x = 0$ is not a local minimum.
\end{remark}

Although not sufficient, necessary conditions may help select potential local minima. Theorem
\ref{th:sufficient-conditions} states sufficient conditions of optimality.

\begin{theorem}[Sufficient condition of optimality]
Let $f: \mathbb{R}^n \rightarrow \mathbb{R}$ be a function differentiable in an open
neighborhood of $\bm{x} \in \mathbb{R}^n$ and twice differentiable at $\bm{x}$. If $\nabla
f(\bm{x}) = 0$ and $\nabla^2 f(\bm{x})$ is positive definite, then $\bm{x}$ is a
\emph{local minimum} of $f$.
\label{th:sufficient-conditions}
\end{theorem}

\subsection{Constrained optimization}
In this section, the feasible set $D \subset \mathbb{R}^n$ is defined by
equality and inequality constraints (Definition \ref{def:constraint}):
\begin{equation}
D = \{\bm{x} \in \mathbb{R}^n ~|~ g(\bm{x}) \le 0 \land h(\bm{x}) = 0 \}
\end{equation}
where $g: \mathbb{R}^n \rightarrow \mathbb{R}^m$ and $h: \mathbb{R}^n \rightarrow \mathbb{R}^p$ are
continuous. A constrained optimization problem is defined in standard form:
\begin{equation}
\begin{aligned}
& \min_{\bm{x} \in \mathbb{R}^n} 	& f(\bm{x}) & \\
& s.t.	& g_j(\bm{x}) \le 0, & \quad j \in \{1, \ldots, m\} \\
&  		& h_j(\bm{x}) = 0, & \quad j \in \{1, \ldots, p\}
\end{aligned}
\tag{$\mathcal{P}$}
\label{eq:constrained-problem}
\end{equation}

\begin{definition}[Constraint, relation]
Let $\mathcal{V} = (x_1, \ldots, x_n)$ be a set of variables and $D$ their domain.
A \emph{constraint} $c$ is a logical expression:
\begin{equation}
c(x_1, \ldots, x_n) \diamond 0
\end{equation}
where $\diamond \in \{\le, \ge, =\}$. Reciprocally, $var(c)$ is the set of
variables that occur in the expression of $c$. The \emph{relation} $\rho_c$ of $c$ is the set
of solutions of $c$. 
\label{def:constraint}
\end{definition}

The first-order necessary condition of optimality in unconstrained optimization
(Theorem \ref{th:necessary-conditions}) does not apply in constrained optimization.
Example \ref{ex:constrained-example} shows that a global minimum that is not stationary
is located on the frontier of a constraint (the constraint is called \emph{active}, see
Definition \ref{def:active}).

\begin{example}
\label{ex:constrained-example}
Consider the following constrained problem:
\begin{equation}
\begin{aligned}
& \min_{x \in \mathbb{R}} & f(x) = x^2 & \\
& s.t.	& x \ge 1
\end{aligned}
\label{eq:constrained-example}
\end{equation}

\begin{figure}[htbp]
\centering
\def\svgwidth{0.5\columnwidth}
\small
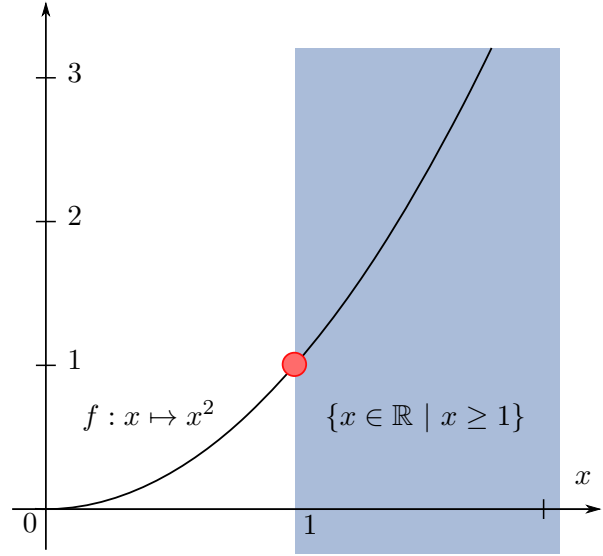
\caption{Global minimum of a constrained problem}
\label{fig:constrained-example}
\end{figure}

In Figure \ref{fig:constrained-example}, the colored rectangle represents the feasible
set $\{x \in \mathbb{R} ~|~ x \ge 1 \}$ of \ref{eq:constrained-example}. The global minimum
$x^* = 1$ is not stationary: $f'(x^*) = 2 \neq 0$.
\end{example}

\begin{definition}[Active/inactive constraint]
An inequality constraint $g \le 0$ is \emph{active} at $\bm{x} \in \mathbb{R}^n$
if $g(\bm{x}) = 0$, and \emph{inactive} if $g(\bm{x}) < 0$.
\label{def:active}
\end{definition}

Necessary conditions of optimality in constrained optimization (Theorem \ref{th:kkt})
rely upon the distinction between active and inactive inequality constraints.
An additional condition, constraint qualification (Definition \ref{def:qualification}), is required.

\begin{definition}[Constraint qualification]
\label{def:qualification}
The constraints of \ref{eq:constrained-problem} are \emph{qualified} at
$\bm{x} \in \mathbb{R}^n$ if the gradients of the active inequality constraints and
the gradients of the equality constraints are linearly independent at $\bm{x}$.
\end{definition}

\begin{theorem}[Karush-Kuhn-Tucker optimality conditions~\cite{Karush1939Minima}]
\label{th:kkt}
Suppose \\
that the functions $f$, $g_j~(j \in \{1, \ldots, m\})$ and $h_j~(j \in \{1, \ldots, p\})$ are continuously differentiable at a point $\bm{x}^* \in \mathbb{R}^n$. If $\bm{x}^*$ is a local minimum and the constraints are qualified at $\bm{x}^*$, then there exist $\lambda_j~(j \in \{1, \ldots, m\})$ and $\mu_j~(j \in \{1, \ldots, p\})$, called Lagrange multipliers, such that the following conditions are satisfied:

Stationarity
\begin{equation}
\nabla f(\bm{x}^*) + \sum_{j=1}^m \lambda_j \nabla g_j(\bm{x}^*) + \sum_{j=1}^p \mu_j \nabla h_j(\bm{x}^*) = 0
\end{equation}

Primal feasibility
\begin{equation}
\begin{aligned}
g_j(\bm{x}^*) & \le 0, \quad & j \in \{1, \ldots, m\} \\
h_j(\bm{x}^*) & = 0, \quad	 & j \in \{1, \ldots, p\}
\end{aligned}
\end{equation}

Dual feasibility
\begin{equation}
\lambda_j \ge 0, \quad j \in \{1, \ldots, m\}
\end{equation}

Complementarity
\begin{equation}
\lambda_j g_j(\bm{x}^*) = 0, \quad j \in \{1, \ldots, m\}
\end{equation}
\end{theorem}

Complementarity means that if $g_j$ is inactive at $\bm{x}^*$ ($g_j(\bm{x}^*) < 0$), then $\lambda_j = 0$. KKT conditions may be interpreted as optimality conditions for a problem in which active inequality constraints have been replaced by equality constraints, and inactive inequality constraints have been ignored (but must be satisfied).
The constraint qualification condition is necessary to guarantee the existence of Lagrange multipliers (Example \ref{ex:no-qualification}).

\begin{example}[Non-qualified constraint]
\label{ex:no-qualification}
Consider the following constrained problem:
\begin{equation}
\begin{aligned}
& \min_{x \in \mathbb{R}} & f(x) & = x & \\
& s.t.	& g(x) & = -x^3 \le 0
\end{aligned}
\label{eq:no-qualification}
\end{equation}
The optimal solution of \ref{eq:no-qualification} is $x^* = 0$. However, no $\lambda \ge 0$
satisfies the KKT condition:
\begin{equation}
f'(x^*) + \lambda g'(x^*) = 0
\end{equation}
The reason is that $g'(x^*) = 0$, that is $g$ is not qualified at $x^*$.
\end{example}

Note that any point that satisfies the KKT conditions is not a local minimum,
much like any stationary point of an unconstrained problem is not a local minimum.

\section{Optimization techniques}
\label{sec:optimization-techniques}

Numerous optimization techniques tackle optimization problems with particular structures.
In the following sections, we briefly introduce linear programming, convex optimization
and nonconvex optimization, as well as nonconvex optimization methods.

\subsection{Linear programming}
Solving a linear problem boils down to minimizing a linear function (Definition
\ref{def:linear_function}) over a convex polytope of $\mathbb{R}^n$ defined by
linear (in)equalities.

\begin{definition}[Linear function]
Let $\bm{x} = (x_1, \ldots, x_n) \in \mathbb{R}^n$.
A \emph{linear function} $f: \mathbb{R}^n \rightarrow \mathbb{R}$ of $\bm{x}$
may be written as a linear combination of the components of $\bm{x}$:
\begin{equation}
f(\bm{x}) = \sum_{i=1}^n a_i x_i = \bm{a}^T \bm{x}
\end{equation}
where $\bm{a} = (a_1, \ldots, a_n) \in \mathbb{R}^n$.
\label{def:linear_function}
\end{definition}

A linear problem (Definition \ref{def:linear-problem}), in which the objective and
constraints are linear, is generally written in canonical form.
The two most popular linear programming techniques are the simplex method, whose worst-case
complexity is exponential but is very efficient in practice, and interior point methods.

\begin{definition}[Linear problem]
Let $\bm{x} \in \mathbb{R}^n$, $\bm{c} \in \mathbb{R}^n$, $\bm{b} \in \mathbb{R}^m$ and $A \in \mathbb{R}^{m \times n}$ a matrix of size $m \times n$. The canonical form of the associated \emph{linear problem} is:
\begin{equation}
\begin{aligned}
& \min_{\bm{x} \in \mathbb{R}^n} & \bm{c}^T \bm{x} & \\
& \text{s.t.} & A\bm{x} & \le \bm{b}
\end{aligned}
\end{equation}
\label{def:linear-problem}
\end{definition}

\subsection{Convex problems}
A convex problem is a problem whose objective function is convex (Definition \ref{def:convex-function}) and whose feasible set is convex (Definition \ref{def:convex-set}).

\begin{definition}[Convex function]
\label{def:convex-function}
A function $f: I \rightarrow \mathbb{R}$, where $I$ is a real interval, is \emph{convex} if:
\begin{equation}
\forall (x, y) \in I^2, \quad \forall t \in [0, 1], \qquad f(tx + (1-t)y) \le tf(x) + (1-t)f(y)
\end{equation}
$f$ is \emph{strictly convex} if:
\begin{equation}
\forall (x, y) \in I^2, x \neq y, \quad \forall t \in ]0, 1[, \qquad f(tx + (1-t)y) < tf(x) + (1-t)f(y)
\end{equation}
\end{definition}

\begin{definition}[Convex set]
\label{def:convex-set}
A set $S$ is \emph{convex} if:
\begin{equation}
\forall (x, y) \in S^2, \quad \forall t \in [0, 1], \qquad tx + (1-t)y \in S
\end{equation}
\end{definition}
It can be easily shown that a local solution of a strictly convex optimization problem is the unique global minimum. Convex problems may be solved using interior point methods, cutting-plane methods, subgradient methods or bundle methods.

\subsection{Nonconvex problems}

Numerous real-world applications are nonconvex problems with multiple local minima. Therefore, convergence towards a local minimum does not guarantee global optimality.

\subsubsection{Local optimization methods}
Local optimization methods explore the neighborhood of an initial guess and offer a good tradeoff between quality of the solution and computational effort. They encompass two families of methods:
\begin{itemize}
\item mathematical programming methods generally exploit high-order information and compute a sequence of iterates in the search space. Gradient descent successively improves an initial solution by computing steps proportional to the negative of the gradient, in order to decrease the objective value. Newton-based methods consist in linearizing the optimality conditions (Sequential Quadratic Programming) or a perturbation thereof (Interior Point Method) ;

\item heuristic techniques are general methods that seek approximate solutions. Nelder-Mead method~\cite{NelderMead1965Simplex} maintains a polytope with $n+1$ vertices in an $n$-dimensional search space, which undergoes simple geometrical transformations, until hopefully reaching a local minimum. Pattern search~\cite{Hooke1961Direct} maintains $2n+1$ points in the search space in a similar fashion.
\end{itemize}

\subsubsection{Global optimization methods}
Global optimization methods seek the global minimum over the whole domain. Among them:
\begin{itemize}
\item metaheuristics are generic methods based on mechanisms such as local memory (taboo search~\cite{Glover1990Tabu}), greedy search (GRASP~\cite{Feo1989Probabilistic}) or random search (population-based algorithms, simulated annealing~\cite{Kirkpatrick1983Optimization}). Note that metaheuristics are equipped with mechanisms that help escape from local minima, but cannot guarantee the optimality of the solution ;

\item deterministic methods explore the search space in an exhaustive manner in order to identify the global minimum. They include branch and bound methods (see Chapter \ref{chap:ia}) and Lipschitz optimization.
\end{itemize}


\subsubsection{Reliable methods}
Deterministic global optimization methods, albeit exhaustive, do not provide a numerical guarantee of optimality of the solution within a given tolerance. A comparison of the main global optimization solvers~\cite{Neumaier2005Comparison} shows that most of them suffer from numerical approximations due to roundoff errors inherent to floating-point arithmetic.

Currently, only interval analysis provides rigorous bounds in numerical computations, even in the presence of roundoff errors. A variety of reliable interval-based techniques is presented in Chapter \ref{chap:ia}.

\subsection{Overview of the optimization techniques}

Table \ref{tab:taxonomy} is a summary of continuous optimization techniques and their characteristics. ``NC'' indicates nonconvex problems, ``$g \le 0$'' constrained problems and ``$n \uparrow$'' large-scale problems.

\begin{table}[htb]
	\centering
	\normalsize
	\caption{Continuous optimization techniques}
	\begin{tabular}{|l|ccc|ccc|}
	\hline
			& \multicolumn{3}{c|}{Problems} & \multicolumn{3}{c|}{Technique} \\
			& NC & $g \le 0$ & $n \uparrow$ & derivative-free & deterministic & global \\
	\hline
	Linear programming 		& & \checkmark & \checkmark & \checkmark & \checkmark & \checkmark \\
	Interior point methods	& & \checkmark & \checkmark & \checkmark & \checkmark & \checkmark \\
	Subgradient methods		& & \checkmark & \checkmark & \checkmark & \checkmark & \checkmark \\
	Quasi-Newton			& \checkmark & & \checkmark & & \checkmark & \\
	Simulated annealing 	& \checkmark & & \checkmark & \checkmark & & \checkmark \\
	Nelder-Mead 			& \checkmark & & \checkmark & \checkmark & \checkmark & \\
	Pop-based algorithms	& \checkmark & \checkmark & \checkmark & \checkmark & & \checkmark \\
	Lipschitz optimization	& \checkmark & \checkmark & \checkmark & \checkmark & \checkmark &
\checkmark \\
	Méthodes d'intervalles	& \checkmark & \checkmark & & \checkmark & \checkmark & \checkmark \\
	\hline
	\end{tabular}
	\label{tab:taxonomy}
\end{table}

\chapter{Evolutionary algorithms}
\label{chap:evolutionary-algorithms}

\minitoc

\begin{figure}[htb]
\centering
\def\svgwidth{0.7\columnwidth}
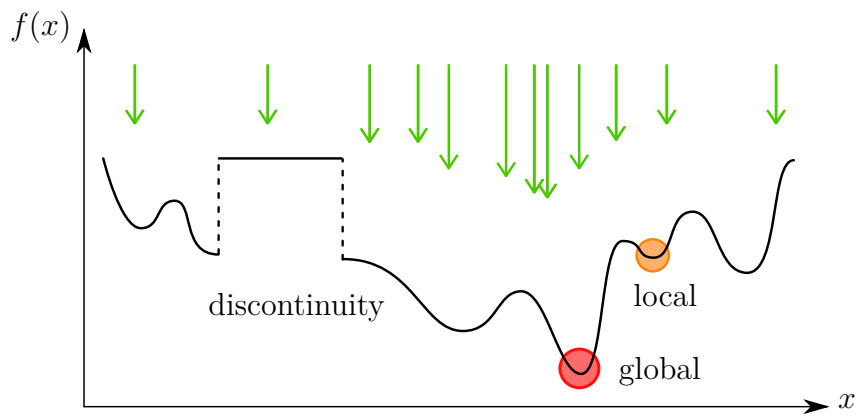
\caption{Sampling mechanism of metaheuristics}
\label{fig:metaheuristics}
\end{figure}

\Glspl{EA} are a subclass of metaheuristics that sample the search space in a stochastic manner, guided by the most promising individuals, and attempt to find a global minimum among the variety of local minima (Figure \ref{fig:metaheuristics}). They maintain a set of candidate solutions (individuals) in the search space, in order to obtain good solutions to an optimization problem. They do not require any regularity hypotheses of the objective function (continuity, differentiability), unlike mathematical programming techniques (gradient method, Newton-based methods). Only an evaluation procedure of the objective function is needed.
Section \ref{sec:evolutionary} presents the general framework of evolutionary algorithms. In particular, we explain genetic algorithms in Section \ref{sec:ga} and differential evolution algorithms in Section \ref{sec:de}.

\section{Evolutionary algorithms}
\label{sec:evolutionary}

\Gls{EA} mimics the evolution of a population in its environment (Algorithm \ref{alg:evolutionary}).
Variation operators bring diversity to the population in order to promote the exploration of the search space, while selection and replacement operators intensify the search in the vicinity of a solution.
Consequently, they are particularly suited to the optimization of difficult, multimodal, black box (whose analytical expression is not known), noisy or dynamic problems, for which other optimization methods fail at finding a satisfactory solution. Although \gls{EA} are equipped with mechanisms that help escape local minima, the optimality of the solution generally cannot be guaranteed.
 
\begin{algorithm}[htbp]
\caption{Evolutionary algorithms}
\label{alg:evolutionary}
\begin{algorithmic}
\State Initialize the initial population $P$
\State Evaluate the individuals
\Repeat
	\State Generate a new population $P'$
	\State Determine the new population by selection in $P \cup P'$
\Until{termination criterion met} \\
\Return best individual in the population
\end{algorithmic}
\end{algorithm}

\subsection{Constraint handling}
\label{sec:constraint-handling}
The ubiquity of constrained optimization problems has motivated the development of mechanisms that handle linear/nonlinear or equality/inequality constraints~\cite{Michalewicz1996Evolutionary,PriceStornLampinen2006, Talbi2009Metaheuristics}. Among the main techniques:

\begin{itemize}
\item rejection strategies do not exploit infeasible individuals, which are discarded. This strategy is relevant only if the feasible set is large relative to the search space ;
\item penalty methods (including logarithmic barrier terms and penalty terms that penalize constraint violations) consist in moving the constraints to the objective function as a weighted sum of penalties $p_i$ with weights $w_i$:
\begin{equation}
f_p(\bm{x}) = f(\bm{x}) + \sum_{j=1}^m w_j p_j(g_j(\bm{x}))
\label{eq:penalites}
\end{equation}
Penalty methods suffer from several drawbacks: if the weights $w_i$ are not properly adjusted, the weighted sum may be dominated by one of the penalties, or the metaheuristic may be trapped in an infeasible region if $f(\bm{x})$ is much smaller than the constraint violations $p_j(g_j(\bm{x}))_i$ ;
\item repair strategies are (usually greedy) problem-specific heuristics that generate a new feasible individual from an infeasible individual ;
\item decoding strategies establish a bijection between the set of representations of the individuals and the feasible set ;
\item direct methods consider an ordering of the individuals: feasible individuals get a better evaluation than infeasible individuals. Exploration is therefore guided towards feasible regions of the search space.
\end{itemize}

\subsection{Termination criteria}
There exist two categories of termination criteria:
\begin{itemize}
\item a \textit{static} criterion is generally based on available material resources (CPU time, number of iterations or number of evaluations) that are known a priori ;
\item a \textit{dynamic} criterion refers to the quality of the current solution (close enough to an optimum known a priori) or the end of convergence (number of consecutive iterations without improvements).
\end{itemize}

\section{Genetic algorithms}
\label{sec:ga}
Among the oldest \glspl{EA}, \glspl{GA} are inspired by the Darwinian theory of natural selection~\cite{Holland1975}:
genes that are the most adapted to the needs of a species in its environment are more likely to remain in a population over time. \Glspl{GA} establish a mapping between:
\begin{itemize}
\item the genotype (the genetic makeup carried by the chromosomes) of an individual and the components of a solution ;
\item the phenotype (the observable characteristics) of an individual and the objective value of a solution.
\end{itemize}

\Glspl{GA} maintain a population of $N$ individuals (a set of candidate solutions) that is partially replaced at each generation (iteration) through natural processes such as heredity, mutation and selection. A possible implementation is given in Algorithm \ref{alg:ga}.
\Glspl{GA} originally solved combinatorial optimization problems and encoded binary genes (0s and 1s). Since then, they have been extended to continuous optimization by adopting real-valued representations and have proven successful for a wide range of applications, such as bioinformatics, economy or chemistry.

\begin{algorithm}[htbp!]
\caption{Genetic algorithm}
\label{alg:ga}
\begin{algorithmic}
\State Initialize the initial population
\State Evaluate the individuals
\Repeat
	\State Select a pool of parents
	\State Combine the parents using crossover and mutation operators
	\State Evaluate the offspring
	\State Retain the best individuals
\Until{termination criterion is met} \\
\Return the best individual
\end{algorithmic}
\end{algorithm}

\subsection{Parent selection}
The theory of evolution states that individuals best adapted to their environments are more likely to survive, reproduce and pass on their gene pool to their offspring, whereas the maladapted die before reproducing. Two parent selection schemes are common in the \gls{GA} literature: roulette wheel selection~\cite{Goldberg1989Genetic} and stochastic remainder without replacement selection~\cite{Goldberg1989Genetic}.

\subsubsection{Roulette wheel selection}
The roulette wheel selection is similar to a roulette wheel in a casino. For a maximization problem, the probability $p_i$ of an individual $i \in \{1, \ldots, N\}$ to be selected as parent is proportional to its objective value $f_i \ge 0$:
\begin{equation}
p_i = \frac{f_i}{\sum_{j=1}^N f_j} \in [0, 1]
\end{equation}

Although individuals best adapted (with a higher objective value) are more likely to reproduce, it is not out of the questions that individuals that have low objective values succeed in reproducing, thus contributing to the next generation.
However, a selection bias may exist for small problems on account of the low number of selections.

\subsubsection{Stochastic remainder without replacement selection}
The stochastic remainder without replacement selection avoids the selection bias inherent to the roulette wheel selection. Each individual $i$ is replicated $\lfloor r_i \rfloor$ times, where:
\begin{equation}
r_i = \frac{N f_i}{\sum_{j=1}^N f_j}
\end{equation}
and $\lfloor \cdot \rfloor$ is the floor function. The roulette wheel selection is then performed on the set of all individuals, with objective values $r_i - \lfloor r_i \rfloor$.
This selection scheme generally produces better results when the population size is low.

\subsection{Crossover and mutation}
\begin{figure}[htbp]
\centering
\def\svgwidth{0.8\columnwidth}
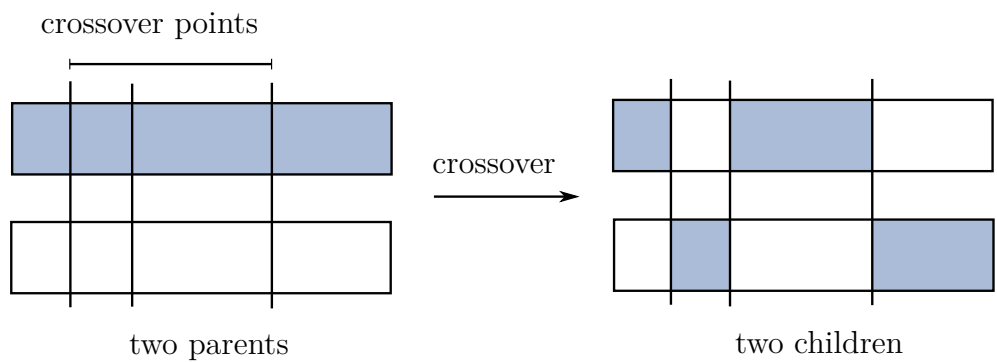
\caption{Crossover between two chromosomes}
\label{fig:ag-croisement}
\end{figure}

Crossover and mutation contribute to the diversification and intensification of the population. However their roles depend on the choice of implementation.
Crossover (or recombination) is the exchange of one or several portions of genetic material between two chromosomes (the parents) with a probability $p_c \in ~]0, 1[$.
Due to the crossover, the offspring (generally two chromosomes) have a different set of genes than their parents do (Figure \ref{fig:ag-croisement}).

\begin{figure}[htbp]
\centering
\def\svgwidth{0.8\columnwidth}
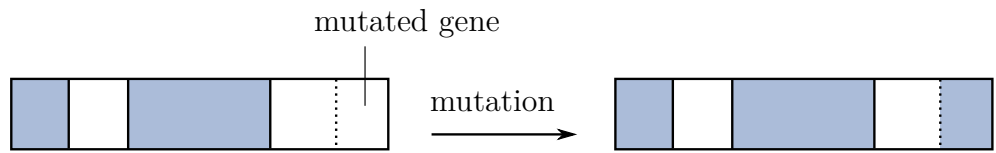
\caption{Mutation of a chromosome}
\label{fig:mutation}
\end{figure}

A gene can randomly mutate during crossover (Figure \ref{fig:mutation}). A chromosome subject to mutation thus possesses a genetic sequence that does not exclusively stem from its parents.
The mutation probability $p_m$ controls the randomness of the search: it allows an individual to avoid converging towards a local minimum by escaping from its vicinity.
Usually, $p_m$ is kept low in order to maintain the natural evolution of the population and to avoid turning the \gls{GA} into a mere random search.

\subsection{Population replacement}
\label{sec:elitism}
In order to determine which individuals will be kept in the next generation, parents and offspring can come face to face in a tournament: the individual with the better objective value wins the game and is kept in the population.
Elitism refers to the systematic conservation of the $k$ best individuals in the next generation. 
Usually, it is counterproductive to discard the individuals with the worse objective values, since they may carry genes that contribute to the elaboration of a satisfactory solution. Keeping "poor" individuals thus garantees the diversification of the population and reduces premature convergence towards local minima.
			
\section{Differential evolution}
\label{sec:de}
The \gls{DE} is among the simplest and most powerful \glspl{EA}~\cite{StornPrice1997}. Initially devised to handle unconstrained problems with continuous variables, \gls{DE} was extended to constrained and mixed problems.
Robust and relatively simple (it has few hyperparameters), it gained fame by solving difficult instances in aerodynamic design~\cite{Rogalsky2000Differential}, neural network training~\cite{Slowik2008Training}, multicriteria optimization, polynomial approximation and scheduling.
Contrary to \glspl{GA} whose new individuals are generated using crossover and mutation operators, \gls{DE} combines the components of existing individuals with a certain probability to build new individuals (Algorithm \ref{alg:de}).

\begin{algorithm}[htbp]
\caption{Differential evolution}
\label{alg:de}
\begin{algorithmic}[]
\Function{DifferentialEvolution}{$f$: objective function, $\mathit{NP}$: population size, $W$: scaling factor, $\mathit{CR}$: crossover rate, $\bm{D}$: domain}
\State $P \gets$ randomly initialized population in $\bm{D}$
\Repeat
	\State $P' \gets \varnothing$
	\Comment new population
	\For{$\bm{x} \in P$}
		\State $(\bm{u}, \bm{v}, \bm{w}) \gets$ \Call{ParentSelection}{$\bm{x}$, $P$}
		\State $\bm{y} \gets$ \Call{Crossover}{$\bm{x}$, $\bm{u}$, $\bm{v}$, $\bm{w}$, $W$, $\mathit{CR}$}
		\Comment generation of a new individual
		\If{$f(\bm{y}) < f(\bm{x})$}
			\State $P' \gets P' \cup \{\bm{y}\}$
			\Comment $\bm{y}$ replaces $\bm{x}$
		\Else 
			\State $P' \gets P' \cup \{\bm{x}\}$
			\Comment $\bm{x}$ is kept
		\EndIf
	\EndFor
	\State $P \gets P'$
	\Comment the temporary population replaces the current population
\Until{termination criteria met} \\
\Return best individual of $P$
\EndFunction
\end{algorithmic}
\end{algorithm}

\subsection{Quaternary crossover operator}
Let $\mathit{NP}$ be the population size, $W > 0$ be the scaling factor and $\mathit{CR} \in [0, 1]$ be the crossover ratio. At each generation, $\mathit{NP}$ new individuals are generated: for each individual $\bm{x} = (x_1, \ldots, x_n)$, three individuals $\bm{u} = (u_1, \ldots, u_n)$ (the base individual), $\bm{v} = (v_1, \ldots, v_n)$ and $\bm{w} = (w_1, \ldots, w_n)$, all different and different from $\bm{x}$, are picked from the population at random.
The components $y_i$ ($i \in \{1, \ldots, n\}$) of the new individual $\bm{y} = (y_1, \ldots, y_n)$ are computed as follows:
\begin{equation}
y_i =
\begin{cases}
u_i + W \times (v_i - w_i) & \text{if } i = R \text{ or } r_i < \mathit{CR} \\
x_i					& \text{otherwise}
\end{cases}
\end{equation}
where $r_i \in [0, 1]$ is a uniformly distributed random number and $R \in \{1, \ldots, n\}$ is a random index that guarantees that at least a component of $\bm{y}$ differs from that of $\bm{x}$. $\bm{y}$ replaces $\bm{x}$ in the population if it improves its objective value ; the selection operator is thus elitist (see Section \ref{sec:elitism}).

\begin{figure}[htbp]
\centering
\def\svgwidth{0.6\columnwidth}
\small
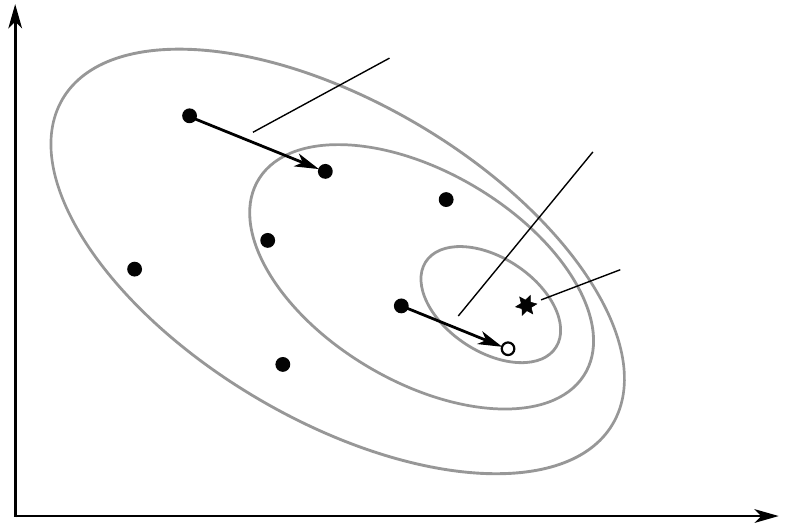
\caption{Quaternary crossover operator of the differential evolution algorithm}
\label{fig:de}
\end{figure}

Figure \ref{fig:de} illustrates the two-dimensional crossover between $\bm{x}$, $\bm{u}$ (the base individual), $\bm{v}$ and $\bm{w}$, and shows the contour map of the objective function. The difference $\bm{v}-\bm{w}$ determines the direction of displacement (an approximation of the direction opposite to the gradient) in which $\bm{u}$ is translated.

\subsection{Basic individual selection}
In a \gls{GA}, the probability that an individual is selected as parent is usually proportional to its objective value. In the \gls{DE}, all individuals are equally likely to become a base individual $\bm{u}$. Two variants were suggested by \cite{PriceStornLampinen2006} in order to guarantee that all individuals of the current population are selected once as base individuals at each generation:
\begin{enumerate}
\item the base vectors are picked from a random permutation of the population ;
\item the index of a base vector $\bm{u}$ is the sum modulo $\mathit{NP}$ of the index of $\bm{x}$ and an offset picked from $\{1, \ldots, \mathit{NP}\}$ at random.
\end{enumerate}

\subsection{Bound constraints}
The components of the newly generated individual $\bm{y}$ that are outside of the domain $\bm{D}$, thus violating the bound constraints, may be handled in two ways:
\begin{itemize}
\item the objective function is penalized: a term, either constant or depending on the number and the magnitude of the bound violations, is added to the objective value. This approach converges slowly when the newly generated individuals tend to often violate the bound constraints ;
\item a new component $y_i$ is generated in the domain $\bm{D}$: fixed at the bound of the domain, picked in the domain at random, or picked between $u_i$ and the bound of $\bm{D}$~\cite{PriceStornLampinen2006}:
\begin{equation}
y_i =
\begin{cases}
u_i + \omega(\overline{D_i} - u_i) & \text{if } y_i > \overline{D_i} \\
u_i + \omega(\underline{D_i} - u_i) & \text{if } y_i < \underline{D_i}
\end{cases}
\end{equation}
where $\omega \in [0, 1]$ is picked at random.
\end{itemize}

\subsection{Direct constraint handling}
Direct constraint handling consists in maintaining separate values for the objective function and the constraints for the comparison of individuals. The base vector $\bm{u}$ may be chosen according to the following rules:
\begin{itemize}
\item $\bm{u}$ is feasible and $\bm{x}$ is infeasible ;
\item $\bm{u}$ and $\bm{x}$ are feasible, and $f(\bm{u}) < f(\bm{x})$ ;
\item $\bm{u}$ and $\bm{x}$ are infeasible, and $\bm{u}$ does not violate the constraints more than $\bm{x}$.
\end{itemize}

\chapter{Interval analysis}
\label{chap:ia}

\minitoc

Interval analysis is a branch of numerical analysis dedicated to bounding roundoff errors. Interval methods are set-oriented enclosure methods that can compute rigorous lower and upper bounds of a function on a given interval, even in the presence of roundoff errors. They are therefore particularly suited to reliable global optimization.

Interval arithmetic is introduced in Section~\ref{sec:interval-computations}. The concept of interval extension is detailed in Section \ref{sec:interval-extensions}.
Interval branch and bound algorithms, dedicated to globally solving continuous optimization problems, are presented in Section~\ref{sec:ibb}. Section \ref{sec:ad} briefly covers automatic differentiation techniques.

\section{Interval computations}
\label{sec:interval-computations}

\Gls{FPA} is an approximate representation of real numbers on computers~\cite{Goldberg1991Every}. A floating-point number $x$ is represented by its sign, its significand (a fractional coefficient) and its exponent. The \glspl{FPU} embedded within computers handle floating-point numbers with fixed-size significands, which leads to approximation errors when real numbers are not exactly representable. For example, the constant $\pi$ rounded to three decimal digits is either 3.141 or 3.142, and its exact value lies somewhere in the interval $[3.141, 3.142]$.

The IEEE Standard for Floating-Point Arithmetic (IEEE 754) defines the floating-point representation of real numbers and the behavior of basic floating-point operations, and has been adopted by most \glspl{FPU} since the first standard in 1985. It defines five rounding modes: round to nearest (two modes), round toward 0, round toward $+\infty$ and round toward $-\infty$.
The IEEE 754 double-precision format provides 15 to 17 significant decimal digits precision. Despite this accurate precision, the accumulation of roundoff errors is the origin of erroneous results in numerically unstable problems (Example~\ref{ex:roundoff}).
\begin{example}[Accumulation of roundoff errors]
\label{ex:roundoff}
A catastrophic accumulation of roundoff errors was illustrated by~\cite{Rump1988Algorithms}. Consider the function:
\begin{equation}
f(x, y) = 333.75y^6 + x^2 (11x^2 y^2 - y^6 - 121y^4 - 2) + 5.5y^8 + \frac{x}{2y}
\end{equation}
When we evaluate $f(77617, 33096)$ in single and double precision, we obtain $1.172603$ and $1.1726039400531$, respectively. However, the exact value is $-\frac{54767}{66192} = -0.827396$.
\end{example}

\subsection{Rounding modes}
The seminal doctoral dissertation of~\cite{Moore1966} laid the foundations of interval computations: his idea was to enclose each step of a numerical computation within an interval that contains the true result. A real number $x$ exactly representable in \gls{FPA} is replaced by the degenerate interval $[x, x]$, while real numbers that cannot be represented exactly in \gls{FPA} (for instance 0.1) are rigorously bounded by an interval with floating-point bounds.

\Gls{IA} reliably extends real arithmetic to intervals. The implementation of elementary operations ($+$, $-$, $\times$, $/$, $\log$, $\exp$, etc.) requires correct (outward) rounding and exploits the rounding modes of the processor: the left (respectively right) bound of each intermediary step of the computation is rounded toward $-\infty$ (respectively $+\infty$). The true value is then numerically guaranteed to belong to the resulting interval.

\subsection{Interval arithmetic}

We adopt the following notations:
\begin{itemize}
\item $\mathbb{R}$ is the set of real numbers ;
\item $\mathbb{F}$ is the set of floating-point numbers ;
\item an interval $X = [\underline{X}, \overline{X}]$ with floating-point bounds defines the set:
\begin{equation}
X \eqdef \{x \in \mathbb{R} ~|~ \underline{X} \le x \le \overline{X} \}
\end{equation}
\item an interval is degenerate when $\underline{X} = \overline{X}$ ;
\item $\mathbb{I}$ is the set of intervals with floating-point bounds:
\begin{equation}
\mathbb{I} \eqdef \{[\underline{X}, \overline{X}] ~|~ (\underline{X}, \overline{X}) \in \mathbb{F}^2 \land \underline{X} \le \overline{X} \}
\end{equation}
\item given the set $D \subset \mathbb{R}$, $\mathbb{I}(D)$ denotes the set of intervals in $D$. The definition can be extended to the multivariate case ;
\item $int(X)$ is the interior of the non-degenerate interval $X$, that is the set:
\begin{equation}
int(X) \eqdef \{x \in \mathbb{R} ~|~ \underline{X} < x < \overline{X} \}
\end{equation}
\item $m(X) \eqdef \frac{1}{2}(\underline{X} + \overline{X})$ is the middle of the interval $X$ ;
\item $w(X) \eqdef \overline{X} - \underline{X}$ is the width of the interval $X$ ;
\item a box $\bm{X} = (X_1, \ldots, X_n)$ is a Cartesian product of intervals ;
\item $m(\bm{X}) \eqdef (m(X_1), \ldots, m(X_n))$ is the middle of the box $\bm{X}$ ;
\item $\displaystyle w(\bm{X}) = \max_{i=\{1, \ldots, n\}} w(X_i)$ is the width of the box $\bm{X}$ ;
\item $\square(X, Y)$ is the convex hull of $X$ and $Y$, that is the smallest interval of $\mathbb{I}$ that contains $X$ and $Y$.
\end{itemize}
In the rest of the document, upper-case letters denote interval quantities and bold letters denote vectors. An interval is thus written $X$, a box $\bm{X}$ and a real-valued vector $\bm{x}$.

The interval counterpart of a binary operator $\diamond \in \{+, -, \times, /\}$ provides the smallest interval that contains the range of the operator:
\begin{equation}
X \diamond Y = \square\{x \diamond y ~|~ x \in X \land y \in Y\}
\end{equation}
The left and right bounds of the results can be computed explicitly as functions of the bounds of the operands:
\begin{equation}
\begin{split}
[a, b] + [c, d] = & [a + c, b+d] \\
[a, b] - [c, d] = & [a-d, b+c] \label{eq:soustraction} \\
[a, b] \times [c, d] = & [\min(ac, ad, bc, bd), \max(ac, ad, bc, bd)] \\
\frac{1}{[a, b]} = & [\frac{1}{b}, \frac{1}{a}] \text{ if } 0 \notin [a, b] \\
\frac{[a, b]}{[c, d]} = & [a, b] \times \frac{1}{[c, d]} \text{ if } 0 \notin [c, d]
\end{split}
\end{equation}
Remember that interval computations must be evaluated using outward rounding (the left bound is rounded toward $-\infty$ and the right bound toward $+\infty$).

The interval counterparts of most unary operations exploit monotonicity:
\begin{equation}
\begin{aligned}
\exp([a, b]) 	& = [\exp(a), \exp(b)] 		& \\
\log([a, b]) 	& = [\log(a), \log(b)] 		& \quad (a > 0) \\
\sqrt{[a, b]} 	& = [\sqrt{a}, \sqrt{b}] 	& \quad (a \ge 0)
\end{aligned}
\end{equation}
Non-monotonic functions on a given interval (for instance, even powers or trigonometric functions on $\mathbb{R}$) are piecewise monotonic and must be studied more thoroughly.

\Gls{IA} exhibits weaker properties than real arithmetic: subtraction is not the inverse operation of addition in $\mathbb{I}$, just as division is not the inverse operation of multiplication (Example \ref{ex:non-inversibility}). Furthermore, the distributive property of multiplication over addition does not apply in $\mathbb{I}$ ; only a weaker subdistributivity property holds:
\begin{equation}
\forall (X, Y, Z) \in \mathbb{I}^3, \quad X(Y + Z) \subset XY + XZ
\end{equation}

\begin{example}[Non-inversibility of addition and multiplication]
\label{ex:non-inversibility}
Let $X = [2, 4]$ and $Y = [3, 5]$. Then:
\begin{equation}
\begin{aligned}
W & = X + Y = [2, 4] + [3, 5] = [5, 9] 		& \quad W - X & = [5, 9] - [2, 4] = [1, 7] \supsetneq Y \\
Z & = X \times Y = [2, 4] \times [3, 5] = [6, 20]	& \quad \frac{Z}{X} & = \frac{[6, 20]}{[2, 4]} = [1.5, 10] \supsetneq Y
\end{aligned}
\end{equation}
\end{example}

Extended interval arithmetic~\cite{Hanson1968Interval,Kahan1968More,Hansen1992} generalizes \gls{IA} to interval with infinite bounds, which allows the definition of the division by an interval that contains zero:
\begin{equation}
\frac{1}{[a, b]} =
\begin{cases}
	\varnothing 						& \text{ if } a = b = 0 \\
	[\frac{1}{b}, \frac{1}{a}] 				& \text{ if } 0 < a \text{ ou } b < 0 \\
	[\frac{1}{b}, +\infty] 					& \text{ if } a = 0 \\
	[-\infty, \frac{1}{a}] 					& \text{ if } b = 0 \\
	[-\infty, \frac{1}{a}] \cup [\frac{1}{b}, +\infty] 	& \text{ if } a < 0 < b \\
	\end{cases}
\end{equation}
or computations with infinite bounds:
\begin{equation}
\sqrt{[4, +\infty]} = [2, +\infty]
\end{equation}

Numerous software libraries implement \gls{IA}:
Profil/BIAS~\cite{Knuppel1994Profil} (a C++ library developed at the Technische Universität Hamburg),
Gaol~\cite{Gaol} (a C++ implementation of interval constraint programming operators),
Boost~\cite{Bronnimann2006Design} (using C++ templates),
MPFI~\cite{Revol2002} (a C and C++ multiprecision library),
Sun~\cite{SUN2001} (a Fortran 95 and C++ library) and Filib~\cite{Filib}.
We implemented a library for interval computations in the functional language OCaml~\cite {Alliot2012Implementing}. For the sake of performance, low-level routines (in C and assembly language) allow a fine control of the rounding.

\section{Interval extensions}
\label{sec:interval-extensions}

Interval extensions (Definition \ref{def:interval-extension}) build upon the conservativity of interval computations to compute rigorous enclosures of factorable functions (Definition \ref{def:factorable-function}).

\begin{definition}[Factorable function]
A factorable function can be recursively written as a finite composition of elementary operations (operators, functions or variables).
\label{def:factorable-function}
\end{definition}

\begin{definition}[Range]
Let $f: D \subset \mathbb{R}^n \rightarrow \mathbb{R}$ and $\bm{X} \in \mathbb{I}(D)$. $f(\bm{X})$ is the \emph{range} of $f$ on $\bm{X}$:
\begin{equation}
f(\bm{X}) \eqdef \{f(\bm{x}) ~|~ \bm{x} \in \bm{X} \}
\end{equation}
\end{definition}

\begin{definition}[Interval extension]
\label{def:interval-extension}
Let $f: D \subset \mathbb{R}^n \rightarrow \mathbb{R}$ and $F: \mathbb{I}(D) \rightarrow \mathbb{I}$. $F$ is an \emph{interval extension} of $f$ if:
\begin{equation}
\begin{aligned}
\forall \bm{X} \subset \mathbb{I}(D)	&, \quad f(\bm{X}) \subset F(\bm{X}) & \text{
(conservativity)} \\
\forall (\bm{X}, \bm{Y}) \in \mathbb{I}(D)^2 	&, \quad \bm{X} \subset \bm{Y}
\implies F(\bm{X}) \subset F(\bm{Y}) & \text{ (inclusion isotonicity)}
\end{aligned}
\end{equation}
\end{definition}

The interval extensions of a real-valued function are not unique: interval extensions of various orders of convergence (Definition \ref{def:convergence-order}) can generally be constructed, that is the overestimation tends to zero at different speeds when the size of the interval tends to zero.
\begin{definition}[Order of convergence]
Let $f: D \subset \mathbb{R}^n \rightarrow \mathbb{R}$ and $F$ be an interval extension of $f$. $F$ has an \emph{order of convergence} $\alpha > 0$ if:
\begin{equation}
w(F(\bm{X})) - w(f(\bm{X})) = O(w(\bm{X})^\alpha), \quad \forall \bm{X} \in \mathbb{I}(D)
\label{eq:convergence-order}
\end{equation}
$w(F(\bm{X})) - w(f(\bm{X}))$ represents the overestimation error of the range of $f$ on $\bm{X}$.
\label{def:convergence-order}
\end{definition}

The most straightforward interval extension of a factorable function is the natural interval extension   (Definition \ref{def:natural-extension}). It has a linear order of convergence, that is the overestimation tends to zero linearly as the size of the interval tends to zero.

\begin{definition}[Natural interval extension]
Let $f: D \subset \mathbb{R}^n \rightarrow \mathbb{R}$. The \emph{natural interval extension} $F_N: \mathbb{I}(D) \rightarrow \mathbb{I}$ of $f$ is obtained by replacing each variable with its domain and each elementary operation with its interval counterpart in the expression of $f$.
\label{def:natural-extension}
\end{definition}

\subsection{Dependency}
\label{sec:dependency}

Evaluating expressions with equivalent syntaxes in real arithmetic can produce interval enclosures with various accuracies. The main reason why \gls{IA} may overestimate -- sometimes dramatically -- the range of a function is known as the dependency problem: the multiple occurrences of a variable are decorrelated and handled as distinct variables.

For instance, the interval $X = [a, b]$ subtracted to itself produces $Y = X - X = [a - b, b - a]$.
When $a < b$, $Y$ is not the degenerate interval $[0, 0]$ (although the exact result 0 belongs to $Y$).
Here, the overestimation error (the width of $Y$) is $2(b-a)$, that is twice that of $X$. The range computed by \gls{IA} is in fact $\{x_1 - x_2 ~|~ x_1 \in X, x_2 \in X\}$. \Gls{IA} ignores the dependency between $x_1$ and $x_2$ and constructs an enclosure of $f(X)$ with too many degrees of freedom.

In his fundamental theorem of \gls{IA} (Theorem \ref{th:moore-fundamental}), \cite{Moore1966} proved that, under particular circumstances, \gls{IA} does not overestimate the range of $f$.
\begin{theorem}[Fundamental theorem of interval arithmetic~\cite{Moore1966}]
If a function $f$ is continuous over a box $\bm{X}$ and all the variables occur at most once in the expression of $f$, the natural interval extension of $f$ provides the optimal image in exact real arithmetic:
\begin{equation}
F_N(\bm{X}) = f(\bm{X})
\end{equation}
\label{th:moore-fundamental}
\end{theorem}

The following example highlights the importance of the continuity assumption in Theorem \ref{th:moore-fundamental}.
\begin{example}[Discontinuous function on an interval]
Let $f(x) = (\frac{1}{x})^2$ and $X = [-1, 1]$. The image of $X$ under $f$ is $f(X) = [1, +\infty]$. Although $x$ occurs only once in the expression of $x$, \gls{IA} evaluates $\frac{1}{X} = [-\infty, +\infty]$, then $(\frac{1}{X})^2 = [-\infty, +\infty]^2 = [0, +\infty] \supset f(X)$. The overestimation is a consequence of the discontinuity of $f$ over $X$.
\end{example}

A straightforward method to reduce the overestimation of \gls{IA} is to rewrite (when it is possible) the expression of the function (Example \ref{ex:rewriting}).

\begin{example}[Rewriting of an expression]
\label{ex:rewriting}
Let $f(x) = \frac{x-y}{x+y}$, $x \in X = [6, 8]$ and $y \in Y = [2, 4]$. The range of $f$ is $f(X, Y) = [\frac{1}{5}, \frac{3}{5}]$. The natural interval extension of $f$ provides the following enclosure:
\begin{equation}
F_N(X, Y) = \frac{X-Y}{X+Y} = \frac{[6, 8]-[2, 4]}{[6, 8]+[2, 4]} = \frac{[2, 6]}{[8, 12]} = [\frac{1}{6}, \frac{3}{4}] \supset f(X, Y)
\end{equation}
Let us rewrite $f$ so that $x$ and $y$ have single occurrences in the expression:
\begin{equation}
\hat{f}(x, y) = \frac{x-y}{x+y} = \frac{x+y-2y}{x+y} = 1-\frac{2y}{x+y} = 1-\frac{2}{1+\frac{x}{y}}
\end{equation}
\Gls{IA} produces the following range:
\begin{equation}
\begin{aligned}
\hat{F}_N(X, Y) 	& = 1-\frac{2}{1+\frac{X}{Y}} = 1-\frac{2}{1+\frac{[6, 8]}{[2, 4]}} = [\frac{1}{5}, \frac{3}{5}] = g(X, Y)
\end{aligned}
\end{equation}
Rewriting the expression of $f$ so that the variables have single occurrences in the expression gets rid of the dependency effect ; \gls{IA} produces the exact range.
\end{example}

The objective function and the constraints of an optimization problem are usually complex and the variable have multiple occurrences. It is not always possible to rewrite the expressions, which often leads to crude enclosures, even when the intervals are small. Solving difficult optimization problems using interval methods is thus an arduous task. However, numerous authors attempted to tackle the overestimation problem.

The first approach is second-order interval extensions. Recall that for an interval extension with an order of convergence $\alpha$, the overestimation of the range of $f$ on $X$ decreases with $w(X)^\alpha$. The interval enclosures therefore tend to be tighter for small intervals. Second-order interval extensions are presented in Section \ref{eq:2nd-order}.

The second approach is monotonicity. If a function is monotonic with respect to some of its variables on a box, the enclosure of the range can be reduced to punctual evaluations at the bounds of the box. Indubitably the most powerful method, detecting monotonicity eliminates the dependency problem related to these variables. The method is described in Section \ref{sec:monotonicity-based-extension}.

The third approach to reduce the overestimation of \gls{IA} is to reduce the width of the intervals. Interval extensions that are inclusion isotonic produce tighter enclosures on small intervals. Partitioning $X$ into $\{X_i\}_{i=1 \ldots P}$ may improve (sometimes strictly) the enclosure of the range of the function. This concept is exploited by branch and bound methods that alternate between interval evaluation and partitioning. The framework is described in Section \ref{sec:ibb}.

\subsection{Second-order extensions}
\label{eq:2nd-order}


Let $f: D \subset \mathbb{R} \rightarrow \mathbb{R}$, $X \in \mathbb{I}(D)$ and $c \in X$ (for example $c = m(X)$). Taylor's theorem states that, when $f$ is $m-1$ times differentiable at $c$ and $m$ times differentiable on the open interval, we have for $x \in X$:
\begin{equation}
\begin{aligned}
f(x) 	& = \sum_{k=0}^{m-1} \frac{f^{(k)}(c)}{k!}(x-c)^k + \frac{f^{(m)}(\xi)}{m!} (x-c)^m \\
		& \in \sum_{k=0}^{m-1} \frac{f^{(k)}(c)}{k!}(x-c)^k + \frac{F^{(m)}(X)}{m!} (x-c)^m
\end{aligned}
\label{eq:taylor}
\end{equation}
where $\xi$ is a real number between $c$ and $x$, and $F^{(m)}$ is an interval extension of $f^{(m)}$. 
This inclusion defines the Taylor interval extension (or Taylor form) of order $m$. The most commonly used Taylor forms are the linear ($m = 1$) and quadratic ($m = 2$) interval extensions. In the rest of the document, we focus on the linear form:
\begin{equation}
F_{mv}(X, c) \eqdef f(c) + F'(X)(X-c)
\label{eq:mean-value-form}
\end{equation}
also known as the \textit{mean value extension} (in reference to the \textit{mean value theorem}~\cite{Jeffreys1999Methods}). When $c = m(X)$, $F_{mv}$ is inclusion isotonic if $F'$ is inclusion isotonic~\cite{Caprani1980Mean}, and has a quadratic convergence if $F'$ is a Lipschitz continuous (Definition \ref{th:lipschitz-continuous}) function~\cite{Krawczyk1982Centered}.

\begin{definition}[Lipschitz continuous function~\cite{Moore1966}]
\label{th:lipschitz-continuous}
Let $f: D \subset \mathbb{R}^n \rightarrow \mathbb{R}$ and $F: \mathbb{I}(D) \rightarrow
\mathbb{I}$ be an interval extension of $f$. $F$ is a \emph{Lipschitz continuous function} if it exists $K > 0$ such that:
\begin{equation}
\forall \bm{X} \in \mathbb{I}(D), \quad w(F(\bm{X})) \le Kw(\bm{X})
\end{equation}
\end{definition}

\begin{remark}[Roundoff errors]
Since $f(c)$ may be subject to roundoff errors, it must be replaced with $F(c) \eqdef F([c, c])$.
\end{remark}

\begin{figure}[htbp]
\centering
\def\svgwidth{0.7\columnwidth}
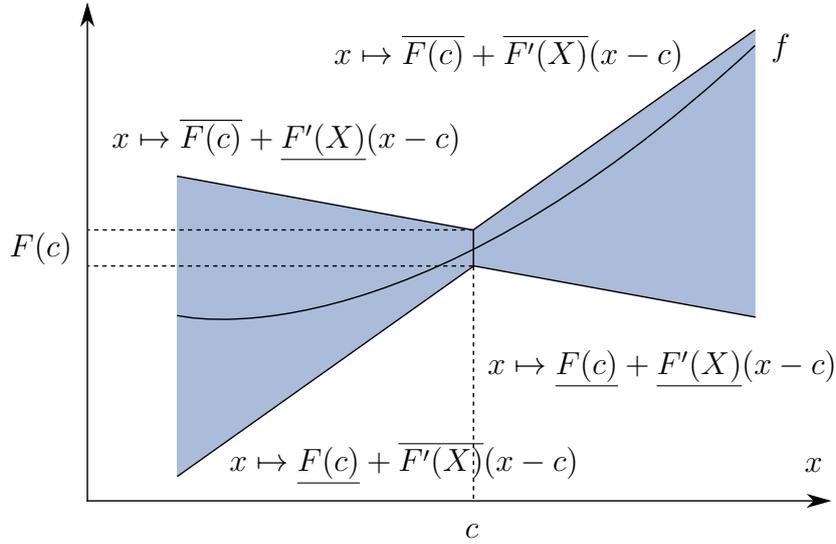
\caption{Mean value extension}
\label{fig:mean-value-extension}
\end{figure}

Figure \ref{fig:mean-value-extension} illustrates the mean value extension of a univariate function. The function is enclosed by a cone that contains the tangents at the point $(c, f(c))$ with all possible slopes in $F'(X)$.

The Taylor extension can be easily extended to a multivariate function $f: D \subset \mathbb{R}^n
\rightarrow \mathbb{R}$:
\begin{equation}
f(\bm{X}) \subset F(\bm{c}) + \sum_{i=1}^n \frac{\partial F}{\partial x_i}(\bm{X}) \cdot
(X_i-x_i)
\label{eq:taylor-multi}
\end{equation}
where $\bm{X} = (X_1, \ldots, X_n) \in \mathbb{I}(D)$, $\bm{c} = (c_1, \ldots, c_n) \in \bm{X}$ and $\frac{\partial F}{\partial x_i}$ is an interval extension of the $i$th partial derivative of $f$.
The $n$ partial derivatives can be evaluated simultaneously using automatic differentiation (Section \ref{sec:ad}).
\cite{Hansen1968Solving} proposed a recursive variant in which the Taylor series is computed variable after variable.

\begin{figure}[htbp]
\centering
\def\svgwidth{0.7\columnwidth}
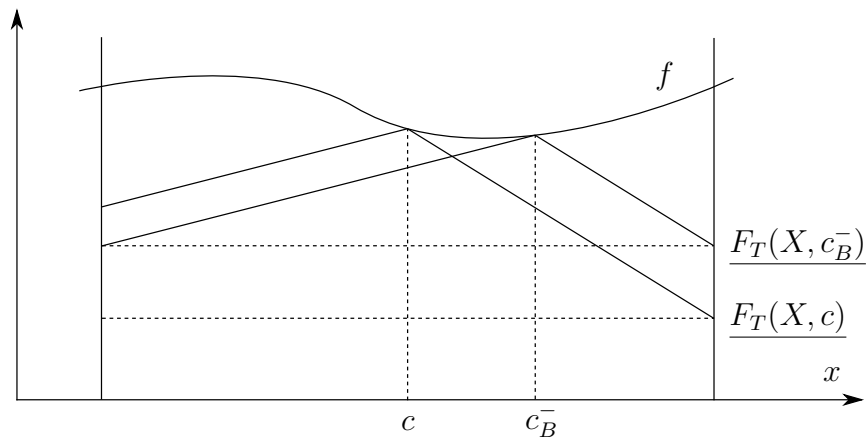
\caption[Mean value extensions with an arbitrary center and with the optimal Baumann center]{Mean value extensions $F_{mv}(X, c)$ with an arbitrary center $c$ and $F_{mv}(X, c_B^-)$ with the optimal Baumann center $c_B^-$}
\label{fig:baumann}
\end{figure}

\cite{Baumann1988} gave the analytical expression of the optimal center $c_B^-$ (\resp{} $c_B^+$) that maximizes the lower bound (\resp{} minimizes the upper bound) of the mean value extension (Figure \ref{fig:baumann}):
\begin{equation}
\label{eq:centre-baumann}
c_B^- \eqdef
\begin{cases}
\underline{X}								& \text{if } 0 \le F'(X) \\
\overline{X}								& \text{if } 0 \ge F'(X) \\
\frac{U\underline{X} - L\overline{X}}{U-L}	& \text{otherwise}
\end{cases}
\quad \quad
c_B^+ \eqdef
\begin{cases}
\overline{X}								& \text{if } 0 \le F'(X) \\
\underline{X}								& \text{if } 0 \ge F'(X) \\
\frac{U\overline{X} - L\underline{X}}{U-L}	& \text{otherwise}
\end{cases}
\end{equation}
where $[L, U] \eqdef F'(X)$.

\subsection{Monotonicity-based extension}
\label{sec:monotonicity-based-extension}

The local monotonicity of a function with respect to some of its variables (Definition \ref{def:variable-monotone}) eliminates the dependency effect related to these variables and computes tighter enclosures than interval extensions with linear or quadratic convergence.

\begin{definition}[Local monotonicity]
Let $f: D \subset \mathbb{R}^n \rightarrow \mathbb{R}$ be a continuous function, $F$ be an interval extension of $f$ and $\bm{X} \in \mathbb{I}(D)$. $f$ is \emph{locally increasing} (\resp{} \emph{decreasing}) with respect to $x_i$ on $\bm{X}$ if $\frac{\partial F}{\partial x_i}(\bm{X})$ is nonnegative (\resp{} nonpositive).
\label{def:variable-monotone}
\end{definition}

\begin{example}[Local monotonicity]
Let $f(x) = x + 2\cos(x)$ and $X = [3, 6]$. The natural interval extension of $f'(x) = 1 - 2\sin(x)$ on $X$ produces (rounded to 3 significant digits) $F'_N(X) = 1 - 2\sin([3, 6]) = 1 - 2[-1, 0.142] = [0.717, 3] \ge 0$. $f$ is therefore locally increasing with respect to $x$ on $X$.
\end{example}

The monotonicity-based extension $F_M$ (Definition \ref{def:monotonicity-extension}) computes an enclosure of the range of $f$ tighter than that of the natural extension $F_N$ when $f$ is detected locally monotonic on $\bm{X}$ with respect to variables that have multiple occurrences in its expression:
\begin{equation}
f(\bm{X}) \subset F_M(\bm{X}) \subset F_N(\bm{X})
\end{equation}

\begin{definition}[Monotonicity-based extension]
\label{def:monotonicity-extension}
Let $f: D \subset \mathbb{R}^n \rightarrow \mathbb{R}$ and $\bm{X} = (X_1, \ldots, X_n)$. $\bm{X}^- = (X_1^-, \ldots, X_n^-)$ and $\bm{X}^+ = (X_1^+, \ldots, X_n^+)$ denote the boxes defined for $i \in \{1, \ldots, n\}$ by:
\begin{equation}
X_i^- \eqdef
\begin{cases}
\underline{X_i} & \text{ if } f \text{ increasing wrt } x_i \\
\overline{X_i} 	& \text{ if } f \text{ decreasing wrt } x_i \\
X_i				& \text{ otherwise}
\end{cases}
\quad \quad
X_i^+ \eqdef
\begin{cases}
\overline{X_i} 	& \text{ if } f \text{ increasing wrt } x_i \\
\underline{X_i} & \text{ if } f \text{ decreasing wrt } x_i \\
X_i				& \text{ otherwise}
\end{cases}
\end{equation}

The \emph{monotonicity-based extension} $F_M$ of $f$ is defined by:
\begin{equation}
F_M(\bm{X}) \eqdef [\underline{F_N(\bm{X}^-)}, \overline{F_N(\bm{X}^+)}]
\end{equation}
\end{definition}

Since the variables with respect to which $f$ is monotonic are replaced by one of their bounds in $\bm{X}^-$ and $\bm{X}^+$, the dependency problem related to these variables disappears in $F_M$ (see Example~\ref{ex:monotonicity-extension}). Consequently, if $f$ is monotonic with respect to all the variables with multiple occurrences, $F_M$ produces an exact enclosure (that may however be overconservative on account of rounding) of the image of $f$ on $\bm{X}$.

\begin{example}[Monotonicity-based extension]
\label{ex:monotonicity-extension}
Let $f(x) = -x_1^2 + x_1 x_2 + x_2 x_3 - 3 x_3$ and $\bm{X} = [6, 8] \times [2, 4] \times [7, 15]$. 
The natural extension of $f$ produces $F_N(\bm{X}) = -[6, 8]^2 + [6, 8] \times [2, 4] + [2, 4] \times 
[7, 15] - 3[7, 15] = [-83, 35]$.
The natural interval extensions of the partial derivatives of $f$ with respect to $x$ are:
\begin{align}
\frac{\partial F}{\partial x_1}(\bm{X}) & = -2X_1 + X_2 = [-14, -8] \le 0 \\
\frac{\partial F}{\partial x_2}(\bm{X}) & = X_1 + X_3 = [13, 23] \ge 0 \\
\frac{\partial F}{\partial x_3}(\bm{X}) & = X_2 - 3 = [-1, 1]
\end{align}
$f$ is decreasing with respect to $x_1$ and increasing with respect to $x_2$ on $\bm{X}$. However, $f$ is not monotonic with respect to $x_3$ on $\bm{X}$. The monotonicity-based extension $F_M$ of $f$ produces:
\begin{equation}
F_M(\bm{X}) = [\underline{F(8, 2, [7, 15])}, \overline{F(6, 4, [7, 15])}] = [-79, 27] \subset
F_N(\bm{X})
\end{equation}
\end{example}

\cite{Araya2010} improved the monotonicity-based extension by computing the partial derivatives of $f$ independently on $\bm{X}^-$ and $\bm{X}^+$: since the variables with respect to which $f$ is monotonic are replaced with one of their bounds, the enclosures of the partial derivatives are tighter.
This new interval extension $F_{Mrec}$, called \textit{recursive monotonicity-based extension}, computes an enclosure that is always at least as good as the monotonicity-based extension $F_M$:
\begin{equation}
f(\bm{X}) \subset F_{Mrec}(\bm{X}) \subset F_M(\bm{X}) \subset F_N(\bm{X})
\end{equation}

\subsection{Affine arithmetic}
\Gls{AA} is an alternative to \gls{IA} to automatically compute an enclosure of a function on a box. Each quantity of a computation is represented by a linear combinations of symbols ; the linear dependencies between the variables are memorized (for example, $X + 2Y - 3X$ is not subject to the dependency problem) and the nonlinear operations are linearized by introducing error terms.

\Gls{AA} has proven very efficient for global optimization~\cite{Messine1997,Ninin2010Reliable}, in particular because the enclosure techniques have a quadratic convergence. However, it remains tricky to implement and was not used in our work.

\section{Interval branch and bound methods}
\label{sec:ibb}

Historically devised to bound rounding errors in numerical computations, set-oriented interval computations have increasingly attracted attention over the last years. New applications include global optimization, robust optimization, constraint satisfaction, root finding and numerical integration.

The generic branch and bound framework for global optimization is presented in Section \ref{sec:bb}, and extended to interval computations in Section \ref{sec:ibb}. Sections \ref{sec:ibb-heuristics} and \ref{sec:ibb-accelerating} detail various heuristics and acceleration techniques.

\subsection{Branch and bound methods}
\label{sec:bb}

\Gls{BB}~\cite{Lawler1966Branch} is a generic framework for solving combinatorial constrained optimization problems: 
\begin{equation}
\begin{aligned}
(\mathcal{P}) \quad & \min_{\bm{x} \in D} 		& f(\bm{x}) & \\
		    & s.t.	& g_j(\bm{x}) \le 0 	& \quad (j \in \{1, \ldots, m\}) \\
		    & 		& h_j(\bm{x}) = 0 	& \quad (j \in \{1, \ldots, p\}) \\
\end{aligned}
\end{equation}
where $D$ is a discrete set.

\begin{remark}
Generally, an equality constraint $h_j = 0$ is considered satisfied if the relaxed constraint $-\varepsilon_= \le h_j \le \varepsilon_=$ (with $\varepsilon_=$ arbitrarily small) is satisfied.
\end{remark}

It is not always possible to enumerate all the elements of $D$, either because there exists no simple algorithm for doing so, or because the cardinality of $D$ precludes it. A \gls{BB} algorithm partitions the search space into subspaces and builds a search tree in which the leaves are punctual solutions. The subspaces that cannot contain the optimal solution are discarded and the algorithm is applied recursively on the remaining subspaces. The best solution of the subspaces is then the solution to the original problem.

Although the worst-case complexity is exponential, discarding large subspaces often avoids the systematic enumeration of all the elements of $D$. \gls{BB} algorithms usually converge in finite time (Theorem \ref{th:convergence-bb}), albeit not necessarily reasonable.

\begin{theorem}[Convergence of a \gls{BB} algorithm]
\label{th:convergence-bb}
If the partitioning of the search space is a finite process, the \gls{BB} algorithm terminates in finite time.
\end{theorem}

\Gls{BB} algorithms proceed by bounding the ranges of the objective function and constraints, and maintaining the best known upper bound $\tilde{f}$ of the global minimum $f^*$ (for example, the evaluation of a feasible point). On each subspace $S_k$, lower bounds $f^-(S_k)$, $g_j^-(S_k)$, $h_j^-(S_k)$ and upper bounds $f^+(S_k)$, $g_j^+(S_k)$, $h_j^+(S_k)$ of the ranges of the objective function $f$ and the constraints $g_j$, $h_j$ are computed. $S_k$ can then be classified according to Definition \ref{def:bb-subspaces}.

\begin{definition}[Classification of subspaces]
\label{def:bb-subspaces}
Let $S_k$ be a subspace of the search space.
\begin{itemize}
\item $S_k$ is \emph{feasible} if all the elements of $S_k$ satisfy $g_j$ and $h_j$:
\begin{equation}
\forall j \in \{1, \ldots, m\}, \quad g_j^+(S_k) \le 0
\end{equation}
\begin{equation}
\forall j \in \{1, \ldots, p\}, \quad -\varepsilon_= \le h_i^-(S_k) \land h_i^+(S_k) \le
\varepsilon_=
\end{equation}
$\tilde{f}$ can be updated with $f^+(S_k)$ ;
\item $S_k$ is \emph{infeasible} if a constraint $g_j$ is violated by all the elements of $S_k$:
\begin{equation}
\exists j \in \{1, \ldots, m\}, \quad 0 < g_j^-(S_k)
\end{equation}
or a constraint $h_j$ is violated by all the elements of $S_k$:
\begin{equation}
\exists j \in \{1, \ldots, p\}, \quad \varepsilon_= < h_j^-(S_k) \lor h_j^+(S_k) < -\varepsilon_=
\end{equation}
\item $S_k$ is \emph{suboptimal} if $S_k$ is feasible and if all the elements of $S_k$ have an objective value worse than $\tilde{f}$:
\begin{equation}
\tilde{f} < f^-(S_k)
\end{equation}
\item $S_k$ is \emph{undecided} if it is neither feasible, nor infeasible, nor suboptimal: $S_k$ may contain a global minimizer.
\end{itemize}
\end{definition}

During the exploration of the search space, subspaces that are infeasible or suboptimal are discarded with the guarantee that they cannot contain an optimal solution. The remaining subspaces are stored in an appropriate data structure ; its order determines the nature of the search space exploration. Several heuristics are mentioned in Section \ref{sec:ibb-heuristics}.
The feasible or undecided subspaces are recursively processed until they are discarded. The optimal solution to $(\mathcal{P})$ is obtained once all the subspaces have been processed.

An extension of the \gls{BB} algorithm, known as spatial \gls{BB}~\cite{Ryoo1995Global}, aims at solving continuous or mixed (involving discrete and continuous variables) problems. A lower bound of the problem on a subspace is computed by solving a convex relaxation. Since floating-point points have a finite representation on a computer, the spatial \gls{BB} algorithm converges in finite time (Theorem \ref{th:convergence-bb}). The conservative properties of \gls{IA} allow to automatically generate rigorous enclosures of the range of a factorable function on a box ; interval-based \gls{BB} algorithms are described in Section \ref{sec:ibb}.

\subsection{Interval branch and bound methods}

The first \gls{IBB} algorithms were introduced by \cite{Moore1976Computing} and \cite{Skelboe1974Computation}: the search space is partitioned into disjoint boxes on which interval extensions of the objective function and the constraints are evaluated. If a box cannot be discarded, it is partitioned into subboxes that are inserted into a priority queue $\mathcal{Q}$ and processed at a later stage. For a given tolerance $\varepsilon$, the algorithm returns a punctual solution $\tilde{\bm{x}}$ with an objective value $\overline{F(\tilde{\bm{x}})} = \tilde{f}$ such that $\tilde{f} - f^* < \varepsilon$, even in the presence of roundoff errors. Subboxes $\bm{X}$ are discarded when no point of $\bm{X}$ can improve $\tilde{f}$ by at least $\varepsilon$, that is when $\tilde{f} - \varepsilon < \underline{F(\bm{X})}$.

A generic \gls{IBB} algorithm is described in Algorithm \ref{alg:ibb}. Advanced implementations embed accelerating techniques (Section \ref{sec:ibb-accelerating}) that eliminate or reduce boxes without losing solutions.

\begin{algorithm}[htbp]
\caption{Interval branch and bound methods}
\label{alg:ibb}
\begin{algorithmic}
\Function{IBB}{$\bm{X}_0$: initial box, $F$: objective function, $\mathcal{C}$: set of constraints}
\State $(\tilde{\bm{x}}, \tilde{f}) \gets (\varnothing, +\infty)$
\Comment Best known solution
\State $\mathcal{Q} \gets \{ \bm{X}_0 \}$
\Comment Priority queue
\While{$\mathcal{Q} \neq \varnothing$}
	\State Extract a box $\bm{X}$ from $\mathcal{Q}$
	\Comment Section \ref{sec:search-heuristics}
	\State Apply accelerating techniques
	\Comment Section \ref{sec:ibb-accelerating}
	\State Evaluate constraints in $\mathcal{C}$
	\If{$\bm{X}$ cannot be discarded}
		\State Update the best known solution $(\tilde{\bm{x}}, \tilde{f})$
		\State Partition $\bm{X}$ into $\{\bm{X}_1, \ldots, \bm{X}_k\}$
		\Comment Section \ref{sec:search-heuristics}
		\State Insert $\{\bm{X}_1, \ldots, \bm{X}_k\}$ into $\mathcal{Q}$
	\EndIf
\EndWhile
\State \Return $(\tilde{\bm{x}}, \tilde{f})$
\EndFunction
\end{algorithmic}
\end{algorithm}

\begin{remark}
Since a box is a vector of closed intervals, bisecting a box produces two subboxes that share a face.
\end{remark}

\subsection{Heuristics}
\label{sec:ibb-heuristics}

Heuristics play a crucial role in the speed of convergence of \gls{IBB} algorithms, and are generally problem-dependent. They can be divided into two categories: exploration heuristics and partitioning heuristics.

\subsubsection{Exploration of the search space}
\label{sec:search-heuristics}

The order in which the subboxes are inserted into and extracted from the priority queue $\mathcal{Q}$ determines how the search space is explored. The most widely used heuristics are:
\begin{itemize}
\item "best-first search": the subbox $\bm{X}$ with the lowest $\underline{F_N(\bm{X})}$ is extracted from $\mathcal{Q}$. This strategy promotes the exploration of the most promising subspaces ;
\item "largest first search": the box with the maximal size is extracted from $\mathcal{Q}$. Also known as "breadth first search", this strategy explores the oldest subspaces first ;
\item "depth-first search": this strategy explores the most recent subspaces first.
\end{itemize}

\subsubsection{Box partitioning}
\label{sec:box-partitioning}

Historically, two partitioning heuristics have stood out:
\begin{itemize}
\item the variable with the largest domain is partitioned ;
\item the variables are partitioned one after the other in a round-robin fashion.
\end{itemize}
In recent years, the Smear heuristic~\cite{Csendes1997Subdivision} has proven a competitive alternative to the aforementioned strategies:
the variable $x_i$ with the largest Taylor extension $\frac{\partial F}{\partial x_i}(\bm{X}) \cdot (X_i - x_i)$ over the box $\bm{X}$ is partitioned.

\subsection{Accelerating techniques}
\label{sec:ibb-accelerating}

\subsubsection{Upper bounding}

The interval upper bound $\overline{F(\bm{X})}$ computed over a feasible box $\bm{X}$ is guaranteed to be an upper bound of the global minimum $f^*$, but is coarse in that it is not necessarily the image of any point of $\bm{X}$. \cite{Ichida1979Interval} suggested that the midpoint $m(\bm{X})$ of a box $\bm{X}$ be systematically evaluated to provide a possibly better upper bound of $f^*$ (Example \ref{ex:ibb}). If $m(\bm{X})$ is a feasible point, the best known upper bound $\tilde{f}$ of $f^*$ may be replaced with $\overline{F(m(\bm{X}))}$. Some implementations perform local search in order to update the best known solution.

\begin{remark}
If $\bm{X}$ is a feasible box, it is not necessary to test the feasibility of $m(\bm{X})$. If $\bm{X}$ however may contain feasible and infeasible points, the set of constraints is evaluated on $m(\bm{X})$ using \gls{IA}.
\end{remark}

\begin{example}[Interval branch and bound]
\label{ex:ibb}
Let $f(x) = x^2\cos(x) + x$. Let us search for the global minimum of $f$ over the interval $X = [-5, 3]$ with a precision $\varepsilon$ (Figure \ref{fig:ibb}). The optimal range of $f$ over $X$ is approximately $f(X) = [-15.311, 2.092]$.

\clearpage

\begin{figure}[htbp]
\centering
\def\svgwidth{0.8\columnwidth}
\scriptsize
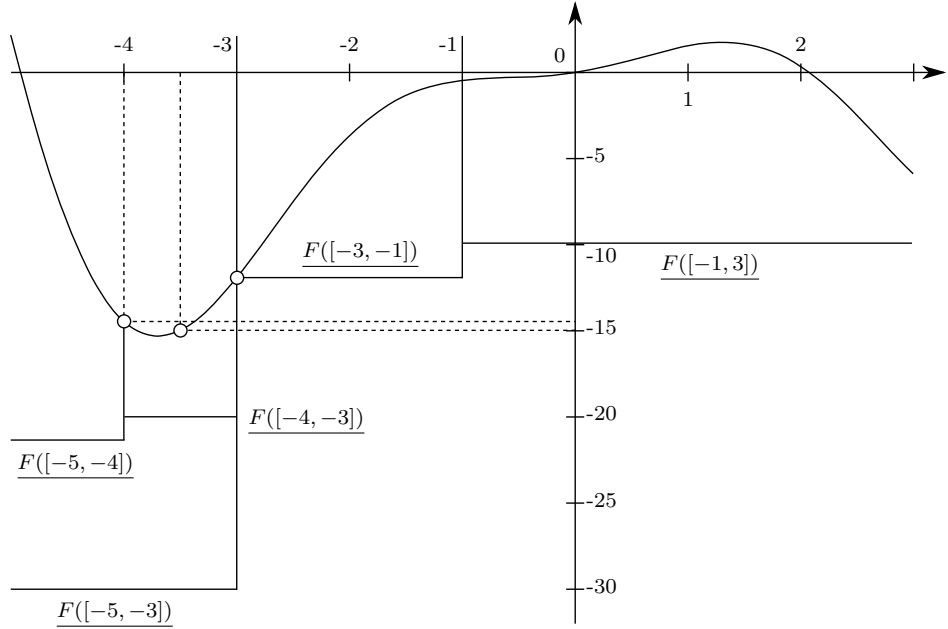
\caption{Global optimization of $f(x) = x^2\cos(x) + x$ over $X = [-5, 3]$ using an interval branch and bound algorithm}
\label{fig:ibb}
\end{figure}

\begin{enumerate}
\item the best known upper bound $\tilde{f}$ of the global minimum is initially set to $+\infty$ ;
\item $X$ is partitioned into $X_1 = [-5, -1]$ and $X_2 = [-1, 3]$ ;
\item process $X_1 = [-5, -1]$: an interval lower bound of $f$ is $\underline{F_N(X_1)} = -30$. $\tilde{f}$ is updated using the midpoint of $X_1$: $\tilde{f} \gets \overline{F(-3)} = -11.91$. $X_1$ is partitioned into $X_3 = [-5, -3]$ and $X_4 = [-3, -1]$ ;
\item process $X_2 = [-1, 3]$: an interval lower bound of $f$ is $\underline{F_N(X_2)} = -9.91 \ge \tilde{f}$. $X_2$ cannot contain a global minimizer and can be safely discarded ;
\item process $X_4 = [-3, -1]$: an interval lower bound of $f$ is $\underline{F_N(X_4)} = -11.91 \ge \tilde{f}$. $X_4$ cannot contain a global minimizer and can be safely discarded ;
\item process $X_3 = [-5, -3]$: an interval lower bound of $f$ is $\underline{F_N(X_3)} = -30$. $\tilde{f}$ is updated using the midpoint of $X_3$: $\tilde{f} \gets \overline{F(-4)} = -14.458$. $X_3$ is partitioned into $X_5 = [-5, -4]$ and $X_6 = [-4, -3]$ ;
\item process $X_6 = [-4, -3]$: an interval lower bound of $f$ is $\underline{F_N(X_6)} = -20$. $\tilde{f}$ is updated using the midpoint of $X_6$: $\tilde{f} \gets \overline{F(-4)} = -14.972$ ;
\item process $X_5 = [-5, -4]$: an interval lower bound of $f$ is $\underline{F_N(X_5)} = -21.341$.
\end{enumerate}

At this stage, the interval lower bounds of $f$ over $X_5$ and $X_6$ do not provide enough information to discard these subboxes. Nonetheless, a global minimizer is guaranteed to lie within $X_5 \cup X_6 = [-5, -3]$. Both subboxes can be subsequently partitioned and recursively processed, in order to update $\tilde{f}$ and refine the computation of interval lower bounds. The algorithm terminates when all subboxes $X_k$ satisfy the cut condition $\tilde{f} - \varepsilon \le \underline{F_N(X_k)}$.
\end{example}

\subsubsection{Lower bounding}

\cite{Kearfott1993Cluster,Du1994Cluster} have shown that the convergence time of an \gls{IBB} algorithm is strongly influenced by the behavior of the function in the neighborhood of the minima, and by the order of convergence of the interval extensions used in the algorithm. Second-order extensions (Section \ref{eq:2nd-order}), monotonicity-based extensions (Section \ref{sec:monotonicity-based-extension}) and \gls{AA} are powerful tools to generate high-quality lower bounds.

\subsubsection{Box reduction}

Advanced techniques that reduce boxes without losing solutions, called filtering operators or contractors, are presented in Chapter \ref{chap:contractors}.

\section{Automatic differentiation}
\label{sec:ad}

\Glspl{CAS} that implement symbolic differentiation generate formulae for the partial derivatives of a function based on its analytical expression. However, the size of the expressions grows rapidly with that of the function, and \glspl{CAS} usually struggle to provide a unique expression of the partial derivatives when the function involves conditional statements, such as if-then-else constructs.

\Gls{AD} is a set of techniques (divided into direct or reverse modes) for evaluating \textit{at a particular point} the partial derivatives of a function described by a finite composition of elementary functions (for example, a computer program). The expression of the function can be modeled as a syntax tree whose nodes are elementary functions and whose leaves are variables or constants. \gls{AD} exploits the chain rule (Theorem \ref{th:chain-rule} and Example \ref{ex:chain-rule}) at a particular point.

\begin{theorem}[Chain rule]
\label{th:chain-rule}
Let $g: \mathbb{R}^m \rightarrow \mathbb{R}^n$ be differentiable at the point $x \in \mathbb{R}^m$, and let $f: \mathbb{R}^n \rightarrow \mathbb{R}$ be differentiable at the point $g(x) \in \mathbb{R}^n$. The Jacobian of $f \circ g$ at the point $x$ can be written as the product of the Jacobians of $f$ and $g$ at the points $g(x)$ and $x$, respectively:
\begin{equation}
J_{f \circ g}(x) = J_f(g(x)) J_g(x)
\end{equation}
\end{theorem}

\begin{example}[Chain rule]
\label{ex:chain-rule}
Let $g: \mathbb{R}^2 \rightarrow \mathbb{R}^2$ be differentiable at the point $x \in \mathbb{R}^2$, and let $f: \mathbb{R}^2 \rightarrow \mathbb{R}$ be differentiable at the point $g(x) \in \mathbb{R}^2$. The derivatives of $f \circ g$ at the point $x$ are given by:
\begin{equation}
\label{eq:chain-rule}
\begin{split}
\frac{d(f \circ g)}{dx_1}(x) & = \frac{\partial f}{\partial g_1}(g(x)) \frac{\partial g_1}{\partial x_1}(x) + \frac{\partial f}{\partial g_2}(g(x)) \frac{\partial g_2}{\partial x_1}(x) \\
\frac{d(f \circ g)}{dx_2}(x) & = \frac{\partial f}{\partial g_1}(g(x)) \frac{\partial g_1}{\partial x_2}(x) + \frac{\partial f}{\partial g_2}(g(x)) \frac{\partial g_2}{\partial x_2}(x)
\end{split}
\end{equation}
\end{example}

\subsection{Direct mode}
The direct mode (1965-1970) evaluates Equation \ref{eq:chain-rule} from right to left: the directional derivatives
\begin{equation}
\nabla(f \circ g)(x) \cdot d = \lim_{h \rightarrow 0} \frac {(f \circ g)(x + h d) - (f \circ g)(x)}{h}
\end{equation}
along direction $d$ are computed in a single bottom-up evaluation phase, starting from the leaves of the syntax tree (the variables) and simultaneously evaluating the intermediary computations of the function and its gradient. The process can be extended to higher-order derivatives. The computation of the $n$ partial derivatives of $f$ using the direct mode has a complexity of $O(ne)$, where $e$ is the number of elementary operations that occur in the expression the $f$.

\subsection{Adjoint mode}
\label{sec:ad-adjoint}
The adjoint (or reverse) mode (1976-1980) evaluates Equation \ref{eq:chain-rule} from left to right, thus exploiting the fact that $\frac{\partial f}{\partial g_1}(g(x))$ et $\frac{\partial f}{\partial g_2}(g(x))$ occur in both $\frac{d(f \circ g)}{dx_1}(x)$ and $\frac{d(f \circ g)}{dx_2}(x)$. A bottom-up evaluation phase computes the values of the intermediary operations (the nodes), then a top-down propagation phase computes the partial derivative of each node with respect to its children (Example \ref{ex:adjoint-ad}). The total derivative of $f$ with respect to a variable $x_i$ is then gathered by summing the partial derivatives of $f$ with respect to all the occurrences of $x_i$ (the leaves of the syntax tree). The computation of the $n$ partial derivatives of $f$ has a complexity of $O(e)$, independent of $n$.

\begin{example}[Automatic differentiation in adjoint mode]
\label{ex:adjoint-ad}
Let $f$ be a two-dimensional function:
\begin{equation}
f(x, y) = \cos^2(x) - xy
\end{equation}

$f$ is a finite composition of elementary functions. The intermediary results (the nodes of the syntax tree) $t_i$ are computed during the evaluation phase:
\begin{equation}
\begin{aligned}
t_1 = x 		& \quad t_4 = t_3^2 \\
t_2 = y 		& \quad t_5 = t_1 t_2 \\
t_3 = \cos(t_1) & \quad t_6 = t_4 - t_5 := f(x, y)
\end{aligned}
\end{equation}

The gradient of $f$ at the point $(x, y)$ is computed in adjoint mode as follows:
\begin{equation}
\begin{aligned}
\frac{\partial t_6}{\partial t_6} & = 1 \qquad \frac{\partial t_6}{\partial t_5} = -1 \qquad \frac{\partial t_6}{\partial t_4} = 1 \\
\frac{\partial t_5}{\partial t_1} & = t_2 \qquad \frac{\partial t_5}{\partial t_2} = t_1 \\
\frac{\partial t_4}{\partial t_3} & = 2 t_3 \\
\frac{\partial t_3}{\partial t_1} & = -\sin(t_1) \\
\frac{\partial f}{\partial x}(x, y) \eqdef \frac{\partial t_6}{\partial t_1} & = \frac{\partial t_6}{\partial t_5} \frac{\partial t_5}{\partial t_1} + \frac{\partial t_6}{\partial t_4} \frac{\partial t_4}{\partial t_3} \frac{\partial t_3}{\partial t_1} = -t_2 - 2t_3 \sin(t_1) = -y - 2\cos(x)\sin(x) \\
\frac{\partial f}{\partial y}(x, y) \eqdef \frac{\partial t_6}{\partial t_2} & = \frac{\partial t_6}{\partial t_5} \frac{\partial t_5}{\partial t_2} = -t_1 = -x
\end{aligned}
\end{equation}

The computation of the gradient of $f$ at the point $(1, 3)$ (in \gls{FPA}) and on the box $([0.9, 1], [2.9, 3.1])$ (using \gls{IA}) is detailed in Table \ref{tab:adjoint-ad}. Notice that, since $(1, 3)$ belongs to $([0.9, 1], [2.9, 3.1])$, all intermediary results using interval analysis enclose the results computed using \gls{FPA}.
\begin{table}[h!]
	\centering
	\caption{Automatic differentiation in adjoint mode}
	\begin{tabular}{|c|cc|}
	\hline
	Nodes & at the point $(1, 3)$ & on the box $[0.9, 1] \times [2.9, 3.1]$ \\
	\hline
	$t_1$ & 1 & [0.9, 1] \\
	$t_2$ & 3 & [2.9, 3.1] \\
	$t_3$ & 0.54 & [0.54, 0.622] \\
	$t_4$ & 0.292 & [0.291, 0.39] \\
	$t_5$ & 3 & [2.61, 3.1] \\
	$t_6$ & -2.71 & [-2.81, -2.22] \\
	\hline
	$\frac{\partial f}{\partial x}(x, y)$ & -3.91 & [-4.15, -3.74] \\
	$\frac{\partial f}{\partial y}(x, y)$ & -1 & [-1, -0.9] \\
	\hline
	\end{tabular}
	\label{tab:adjoint-ad}
\end{table}
\end{example}

\section{Conclusion}
Interval analysis is the method of choice for solving global optimization problems in the presence of roundoff errors. Due to the inherent dependency problem however, interval overestimation strongly hinders the efficiency of optimization algorithms. Advanced techniques such as higher-order interval extensions, affine arithmetic and constraint propagation are nowadays systematically embedded within state-of-the-art solvers. In the next chapter, we investigate filtering algorithms (or contractors) that reduce boxes by removing inconsistent values.

\chapter{Contractors}
\label{chap:contractors}

\minitoc

\Gls{IBB} methods are nowadays endowed with contraction procedures in order to reduce the domains of the variables with respect to individual constraints (local consistency) or all the constraints simultaneously (global consistency). The resulting framework, called \gls{IBC}, alternates between contraction (and evaluation) phases and branching phases. This chapter compares various local and global consistencies and the associated contractors. Contraction procedures stem from the interval analysis and \gls{ICP} communities:
\begin{itemize}
\item the interval lower bound of the objective function of a constrained problem is generally coarse, since it does not take the feasible set into account. Convexification-based contractors generate an outer linearization of the objective function and the constraints, then compute a lower bound of the objective function over the polyhedral feasible region and/or reduce the ranges of the variables ;
\item interval constraint programming, inspired by constraint programming~\cite{Mackworth1977Consistency}, encompasses a set of techniques that reduce the ranges of the variables by enforcing consistencies in a fixed-point algorithm.
\end{itemize}

\section{Partial consistency operators}
\label{sec:consistencies}

Any numerical optimization problem can be reformulated as a \gls{NCSP} (Definition \ref{def:ncsp}) in which a dynamic constraint $f \le \tilde{f}$ on the best known upper bound of the global minimum is maintained.

\begin{definition}[Constraint satisfaction problem]
A \gls{CSP} is defined as a triple $P = (\mathcal{V}, \mathcal{C}, D)$, where $\mathcal{V}$ is the set of variables, $\mathcal{C}$ is the set of constraints and $D$ is the set of domains. Solving $P$ boils down to finding one (or all) instantiation of $\mathcal{V}$ in $D$ that satisfies $\mathcal{C}$. A \gls{NCSP} is a \gls{CSP} for which $D$ is a subset of $\mathbb{R}^n$.
\label{def:ncsp}
\end{definition}

Solving an \gls{NCSP} resorts to the notion of partial consistency, a local property related to the consistency of variables and constraints ; inconsistent values -- values of the domain that are not solutions of a constraint -- are discarded by contraction (or filtering) operators (Definition \ref{def:contractor}).

\begin{figure}[htbp]
\centering
\def\svgwidth{0.4\columnwidth}
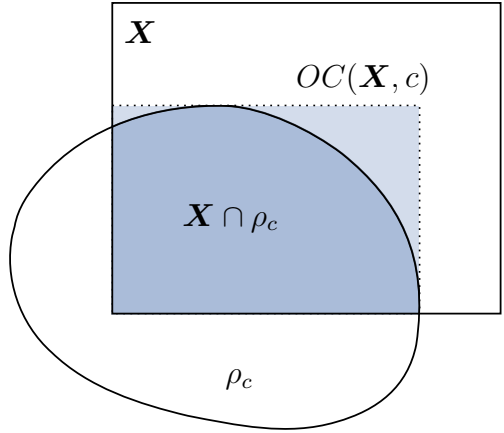
\caption{Contraction of a box $\bm{X}$ with respect to a constraint $c$}
\label{fig:contractor}
\end{figure}

\begin{definition}[Contractor~\cite{ChabertJaulin2009}]
\label{def:contractor}
Let $\bm{X} \in \mathbb{I}^n$ be a box, $c$ a constraint and $\rho_c$ the relation defined by $c$. An outer \emph{contractor} associated with $c$ is a mapping $OC$ that satisfies the correction property (Figure \ref{fig:contractor}), that is it reduces $\bm{X}$ by discarding values that are inconsistent with respect to $c$:
\begin{equation}
\bm{X} \cap \rho_c \subset OC(\bm{X}, c) \subseteq \bm{X}
\label{eq:contractor-correction}
\end{equation}
A contractor that enforces a partial consistency $\phi$ is called $\phi$-consistency operator and is denoted by $OC_{\phi}$.
\end{definition}

Contractors may exhibit the following properties:
\begin{align}
\text{Monotonic: } 	& \bm{X} \subset \bm{Y} \Rightarrow OC(\bm{X}, c) \subset OC(\bm{Y}, c) \\
\text{Convergent: } & \bm{X}_k \underset{k \to +\infty}{\longrightarrow} \bm{x} \Rightarrow OC(\bm{X}_k, c) \underset{k \to +\infty}{\longrightarrow} \bm{x} \cap \rho_c \\
\text{Idempotent: } & \forall \bm{X} \in \mathbb{I}^n, OC(OC(\bm{X}, c), c) = OC(\bm{X}, c) \\
\text{Minimal: } & \forall \bm{X} \in \mathbb{I}^n, OC(\bm{X}, c) = \square(\bm{X} \cap \rho_c) \\
\text{Thin: } & \forall \bm{x} \in \mathbb{R}^n, OC(\bm{x}, c) = \bm{x} \cap \rho_c
\end{align}

Contractor programming~\cite{ChabertJaulin2009} boils down to defining operations (intersection, union, composition, repetition) in order to build more complex contractors (Definition \ref{def:prog-contractors}).
\begin{definition}[Intersection and composition of contractors]
\label{def:prog-contractors}
Let $\bm{X}$ be a box, $c_1$ and $c_2$ two constraints and $OC_1$ and $OC_2$ two contractors associated with $c_1$ and $c_2$, respectively. The intersection and composition of $OC_1$ and $OC_2$ are defined by:
\begin{align}
\text{Intersection: } & OC_1(\bm{X}, c_1) \cap OC_2(\bm{X}, c_2) \\
\text{Composition: } & OC_2(OC_1(\bm{X}, c_1), c_2)
\end{align}
\label{def:algebra-contractors}
\end{definition}
The composition of contractors has a higher filtering power than a mere intersection and is in practice always implemented.

Section \ref{sec:fixed-point} details the fixed-point algorithm, a propagation loop that contracts a box with respect to a system of constraints. Section \ref{sec:local-consistencies} introduces two types of partial consistency based on the arc consistency for discrete domains: the 2B (or hull) consistency and the box consistency. The partial consistencies 3B and CID, stronger than 2B and box consistencies, are detailed in Section \ref{sec:strong-consistencies}. Global consistency algorithms, based on linearization techniques, are presented in Section \ref{sec:global-consistencies}.

\section{Fixed-point algorithm}
\label{sec:fixed-point}
A fixed-point algorithm is an idempotent propagation loop that contracts a box $\bm{X}$ with respect to a system of constraints $\mathcal{C}$ (Algorithm \ref{algo:fixed-point}). A contractor $OC$, called revising procedure, handles individual constraints of $\mathcal{C}$.
The list $\mathcal{Q}$ originally contains the constraints of $\mathcal{C}$. At each iteration, a constraint $c_i$ is extracted from $\mathcal{Q}$. If $\bm{X}$ is contracted with respect to $c_i$, all the constraints that involve the contracted variables of $\bm{X}$ are "waken up" and inserted into $\mathcal{Q}$. Otherwise, the next constraint is handled by the revising procedure.

\begin{algorithm}[htbp]
\caption{Fixed-point algorithm}
\label{algo:fixed-point}
\begin{algorithmic}
\Function{FixedPoint}{$\bm{X}, \mathcal{C}, OC$}
\State $\mathcal{Q} \gets \mathcal{C}$
\Repeat
	\State Extract a contract $c_i$ from $\mathcal{Q}$
	\State $\bm{X}' \gets OC(\bm{X}, c_i)$
	\Comment contraction
	\If{$\bm{X}' \neq \bm{X}$}
		\State $\mathcal{Q} \gets \mathcal{Q} \cup \{c_j ~|~ c_j \in \mathcal{C} \land \exists x_k \in var(c_j), X_k' \neq X_k \}$
		\State $\bm{X} \gets \bm{X}'$
	\EndIf
	\State $\mathcal{Q} \gets \mathcal{Q} \setminus \{c_i\}$
\Until{$\mathcal{Q} = \varnothing$}
\EndFunction
\end{algorithmic}
\end{algorithm}

\newpage

\begin{remark}
To avoid slow convergence, constraints are woken up when variables are \emph{sufficiently} contracted, for example when the ratio between the size of the contracted interval and the size of the initial interval is lower than a threshold $\eta \in [0, 1]$:
\begin{equation}
\exists x_k \in var(c_j), \quad w(X_k') < \eta w(X_k) 
\end{equation}
\end{remark}

\section{Local consistencies}
\label{sec:local-consistencies}

\subsection{2B consistency}
The 2B or hull consistency enforces the property of arc consistency for each bound of the variables that occur in a constraint (Definition \ref{def:2b}). Geometrically speaking, each face of a box that is 2B-consistent with respect to a system of constraints $\mathcal{C}$ intersects all constraints of $\mathcal{C}$.
\begin{definition}[2B consistency]
Let $c$ be an $n$-ary constraint and $\bm{X}$ a box. $\bm{X}$ is 2B-consistent with respect to $c$ if:
\begin{equation}
\forall i \in \{1, \ldots, n\}, \quad X_i = \square{\{x_i \in X_i ~|~ \exists x_1 \in X_1, \ldots,
\exists x_n \in X_n, c(x_1, \ldots, x_i, \ldots, x_n) \}}
\end{equation}
$\bm{X}$ is 2B-consistent with respect to a system of constraints $\mathcal{C}$ if it is 2B-consistent with respect to each constraint of $\mathcal{C}$.
\label{def:2b}
\end{definition}

\subsection{Forward-backward propagation}
\label{sec:hc4revise}

The evaluation-propagation algorithm~\cite{Messine1997}, also knows as HC4Revise~\cite{Benhamou1999} and FBBT~\cite{Belotti2009Branching}, draws its inspiration from the work of~\cite{Cleary1987Logical} on relational \gls{IA}. It computes an approximation of the 2B consistency for an explicit constraint by carrying out a double traversal of its syntax tree, in order to contract the domain of each occurrence of the variables (Example~\ref{ex:hc4}). HC4Revise is the revising procedure of the fixed-point algorithm HC4.

The bottom-up evaluation phase computes the value of each intermediary node using \gls{IA}. The result at the root of the tree -- the evaluation of the constraint -- is then intersected with the right-hand side of the constraint. The top-down propagation phase propagates the result by exploiting inverse (or projection) operations at each node~\cite{Goualard2008Interval}. Whenever an inconsistency is detected, the box is discarded, since it cannot satisfy the constraint. Otherwise, the variables (the leaves of the tree) may be contracted.

\begin{example}[HC4Revise algorithm]
\label{ex:hc4}
Let:
\begin{itemize}
\item $2x = z - y^2$ be an equality constraint ;
\item $X = [0, 20]$, $Y = [-10, 10]$ and $Z = [0, 16]$ be the domains of $x$, $y$ and $z$.
\end{itemize}
The constraint can be written as a finite composition of elementary operations:
\begin{equation}
\begin{aligned}
n_1 & \eqdef 2x 	& n_3 & \eqdef z-n_2 \\
n_2 & \eqdef y^2 	& n_4 & = n_3 \\
\end{aligned}
\end{equation}

The bottom-up evaluation phase (Figure \ref{fig:hc4-forward}) evaluates the nodes of the syntax tree:
\begin{equation}
\begin{aligned}
N_1 & = 2X = 2 \times [0, 20] = [0, 40] \\
N_2 & = Y^2 = [-10, 10]^2 = [0, 100] \\
N_3 & = Z - N_2 = [0, 16] - [0, 100] = [-100, 16]
\end{aligned}
\end{equation}

\begin{figure}[htbp]
\centering
\def\svgwidth{\columnwidth}
\begin{tikzpicture}[
	>=stealth,
	level 1/.style={sibling distance=8cm},
	level 2/.style={sibling distance=4cm},
	elem/.style={circle, draw=black, thin, minimum size = 0.5cm},
	var/.style={circle, draw=black, fill=gray!30},
	cst/.style={draw=none},
	line/.style = {->, draw=black},
	eval/.style = {->, dotted, red, very thick}
	]

\node [elem] (Root) {$=$} 
	child {node [elem] (n4) {$\times$}
		child {node [cst] (n2) {$2$}}
		child {node [var] (n3) {$x$}}
	}
	child {node [elem] (zero) {$-$}
		child {node [var] (n2) {$z$}}
		child {node [elem] (n3) {$\cdot^2$}
			child {node [var] (n4) {$y$}}
		}
	}
;
\node at (-1, -3) {$[0, 20]$}; 
\node at (-2.6, -1.5) {$[0, 40]$}; 
\node at (3, -3) {$[0, 16]$}; 
\node at (7.3, -4.5) {$[-10, 10]$}; 
\node at (7.3, -3) {$[0, 100]$}; 
\node at (5.5, -1.5) {$[-100, 16]$}; 
\path [eval] (-1.1, -2.7) edge (-2.5, -1.8);
\path [eval] (-5.8, -2.7) edge (-2.7, -1.8);
\path [eval] (7.3, -4.2) edge (7.3, -3.2);
\path [eval] (3.2, -2.7) edge (5.5, -1.8);
\path [eval] (7.2, -2.7) edge (5.7, -1.8);
\end{tikzpicture}
\caption{HC4Revise: bottom-up evaluation phase}
\label{fig:hc4-forward}
\end{figure}
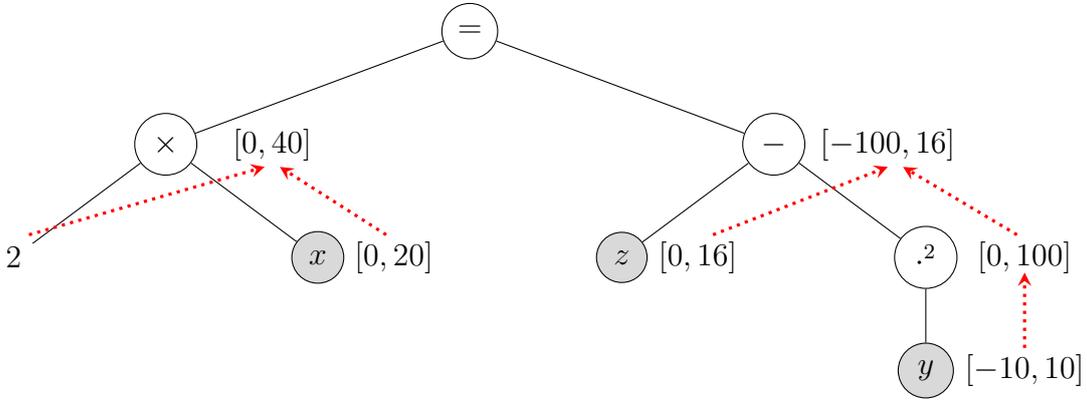

The top-down propagation phase (Figure \ref{fig:hc4-backward}) intersects the ranges of the left-hand ($N_1 = [0, 40]$) and right-hand ($N_3 = [-100, 16]$) sides, then propagates the constraint downwards by evaluating projection functions at each node:
\begin{equation}
\begin{aligned}
N'_1 & = N'_3 = N_1 \cap N_3 = [0, 40] \cap [-100, 16] = [0, 16] \\
X' & = X \cap \frac{N'_1}{2} = [0, 20] \cap [0, 8] = [0, 8] \\
Z' & = Z \cap (N_2 + N'_3) = [0, 16] \cap ([0, 100] + [0, 16]) = [0, 16] \\
N'_2 & = N_2 \cap (Z' - N'_3) = [0, 100] \cap ([0, 16] - [0, 16]) = [0, 16] \\
Y' & = \square \left(Y \cap (-\sqrt{N'_2}), Y \cap \sqrt{N'_2} \right) = \square([-4, 0], [0, 4])
= [-4, 4]
\end{aligned}
\end{equation}

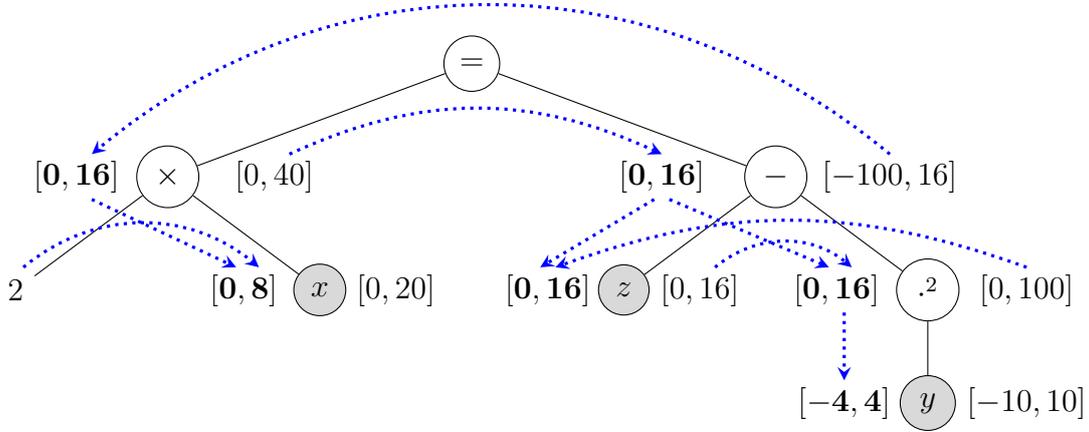
\begin{figure}[htbp]
\centering
\def\svgwidth{\columnwidth}
\begin{tikzpicture}[
	>=stealth,
	level 1/.style={sibling distance=8cm},
	level 2/.style={sibling distance=4cm},
	elem/.style={circle, draw=black, thin, minimum size = 0.5cm},
	var/.style={circle, draw=black, fill=gray!30},
	cst/.style={draw=none},
	line/.style = {->, draw=black},
	prop/.style = {->, dotted, blue, very thick}
	]

\node [elem] (Root) {$=$} 
	child {node [elem] (n4) {$\times$}
		child {node [cst] (n2) {$2$}}
		child {node [var] (n3) {$x$}}
	}
	child {node [elem] (zero) {$-$}
		child {node [var] (n2) {$z$}}
		child {node [elem] (n3) {$\cdot^2$}
			child {node [var] (n4) {$y$}}
		}
	}
;
\node at (-1, -3) {$[0, 20]$}; 
\node at (-2.6, -1.5) {$[0, 40]$}; 
\node at (3, -3) {$[0, 16]$}; 
\node at (7.3, -4.5) {$[-10, 10]$}; 
\node at (7.3, -3) {$[0, 100]$}; 
\node at (5.5, -1.5) {$[-100, 16]$}; 
\node at (-3, -3) {$\mathbf{[0, 8]}$}; 
\node at (-5.2, -1.5) {$\mathbf{[0, 16]}$}; 
\node at (4.9, -4.5) {$\mathbf{[-4, 4]}$}; 
\node at (4.8, -3) {$\mathbf{[0, 16]}$}; 
\node at (2.5, -1.5) {$\mathbf{[0, 16]}$}; 
\node at (1, -3) {$\mathbf{[0, 16]}$}; 
\path [prop] (-5, -1.8) edge (-3.1, -2.7);
\path [prop, bend left=40] (-5.9, -2.7) edge (-2.8, -2.7);
\path [prop, bend left=25] (-2.4, -1.2) edge (2.5, -1.2);
\path [prop, bend right=40] (5.5, -1.2) edge (-5, -1.2);
\path [prop] (2.4, -1.8) edge (0.9, -2.7);
\path [prop] (2.6, -1.8) edge (4.7, -2.7);
\path [prop] (4.9, -3.3) edge (4.9, -4.2);
\path [prop, bend right=20] (7.3, -2.7) edge (1.1, -2.7);
\path [prop, bend left=40] (3.2, -2.7) edge (5, -2.7);
\end{tikzpicture}
\caption{HC4Revise: top-down propagation phase}
\label{fig:hc4-backward}
\end{figure}

The initial box $[0, 20] \times [-10, 10] \times [0, 16]$ was contracted to $[0, 8] \times [-4, 4] \times [0, 16]$ without losing solutions of the constraint.
\end{example}

\subsection{Box consistency}
\label{sec:box-consistency}
The box consistency~\cite{Benhamou1994Clp,Collavizza1999Comparing} defines a coarser consistency than the hull consistency (Definition \ref{def:box}).
However, the algorithms that enforce box consistency tend to deliver a more powerful filtering when the constraints contain several occurrences of the variables. In particular, they are optimal when $c: X \subset \mathbb{R} \rightarrow \mathbb{R}$ is continuous on $X$ with respect to a unique variable $x$ with multiple occurrences. On the contrary, 2B operators are extremely efficient when the variables have a single occurrence.

\begin{definition}[Box consistency]
\label{def:box}
Let $c$ be an $n$-ary constraint, $C$ an interval extension of $c$ and $\bm{X}$ a box. $\bm{X}$ is box-consistent with respect to $c$ if:
\begin{equation}
\forall i \in \{1, \ldots, n\}, \quad X_i = \square{\{x_i \in X_i ~|~ C(X_1, \ldots, [x_i, x_i],
\ldots, X_n) \}}
\end{equation}
$\bm{X}$ is box-consistent with respect to a system of constraints $\mathcal{C}$ if it is box-consistent with respect to each constraint of $\mathcal{C}$.
\end{definition}

\subsection{Interval Newton method}
\label{sec:newton}
Contrary to the classical Newton method, the interval Newton method can provide rigorous bounds on \textit{all the zeros} of a continuous function on an interval $X$. It approximates the box consistency.

Let $z$ be a zero of a function $f$, continuous on an interval $X$ and differentiable on $int(X)$. For all $c \in X$, the mean value theorem states that there exists $\xi$ strictly between $z$ and $c$ such that:
\begin{equation}
0 = f(z) = f(c) + f'(\xi)(z-c)
\end{equation}
The (unknown) value $f'(\xi)$ is rigorously enclosed in $F'(X)$, where $F'$ is an interval extension of $f'$. We get:
\begin{equation}
z \in N_f(X, c) \eqdef c - \frac{f(c)}{F'(X)}
\end{equation}
where $N_f(X, c)$ is called the Newton operator. We build the recurrence relation:
\begin{equation}
\begin{cases}
X_0 	& = X \\
X_{k+1} & = X_k \cap N_f(X_k, c_k), \quad \forall k \ge 0
\end{cases}
\end{equation}

$X_{k+1} = \varnothing$ implies that no zero of $f$ exists in $X_k$, and $X_{k+1} \subset int(X_k)$ proves the existence of a single zero of $f$ in $X_{k+1}$~\cite{Neumaier1990Interval}.

Distinct zeros of $f$ may be automatically separated using extended division (see Example \ref{ex:newton}). In practice however, a fixed point $X_{k+1} = X_k$ is often reached on account of the surestimation of $F'(X)$. \cite{Hansen1992} suggests to bisect $X_k$, then to iterate on both subintervals. All zeros of $f$ can then be bounded with an arbitrary precision.

\begin{remark}
In practice, $f(c)$ is subject to roundoff errors. In order to maintain the conservativity of the computations, $F(c)$ must be used instead.
\end{remark}

If there exists a unique zero $z$ of $f$ in $X$, if $c = m(X)$ and $f$ is monotonic on $X$, the interval Newton method converges Q-quadratically (Definition \ref{def:q-convergence}) to $z$, that is the number of correct digits doubles at each iteration asymptotically~\cite{Hansen1992}.

\begin{definition}[Q-quadratic convergence]
\label{def:q-convergence}
Let $(u_k)_{k \in \mathbb{N}}$ be a sequence and $l$ a real number. $u$ converges Q-quadratically to $l$ if there exists $M > 0$ such that:
\begin{equation}
\lim_{k \rightarrow +\infty} \frac{|u_{k+1} - l|}{|u_k - l|^2} \le M
\end{equation}
\end{definition}

\begin{example}[Interval Newton]
\label{ex:newton}
Let $f(x) = x^2-2$ and $F'(X) = 2X$. We seek the zeros of $f$ on the interval $X_0 = [-3, 2]$ and choose $c_0 = m(X_0) = -0.5$. The first Newton iteration yields:
\begin{equation}
\begin{split}
N_f(X_0, c_0) 	& = c_0 - \frac{f(c_0)}{F'(X_0)} = c_0 - \frac{c_0^2-2}{2X_0}
		= -0.5 - \frac{(-0.5)^2-2}{2[-3, 2]}
		= -0.5 - \frac{-1.75}{[-6, 4]} \\
		& = (-0.5 + [-\infty, -\frac{7}{24}]) \cup (-0.5 + [\frac{7}{16}, +\infty]) \\
	    & = [-\infty, -\frac{19}{24}] \cup [-\frac{1}{16}, +\infty] \\
\end{split}
\end{equation}

\begin{figure}[htbp]
\centering
\small
\def\svgwidth{0.8\columnwidth}
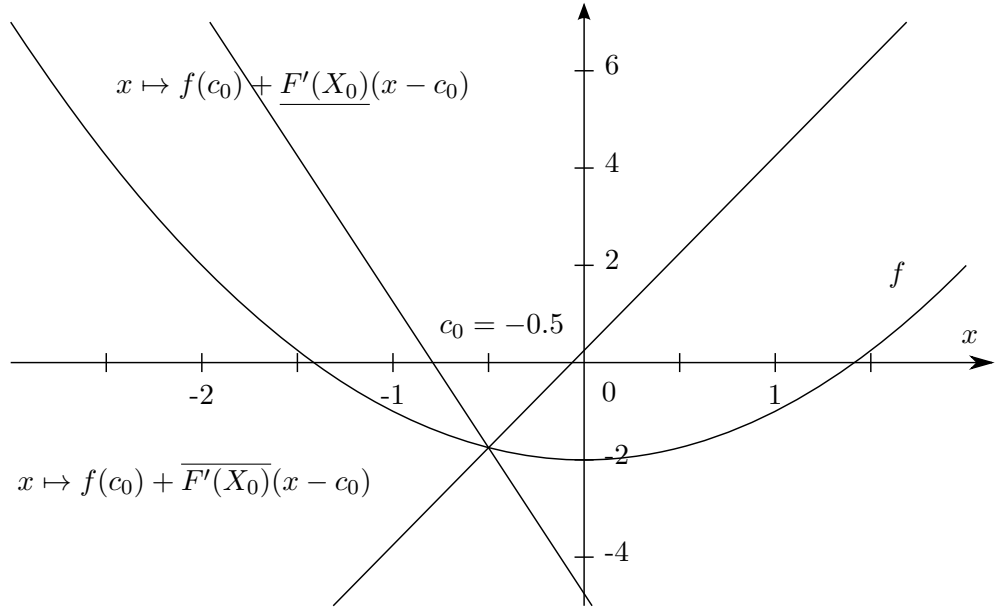
\caption{Interval Newton method with $f: x \mapsto x^2-2$ and $X_0 = [-3, 2]$}
\label{fig:newton}
\end{figure}

$N_f(X_0, c_0)$, composed of two subintervals, represents the intersection of the x-axis and the cone of all tangents (Figure \ref{fig:newton}). Its intersection with the initial interval $X_0$ is:
\begin{equation}
X_0 \cap N_f(X_0, c_0) = [-3, -\frac{19}{24}] \cup [-\frac{1}{16}, 2]
\end{equation}

We now note $X_1 = [-3, -\frac{19}{24}]$ and $X_5 = [-\frac{1}{16}, 2]$ and apply the interval Newton method on both intervals recursively in Table \ref{tab:newton} with a precision of $10^{-4}$.

\begin{table}[htbp]
	\centering
	\footnotesize
	\caption{Interval Newton method with $f: x \mapsto x^2-2$ and $X_0 = [-3, 2]$}
	\begin{tabular}{|c|cccc|c|}
	\hline
	$k$ & $X_k$ & $c_k$ & $F'(X_k)$ & $N_f(X_k, c_k)$ & $X_k \cap N_f(X_k, c_k)$ \\
	\hline
	1 & $[-3, -0.7916]$ & -1.8958 & $[-6, -1.5833]$ & $[-1.6302, -0.8889]$ & 
	$X_2 \eqdef \mathbf{[-1.6302, -0.8889]}$ \\
	2 & $[-1.6302, -0.8889]$ & -1.2596 & $[-3.2603, -1.7779]$ & $[-1.4922, -1.3863]$ &
	$X_3 \eqdef [-1.4922, -1.3863]$ \\
	3 & $[-1.4922, -1.3863]$ & -1.4393 & $[-2.9843, -2.7727]$ & $[-1.4154, -1.4134]$ &
	$X_4 \eqdef [-1.4154, -1.4134]$ \\
	4 & $[-1.4154, -1.4134]$ & -1.4144 & $[-2.8307, -2.8269]$ & $[-1.4143, -1.4142]$ &
	$[-1.4143, -1.4142]$ \\
	\hline
	5 & $[-0.0625, 2]$ & 0.9688 & $[-0.125, 4]$ & $[-\infty, -7.5234] ~\cup$ & $X_6 \eqdef
	[1.2341, 2]$ \\
	& & & & $[1.2341, +\infty]$ & \\
	6 & $[1.2341, 2]$ & 1.6171 & $[2.4682, 4]$ & $[1.3679, 1.4634]$ & $X_7 \eqdef
	\mathbf{[1.3679, 1.4634]}$ \\
	7 & $[1.3679, 1.4634]$ & 1.4156 & $[2.7358, 2.9267]$ & $[1.4141, 1.4143]$ & $X_8 \eqdef
	[1.4141, 1.4143]$ \\
	8 & $[1.4141, 1.4143]$ & 1.4142 & $[2.8283, 2.8286]$ & $[1.4142, 1.4143]$ & $[1.4142,
	1.4143]$ \\
	\hline
	\end{tabular}
	\label{tab:newton}
\end{table}

The two zeros of $f$ on $X_0$ were bounded by the intervals $[-1.414213693, -1.414213455]$ (width $2.4 \cdot 10^{-7}$) and $[1.414213562, 1.414213562]$ (width $6.1 \cdot 10^{-11}$) respectively.
Their existence is guaranteed during iterations 1 and 6 (in bold): $X_2 \eqdef X_1 \cap N_f(X_1, c_1) \subset int(X_1)$ and $X_7 \eqdef X_6 \cap N_f(X_6, c_6) \subset int(X_6)$.
\end{example}

\subsection{Monotonicity-based contractors}
\label{sec:mohc}
The adaptive algorithm Mohc (MOnotonic Hull Consistency)~\cite{Araya2010} exploits the local monotonicity in order to more efficiently contract a box with respect to a system of equations where the variables have multiple occurrences. It exploits the following enclosure (see Definition \ref{def:monotonicity-extension}):
\begin{equation}
c(\bm{X}) \subset C_M(\bm{X}) \eqdef [\underline{C(\bm{X}^-)}, \overline{C(\bm{X}^+)}]
\end{equation}

Consequently, the equation $c = 0$ can be decomposed into two constraints $C(\bm{X}^-) \le 0$ and $C(\bm{X}^+) \ge 0$. Mohc is a fixed-point algorithm whose revised procedure MohcRevise (Algorithme \ref{algo:mohc}) invokes:
\begin{itemize}
\item HC4Revise($C(\bm{X}^-) \le 0$) in order to contract the variables of $\mathcal{V}_{nm}$ ;
\item HC4Revise($C(\bm{X}^+) \ge 0$) in order to contract the variables of $\mathcal{V}_{nm}$ ;
\item a version of the interval Newton method (see Section \ref{sec:newton}), called MonotonicBoxNarrow, in order to contract the variables of $\mathcal{V}_m$.
\end{itemize}
where $\mathcal{V}$ is the set of all the variables of $c$, $\mathcal{V}_m \subset \mathcal{V}$ the set of variables with multiple occurrences with respect to which $c$ is monotonic on $\bm{X}$ and $\mathcal{V}_{nm} = \mathcal{V} \setminus \mathcal{V}_m$ the rest of the variables.

\begin{algorithm}[htbp]
\caption{MohcRevise}
\label{algo:mohc}
\begin{algorithmic}[]
\Function{MohcRevise}{$c = 0$: equality constraint, \inout{} $\bm{X}$: box, $\mathcal{V}$: variables}
\State \Call{HC4Revise}{$C(\bm{X}) = 0$}
\If{$\mathcal{V}_{nm} \neq \varnothing$}
	\State \Call{HC4Revise}{$C(\bm{X}^-) \le 0$}
	\State \Call{HC4Revise}{$C(\bm{X}^+) \ge 0$}
	\For{$x_i \in \mathcal{V}_m$}
		\State $X_i \gets$ \Call{MonotonicBoxNarrow}{$C, \bm{X}, i$}
	\EndFor
\EndIf
\EndFunction
\end{algorithmic}
\end{algorithm}

Octum, an algorithm independently devised by~\cite{Chabert2009Octum}, is identical to MonotonicBoxNarrow when the function is monotonic with respect to all its variables.

\section{Strong consistencies}
\label{sec:strong-consistencies}
The 2B and box consistencies are so-called weak consistencies, since they define a consistency on the bounds of the variables for an individual constraint. Consequently, the resulting filtering with respect to the system of constraints may be poor. The stronger partial consistencies 3B and CID invoke the 2B and box operators as subcontractors during a shaving process and produce tighter contractions. Shaving is the temporary assignment to a variable of a small portion (a slice) of its domain ; slices can then be contracted or discarded from the domain.
Unlike branching algorithms that have an exponential complexity, shaving is a polynomial refutation technique. However, it is not incremental, in that the contraction of the domain of a variable requires that the filtering be performed again on the other variables.

The 3B algorithm (Section \ref{sec:3b}) consists in discarding the extremal inconsistent slices. The process is interrupted when a bound cannot be reduced.
The CID algorithm (Section \ref{sec:cid}) handles all the slices in order to discard the values that are inconsistent in all the slices simultaneously.

\subsection{3B consistency}
\label{sec:3b}
The 3B consistency (Definition \ref{def:3b-consistency}) is a shaving-based relaxation of the \textit{singleton arc consistency} that invokes a 2B subcontractor in order to discard the slices~\cite{Lhomme1993Consistency}. It enforces the consistency of the bounds of the variables with respect to the system of constraints.
The 3B consistency can be recursively extended to the $k$B consistency ($k > 2$) by invoking a $(k-1)$B contractor on the slices.

\begin{definition}[3B($s_{3B}$) consistency]
Let:
\begin{itemize}
\item $P = (\mathcal{V}, \mathcal{C}, \bm{X})$ be a \gls{NCSP} ;
\item $OC_{2B}$ be a 2B contractor ;
\item $s_{3B}$ be the number of slices ;
\item $\bm{X_i^1}$ be the leftmost slice of the domain of $x_i \in \mathcal{V}$, that is the subbox of $\bm{X}$ whose $i$th component is $\displaystyle [\underline{X_i}, \underline{X_i} + \frac{w(X_i)}{s_{3B}}]$ ;
\item $\bm{X_i^{s_{3B}}}$ be the rightmost slice of the domain of $x_i \in \mathcal{V}$, that is the subbox of $\bm{X}$ whose $i$th component is $\displaystyle [\overline{X_i} - \frac{w(X_i)}{s_{3B}}, \overline{X_i}]$.
\end{itemize}
$\bm{X}$ is 3B($s_{3B}$)-consistent with respect to $c \in \mathcal{C}$ if:
\begin{equation}
\forall i \in \{1, \ldots, n\}, \quad OC_{2B}(\bm{X_i^1}, c) \neq \varnothing \land
OC_{2B}(\bm{X_i^{s_{3B}}}, c) \neq \varnothing
\end{equation}
$\bm{X}$ is 3B($s_{3B}$)-consistent with respect to $\mathcal{C}$ if it is 3B($s_{3B}$)-consistent with respect to each constraint of $\mathcal{C}$.
\label{def:3b-consistency}
\end{definition}

\begin{figure}[htb]
\centering
\small
\def\svgwidth{0.8\columnwidth}
\begingroup%
  \makeatletter%
  \providecommand\color[2][]{%
    \errmessage{(Inkscape) Color is used for the text in Inkscape, but the package 'color.sty' is not loaded}%
    \renewcommand\color[2][]{}%
  }%
  \providecommand\transparent[1]{%
    \errmessage{(Inkscape) Transparency is used (non-zero) for the text in Inkscape, but the package 'transparent.sty' is not loaded}%
    \renewcommand\transparent[1]{}%
  }%
  \providecommand\rotatebox[2]{#2}%
  \newcommand*\fsize{\dimexpr\f@size pt\relax}%
  \newcommand*\lineheight[1]{\fontsize{\fsize}{#1\fsize}\selectfont}%
  \ifx\svgwidth\undefined%
    \setlength{\unitlength}{515.80513bp}%
    \ifx\svgscale\undefined%
      \relax%
    \else%
      \setlength{\unitlength}{\unitlength * \real{\svgscale}}%
    \fi%
  \else%
    \setlength{\unitlength}{\svgwidth}%
  \fi%
  \global\let\svgwidth\undefined%
  \global\let\svgscale\undefined%
  \makeatother%
  \begin{picture}(1,0.29573185)%
    \lineheight{1}%
    \setlength\tabcolsep{0pt}%
    \put(0,0){\includegraphics[width=\unitlength,page=1]{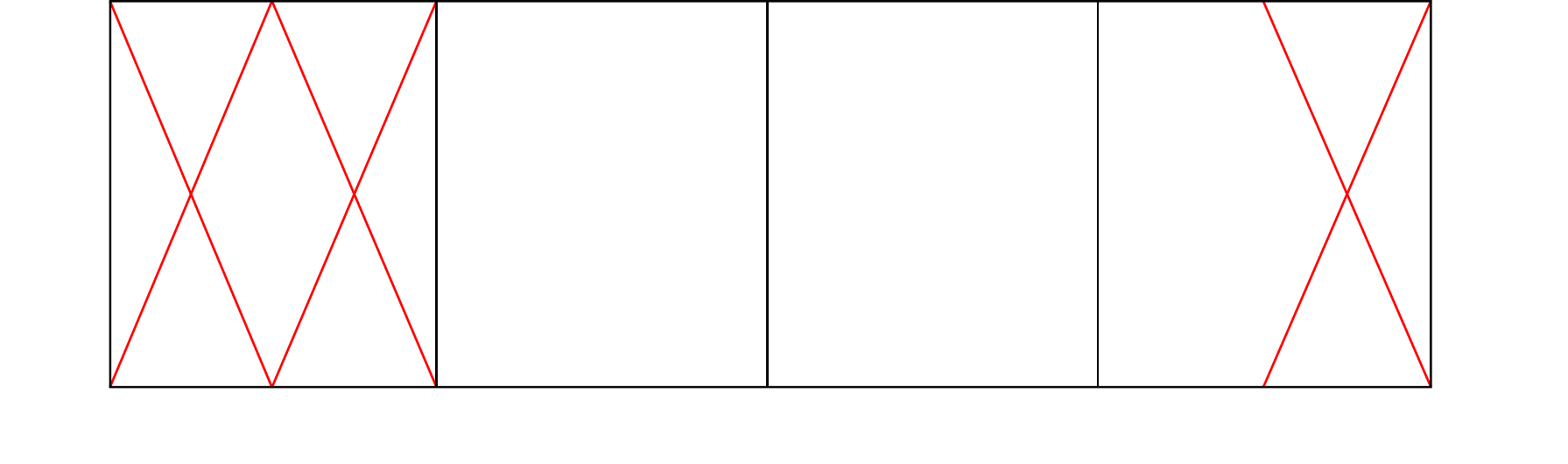}}%
    \put(0.10491335,0.06645341){\color[rgb]{0,0,0}\makebox(0,0)[lt]{\lineheight{0}\smash{\begin{tabular}[t]{l}$\mathbf{X_1^1}$\end{tabular}}}}%
    \put(0,0){\includegraphics[width=\unitlength,page=2]{figures/3b.pdf}}%
    \put(0.47871707,0.00272631){\color[rgb]{0,0,0}\makebox(0,0)[lt]{\lineheight{0}\smash{\begin{tabular}[t]{l}$X_1$\end{tabular}}}}%
    \put(-0.00095967,0.17081528){\color[rgb]{0,0,0}\makebox(0,0)[lt]{\lineheight{0}\smash{\begin{tabular}[t]{l}$X_2$\end{tabular}}}}%
    \put(0.83746652,0.06645387){\color[rgb]{0,0,0}\makebox(0,0)[lt]{\lineheight{0}\smash{\begin{tabular}[t]{l}$\mathbf{X_1^8}$\end{tabular}}}}%
  \end{picture}%
\endgroup%

\caption{3B shaving with $s_{3B} = 8$}
\label{fig:3b}
\end{figure}

The algorithm that enforces the 3B consistency refutes slices at the bounds of the variables (Figure \ref{fig:3b}). The domain $X_i$ of the variable $x_i$ is temporarily instantiated to a subinterval of $X_i$ (a slice) ; a subcontractor is then invoked on the subproblem with respect to the system of constraints. If an inconsistency is detected, the slice can be removed from the domain $X_i$ (Algorithm \ref{algo:3b}).

\begin{algorithm}[htb]
\caption{3B shaving}
\label{algo:3b}
\begin{algorithmic}
\Function{Variable3B}{\inout{} $\bm{X}$: box, $i$: composante, $\mathcal{C}$: system of constraints, $OC_{2B}$: 2B contractor, $s_{3B}$: number of slices}
\State Invoke $OC_{2B}$ to discard the leftmost slices of $X_i$
\State Invoke $OC_{2B}$ to discard the rightmost slices of $X_i$
\EndFunction

\Function{3BShaving}{\inout{} $\bm{X}$: box, $\mathcal{C}$: system of constraints, $OC_{2B}$: 2B contractor, $s_{3B}$: number of slices}
\For{$i \in \{1, \ldots, n\}$}
\Comment for each variable
	\State \Call{Variable3B}{$\bm{X}, i, \mathcal{C}, OC_{2B}, s_{3B}$}
\EndFor
\EndFunction

\Function{3B}{\inout{} $\bm{X}$: box, $\mathcal{C}$: system of constraints, $OC_{2B}$: 2B contractor, $s_{3B}$: number of slices}
\State \Call{FixedPoint}{$\bm{X}, \mathcal{C}$, \Call{3BShaving}{$\bm{X}, \mathcal{C}, OC_{2B},
s_{3B}$}}
\EndFunction
\end{algorithmic}
\end{algorithm}

\subsection{CID consistency}
\label{sec:cid}
Constructive disjunction on \gls{CSP} handles a disjunction of constraints $c_1 \lor \ldots \lor c_m$ by branching alternately on each constraint $c_i$ and by discarding values that are inconsistent with respect to all the constraints. \Gls{CID} (Definition \ref{def:cid}) is a shaving-based technique that discards values that are inconsistent in all the slices~\cite{Trombettoni2007Constructive}.

\newpage

\begin{definition}[CID($s_{CID}$) consistency]
Let:
\begin{itemize}
\item $P = (\mathcal{V}, \mathcal{C}, \bm{X})$ be a \gls{NCSP} ;
\item $OC$ be a partial consistency operator ;
\item $s_{CID}$ be the number of slices ;
\item $\bm{X_i^k}$ be the $k$th slice of the domain of $x_i \in \mathcal{V}$, that is the subbox of $\bm{X}$ whose $i$th component is $[\underline{X_i} + (k-1)\frac{w(X_i)}{s_{CID}}, \underline{X_i} + k\frac{w(X_i)}{s_{CID}}]$.
\end{itemize}
$x_i$ is \emph{CID($s_{CID}$)-consistent} with respect to $P$ and $OC$ if:
\begin{equation}
\bm{X} = \hull_{k = 1}^{s_{CID}} OC(\bm{X_i^k}, \mathcal{C})
\end{equation}
$P$ is \emph{CID($s_{CID}$)-consistent} if all the variables of $\mathcal{V}$ are CID($s_{CID}$)-consistent.
\label{def:cid}
\end{definition}

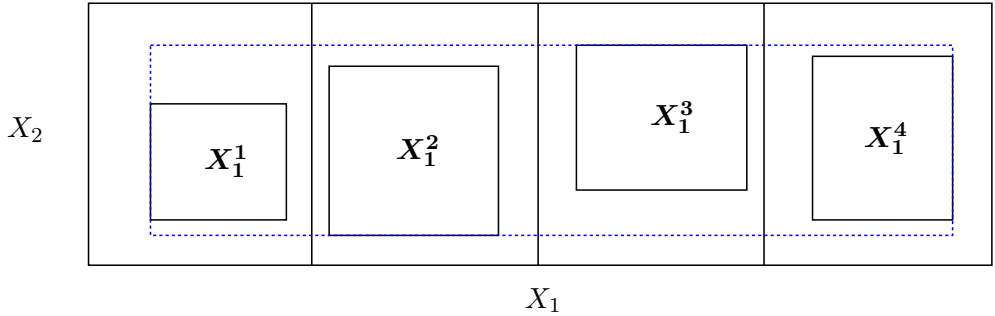
\begin{figure}[htbp]
\centering
\small
\def\svgwidth{0.8\columnwidth}
\begingroup%
  \makeatletter%
  \providecommand\color[2][]{%
    \errmessage{(Inkscape) Color is used for the text in Inkscape, but the package 'color.sty' is not loaded}%
    \renewcommand\color[2][]{}%
  }%
  \providecommand\transparent[1]{%
    \errmessage{(Inkscape) Transparency is used (non-zero) for the text in Inkscape, but the package 'transparent.sty' is not loaded}%
    \renewcommand\transparent[1]{}%
  }%
  \providecommand\rotatebox[2]{#2}%
  \newcommand*\fsize{\dimexpr\f@size pt\relax}%
  \newcommand*\lineheight[1]{\fontsize{\fsize}{#1\fsize}\selectfont}%
  \ifx\svgwidth\undefined%
    \setlength{\unitlength}{473.53591919bp}%
    \ifx\svgscale\undefined%
      \relax%
    \else%
      \setlength{\unitlength}{\unitlength * \real{\svgscale}}%
    \fi%
  \else%
    \setlength{\unitlength}{\svgwidth}%
  \fi%
  \global\let\svgwidth\undefined%
  \global\let\svgscale\undefined%
  \makeatother%
  \begin{picture}(1,0.31523406)%
    \lineheight{1}%
    \setlength\tabcolsep{0pt}%
    \put(0,0){\includegraphics[width=\unitlength,page=1]{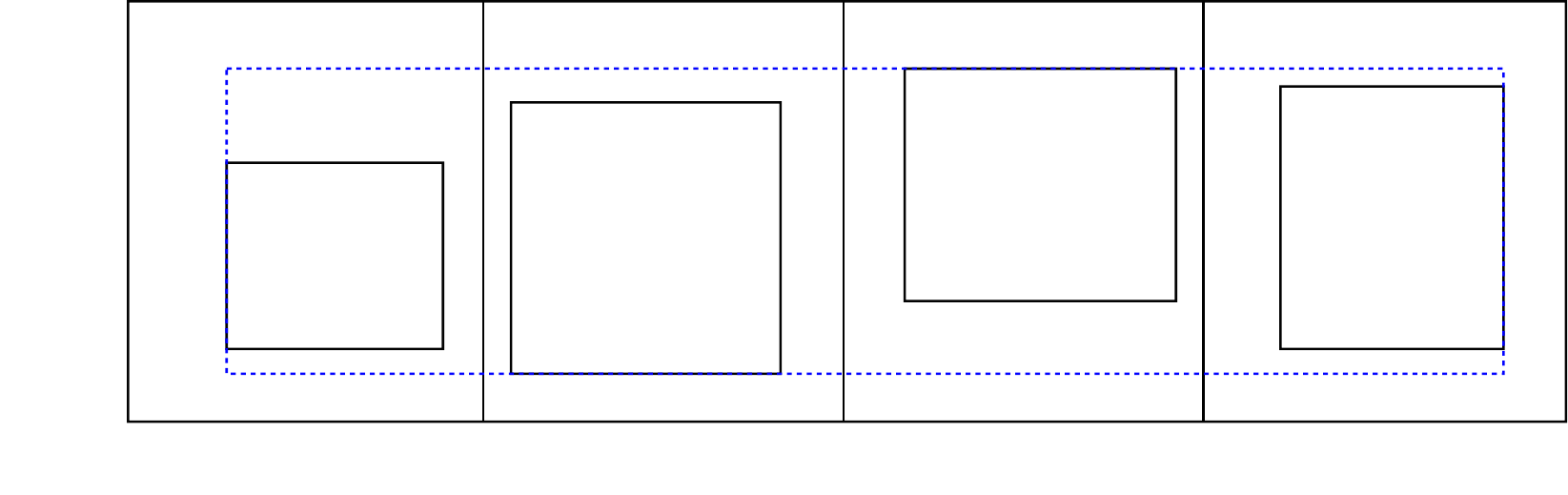}}%
    \put(0.52488649,0.00296967){\color[rgb]{0,0,0}\makebox(0,0)[lt]{\lineheight{0}\smash{\begin{tabular}[t]{l}$X_1$\end{tabular}}}}%
    \put(0.19927283,0.1446581){\color[rgb]{0,0,0}\makebox(0,0)[lt]{\lineheight{0}\smash{\begin{tabular}[t]{l}$\bm{X_1^1}$\end{tabular}}}}%
    \put(0.39375357,0.15497136){\color[rgb]{0,0,0}\makebox(0,0)[lt]{\lineheight{0}\smash{\begin{tabular}[t]{l}$\bm{X_1^2}$\end{tabular}}}}%
    \put(0.65220314,0.1893694){\color[rgb]{0,0,0}\makebox(0,0)[lt]{\lineheight{0}\smash{\begin{tabular}[t]{l}$\bm{X_1^3}$\end{tabular}}}}%
    \put(0.86908739,0.16777063){\color[rgb]{0,0,0}\makebox(0,0)[lt]{\lineheight{0}\smash{\begin{tabular}[t]{l}$\bm{X_1^4}$\end{tabular}}}}%
    \put(-0.00104533,0.17707589){\color[rgb]{0,0,0}\makebox(0,0)[lt]{\lineheight{0}\smash{\begin{tabular}[t]{l}$X_2$\end{tabular}}}}%
  \end{picture}%
\endgroup%

\caption{CID shaving with $s_{CID} = 4$}
\label{fig:cid}
\end{figure}

The 3B algorithm refutes inconsistent values at the bounds of the domain of $x_i$, however the effort spent to contract the other variables is wasted. The substantial advantage of \gls{CID} (Algorithm \ref{algo:cid}) over 3B is that it partially exploits the information of all the contracted domains: the initial box is replaced by the convex hull of the $s_{CID}$ contracted slices. Figure \ref{fig:cid} illustrates the contraction of 4 slices on the domain of $x_1$. Unlike a 3B shaving, the convex hull operation of the contracted slices partially preserves the contraction of $x_2$.

The CID($s_{CID}$) consistency is difficult to obtain on account of the slow convergence of the fixed-point algorithm. In practice, only a quasi-fixed point with a precision $\eta$ is computed: $x_i$ is CID($s_{CID}, \eta$)-consistent if the convex hull of the contracted slices does not contract any variable more than $\eta$. However, experimental results suggest that the fixed-point algorithm is of little benefit and that $s_{CID} = 2$ generally produces satisfactory results.


\newpage

\begin{algorithm}[htbp]
\caption{CID shaving}
\label{algo:cid}
\begin{algorithmic}
\Function{VariableCID}{\inout{} $\bm{X}$: box, $i$: composante, $\mathcal{C}$: system of constraints, $OC$: contractor, $s_{CID}$: number of slices}
\State $\bm{X}_{\square} \gets \varnothing$
\For{$k \in \{1, \ldots, s_{CID}\}$}
\Comment for each slice
	\State $\bm{B}_k \gets$ \Call{slice}{$\bm{X}, i, k, s_{CID}$}
	\State $\bm{X}_{\square} \gets \square(\bm{X}_{\square}, OC(\bm{B}_k, \mathcal{C}))$
	\Comment convex hull
\EndFor
\State $\bm{X} \gets \bm{X}_{\square}$
\EndFunction

\Function{CIDShaving}{\inout{} $\bm{X}$: box, $\mathcal{C}$: system of constraints, $OC$: contractor, $s_{CID}$: number of slices}
\For{$i \in \{1, \ldots, n\}$}
\Comment for each variable
	\State \Call{VariableCID}{$\bm{X}, i, \mathcal{C}, OC, s_{CID}$}
\EndFor
\EndFunction

\Function{CID}{\inout{} $\bm{X}$: box, $\mathcal{C}$: system of constraints, $OC$: contractor, $s_{CID}$: number of slices}
\State \Call{FixedPoint}{$\bm{X}, \mathcal{C}$, \Call{CIDShaving}{$\bm{X}, \mathcal{C}, OC,
s_{CID}$}}
\EndFunction
\end{algorithmic}
\end{algorithm}

The hybrid algorithm 3BCID~\cite{Trombettoni2007Constructive} combines the CID and 3B shavings and reconciles the high number of 3B slices $s_{3B}$ and the low number of CID slices $s_{CID}$ (in general, $s_{CID} \ll s_{3B}$). Its filtering power is greater than those of 3B and CID.

\section{Global consistencies}
\label{sec:global-consistencies}

The contractors than handle the constraints simultaneously are generally linearization techniques:
\begin{itemize}
\item optimization problems subject to a square system of equality constraints can be solved using a multivariate interval Newton algorithm (Section \ref{sec:multivariate-newton}) ;
\item optimization problems subject to (in)equality constraints can be convexified (Section \ref{sec:convexification}).
\end{itemize}

\subsection{Multivariate interval Newton method}
\label{sec:multivariate-newton}

The interval Newton method (Section \ref{sec:newton}) can be extended to a square system of $n$ equations and $n$ variables. Let $g: \mathbb{R}^n \rightarrow \mathbb{R}^n$ be a differentiable vector-valued function, $\bm{X}$ a box and $\bm{x_k} \in \bm{X}$ a point.
The equation $g(x) = 0$ can be linearized at the point $\bm{x_k}$ using the mean value form (Equation \ref{eq:mean-value-form}):
\begin{equation}
g(\bm{x_k}) + J_g(\bm{X})(\bm{X}-\bm{x_k}) = 0
\label{eq:newton-systeme}
\end{equation}
where $J_g$ is an interval extension of the Jacobian of $g$ (the matrix of its first-order partial derivatives) on $\bm{X}$. By solving for $\bm{Y_k} = \bm{X} - \bm{x_k}$, we obtain:
\begin{equation}
\bm{Y_k} = -J_g(\bm{X})^{-1} g(\bm{x_k})
\end{equation}
In practice, $\bm{Y_k}$ is obtained by solving a preconditioned version of the interval linear system $J_g(\bm{X}) \bm{Y_k} = -g(\bm{x_k})$.

\subsection{Convexification}
\label{sec:convexification}

\begin{figure}[htbp]
\centering
\def\svgwidth{0.6\columnwidth}
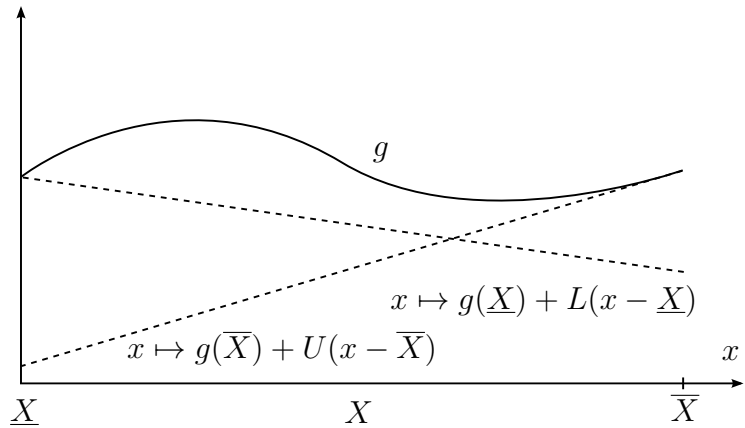
\caption[Convex linear underestimators of a function using corners as expansion points]{Convex linear underestimators $G_{mv}(X, \underline{X})$ and $G_{mv}(X, \overline{X})$ of a function $g$ using corners of $X$ as expansion points of the mean value form}
\label{fig:lp}
\end{figure}

Linearizing a function using the mean value form (Equation \ref{eq:mean-value-form}) at an interior expansion point $x_k \in int(X)$ produces underestimators that are not convex. Fortunately, the underestimators become convex when the expansion point is a corner of the box. \cite{Mentzer1991Lp} thus generates a linear program by computing convex linear underestimators of the objective function and the constraints on a box $\bm{X}$ (Figure \ref{fig:lp}).

Let $\frac{\partial G}{\partial x_i}(\bm{X})$ be an interval enclosure of the $i$th partial derivative of $g$ on $\bm{X}$. If we take the leftmost corner of $\bm{X}$ as an expansion point, we obtain:
\begin{equation}
\forall \bm{x} \in \bm{X}, \quad g(\underline{\bm{X}}) + \sum_{i=1}^n \underline{\frac{\partial G}{\partial x_i}(\bm{X})} \cdot (x_i - \underline{X_i}) \le g(\bm{x})
\label{eq:minoration-gauche}
\end{equation}

We can compute a convex linear underestimator for the objective function and the constraints of the problem, which results in the following linear problem:
\begin{equation}
\begin{aligned}
(\mathcal{P}_{lb}) \quad & \min_{\bm{x} \in \bm{X}} & \sum_{i=1}^n \underline{\frac{\partial F}{\partial x_i}(\bm{X})}
\cdot x_i & & \\
	& \text{s.t.} 	& g_j(\underline{\bm{X}}) + \sum_{i=1}^n \underline{\frac{\partial G_j}{\partial
		x_i}(\bm{X})} \cdot (x_i - \underline{X_i}) \le 0, & \quad \forall j \in \{1, \ldots, m\}
\end{aligned}
\label{eq:convex-relaxation}
\end{equation}

The relaxed problem $\mathcal{P}_{lb}$ (Equation \ref{eq:convex-relaxation}) can be solved using linear programming techniques. If it is infeasible, there exists no solution to the original problem, since the polytope of the relaxed constraints contains the feasible set of the original problem. If it is feasible, the optimum is a lower bound of the initial problem on $\bm{X}$.

\cite{Sotiropoulos2005Optimal} combined a similar approach with the dynamic constraint $f \le \tilde{f}$ to reduce the ranges of the variables. \cite{Yamamura1998Interval,Yamamura2002Finding} replaced the nonlinear terms in the constraints with their interval ranges, then reduced the ranges of the variables with $2n$ calls to the simplex method.
The X-Newton algorithm~\cite{Araya2012Contractor} exploits the two aforementioned approaches with $2n+1$ calls to the simplex method: it computes a lower bound of the original problem and contracts the domains of the variables.
Each inequality constraint is relaxed by two underestimators whose expansion points are a randomly selected corner and its opposite corner. This algorithm converges linearly to the solution of the original problem.
\cite{Ninin2010Reliable} relaxed the initial problem using \gls{AA} and computed a lower bound of the initial problem with one call to the simplex method.

\chapter{Charibde}
\label{chap:charibde}

\minitoc

The aforementioned global optimization methods seem to handle distinct classes of problems: interval methods provide a numerical guarantee of optimality, however are limited to instances with a few dozen variables at most, while \glspl{EA} shine on large multimodal instances for which traditional methods struggle to converge.

In this chapter, we discuss synergetic schemes to reconcile both approaches in order to:
\begin{itemize}
\item quickly explore the search space ;
\item prevent premature convergence towards local minima ;
\item discard suboptimal or infeasible subspaces ;
\item certify the optimality of the solution.
\end{itemize}
and address problems that were up to now deemed as intractable.

Section \ref{sec:hybridization} is a summary of existing hybrid methods and an introduction of Charibde, a cooperative hybrid solver that combines a deterministic method and a stochastic method.
Sections \ref{sec:charibde-ibc} and \ref{sec:charibde-de} describe in detail how both components are implemented in Charibde. In Section \ref{sec:coconut-comparison}, we provide a preliminary comparison of Charibde against state-of-the-art rigorous and non rigorous solvers on a subset of difficult problems.

\section{Hybridization of optimization techniques}
\label{sec:hybridization}
Hybridizing optimization techniques boils down to dividing up exploration tasks between several algorithms (usually two). In the following, we describe heterogeneous hybridizations (a metaheuristic and another algorithm).
The taxonomy of hybrid metaheuristics was addressed by several authors~\cite{Talbi2002Taxonomy, Puchinger2005Combining, Alba2005Parallel, Raidl2006Unified, Jourdan2009Hybridizing} ; we present in this section a hierarchical classification that summarizes possible hybridization strategies between discrete or continuous optimization techniques.

\cite{Puchinger2005Combining} list three distinguishing criteria: the nature of the hybridization (with a metaheuristic, an artificial intelligence technique or an operations research technique), the level of hybridization (the components that are exchanged) and the order of execution (sequential or parallel executions). These criteria suggest the following classification:
\begin{enumerate}
\item integrative or low-level methods: a particular operator of one algorithm is replaced with another algorithm (such as a local search or an exact method) following a master-slave scheme. There are two different schemes:
	\begin{itemize}
	\item the metaheuristic is the master ;
	\item the metaheuristic is the slave.
	\end{itemize}
\item cooperative or high-level methods: the algorithms are not embedded within one another, but are independent and exchange information. There are two different schemes:
	\begin{itemize}
	\item sequential executions (HRH): one of the methods is run before the other (much like a preprocessing step) ;
	\item parallel or intertwined executions (HCH): the two methods run in parallel.
	\end{itemize}
\end{enumerate}

The hybridization is called:
\begin{itemize}
\item \textit{global} when the methods explore the whole search space, and \textit{partial} when they are limited to a particular subspace ;
\item \textit{general} if both methods attempt to solve the same problem, and \textit{specialized} if they handle different problems.
\end{itemize}

Sections \ref{sec:integrative-methods} and \ref{sec:cooperative-methods} list integrative and cooperative hybridizations from the literature, respectively.

\subsection{Integrative methods}
\label{sec:integrative-methods}
Integrative methods embed one algorithm within the other by substituting a particular operator.

\subsubsection{Slave metaheuristic}
\cite{Zhang2007New} embed a \gls{GA} within an \gls{IBB} algorithm. The \gls{GA} provides the direction along which the boxes are partitioned, and an individual is generated within each subspace. When a box is discarded, the corresponding individuals are removed from the population. The evaluation of the best individual updates the best known upper bound of the global minimum at each generation.

\subsubsection{Master metaheuristic}
Memetic algorithms~\cite{Moscato2004Memetic} are \glspl{EA} that embed a local search algorithm (gradient descent, tabu search, Nelder-Mead simplex, simulated annealing) in order to improve individuals locally.

In \cite{Cotta2003Embedding}, the crossover operator is replaced with a \gls{BB} algorithm that considers all possible offspring that can be generated from two given parents, and selects the best combination.

\subsection{Cooperative methods}
\label{sec:cooperative-methods}
Cooperative methods are characterized by the granularity of the hybridization, the software implementation, the memory (shared or distributed) and the synchronization of processes.

\subsubsection{Sequential executions}
\cite{Feltl2004Improved} solve the generalized assignment problem using a hybrid \gls{GA}. They first solve a linear relaxation of the problem with CPLEX, then they round the optimal solution to construct promising individuals. Stochastic repair operators are designed to generate feasible individuals.

\cite{Sotiropoulos1997New} combine an \gls{IBB} algorithm and a \gls{GA}. The \gls{IBB} generates a list $\mathcal{L}$ of candidates boxes of size at most $\delta$. The \gls{GA} population is initialized by generating an individual in each of the boxes of $\mathcal{L}$. In order to keep the population size below 50, numerical results suggest that $\delta \in [10^{-3}, 10^{-1}]$. The authors do not use the midpoint test: the best known upper bound $\tilde{f}$ of the global minimum is updated using the coarse upper bound $\overline{F(\bm{X})}$.

\subsubsection{Parallel executions}
\cite{Gallardo2007Hybridization} (hybridization of a \gls{BB} and a memetic algorithm) and \cite{Blum2011Hybrid} (hybridization of a beam search and a memetic algorithm) describe similar parallel strategies: the deterministic method identifies promising subspaces of the search space, that are then explored by the metaheuristic. The memetic algorithm provides the tree search algorithm with an upper bound of the global minimum in order to reduce the search space, and is in return presented with promising regions of the search space.

\cite{Cotta1995Hybridizing} mentions solving the travelling salesman problem by combining a \gls{BB} and a \gls{GA}. The \gls{GA} provides the \gls{BB} with an upper bound of the global minimum, in order to discard suboptimal subspaces. The \gls{BB} injects promising paths into the \gls{GA} population. The authors however are critical of the approach: exchanging information between processes that run at different speeds is problematic. In particular, injecting promising solutions at the start of the \gls{GA} may result in "superindividuals", which may reduce the diversity within the population. Two alternatives are considered: an integrative hybridization (substitute the crossover operator with a tree search algorithm) or a master-slave scheme (a \gls{BB} and $m$ \gls{GA} in parallel).

\cite{Alliot2012Finding} propose a bound constrained solver that combines a \gls{GA} and an \gls{IBB} that run independently and communicate through shared memory. The \gls{GA} performs a fast exploration of the search space in order to discover a promising solution ; the corresponding upper bound of the global minimum is then sent to the \gls{IBB} to intensify the pruning. Whenever the \gls{IBB} finds a punctual solution that improves the best known upper bound of the global minimum, it is injected into the \gls{GA} population in order to prevent convergence towards local minima. A third process, triggered periodically, projects out-of-domain individuals (that do not belong to any box of the \gls{IBB}) into the closest box.
Their generic framework does not require continuity, differentiability or even factorability of the objective function ; it must simply be computable with \gls{IA}, that is an interval evaluation procedure must be available (possibly a black box). The authors present new optimal results for the Michalewicz function ($n = 12$) and the rotated Griewank function ($n = 6$).

\subsection{Charibde: a cooperative approach}
Our work exploits the cooperative scheme introduced in~\cite{Alliot2012Finding}. However, the efficiency and reliability of their approach is limited:
\begin{itemize}
\item the solver is limited to bound constrained optimization ;
\item their interval techniques remain naive and are not competitive with state-of-the-art solvers ;
\item the projection of out-of-domain individuals is not efficient ;
\item the \gls{GA} is not rigorous, that is the evaluation of the best individual sent to the \gls{IBB} may be subject to roundoff errors ;
\item the boxes whose diameter is lower than a given threshold are discarded without being further explored: solutions may therefore be lost.
\end{itemize}

Our hybrid algorithm \textbf{Charibde} (Cooperative Hybrid Algorithm using Reliable Interval-Based methods and Differential Evolution)~\cite{Vanaret2013Preventing,Vanaret2015Hybridization} combines a \gls{DE} algorithm and an \gls{IBC} algorithm. Although it embeds stochastic components, Charibde is a \textbf{fully reliable} solver.

\begin{figure}[htbp!]
\centering
\def\svgwidth{0.85\columnwidth}
\input{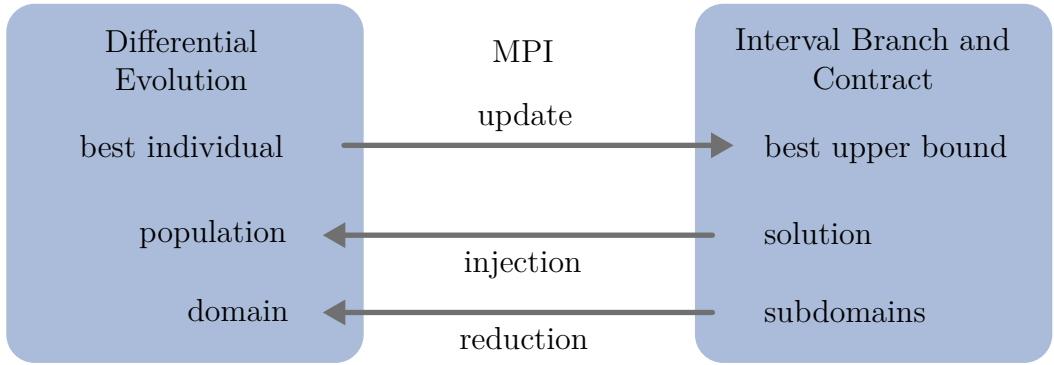}
\caption{Cooperative hybridization in Charibde}
\label{fig:charibde}
\end{figure}

The \gls{DE} algorithm was chosen over the \gls{GA} on account of its convincing performances on continuous problems and its low number of hyperparameters. It communicates with an \gls{IBC} algorithm that benefits from state-of-the-art \gls{ICP} techniques. A new exploration heuristic periodically reduces the domain of the \gls{DE}. Bounds, solutions and domain are exchanged using MPI (Figure \ref{fig:charibde}). 


Sections \ref{sec:charibde-ibc} and \ref{sec:charibde-de} describe our implementations of the \gls{IBC} algorithm and the \gls{DE} algorithm, respectively. A comparison of Charibde against state-of-the-art solvers is given in Section \ref{sec:coconut-comparison}.

\section{Interval branch and contract algorithm}
\label{sec:charibde-ibc}

\subsection{Main framework}

Algorithm \ref{alg:charibde-ibc} describes the \gls{IBC} algorithm implemented in OCaml in Charibde. It embeds MPI routines to exchange information with the \gls{DE}: at each iteration, the \gls{IBC} recovers the best feasible solution of the \gls{DE} if it is available, and possibly updates the best known upper bound $\tilde{f}$ of the global minimum (in red).

\begin{algorithm}[htbp!]
\caption{Charibde: interval branch and contract algorithm}
\label{alg:charibde-ibc}
\begin{algorithmic}
\Function{IntervalBranchAndContract}{$F$: objective function, $\mathcal{C}$: set of constraints, $\bm{D}$: domain}
\State $(\tilde{\bm{x}}, \tilde{f}) \gets (\varnothing, +\infty)$
\Comment best known upper bound
\State $\mathcal{Q} \gets \{ (\bm{D}, -\infty, \mathcal{C}) \}$
\Comment priority queue
\While{$\mathcal{Q} \neq \varnothing$}
	\textcolor{amred}{
	\State $(\bm{x}_{ED}, f_{ED}) \gets$ \texttt{MPI\_ReceiveFromDE}$()$
	\If{$f_{ED} < \tilde{f}$}
		\State $(\tilde{\bm{x}}, \tilde{f}) \gets (\bm{x}_{ED}, f_{ED})$
	\EndIf
	}
	\State Extract an element $(\bm{X}, lb_{\bm{X}}, \mathcal{C}_{\bm{X}})$ from $\mathcal{Q}$
	\State \Call{OptimalityRefutationTest}{$\bm{X}, lb_{\bm{X}}$}
	\Comment Algorithm \ref{alg:test-global-optimality}
	\State $lb_C \gets$ \Call{Contraction}{$\bm{X}$, $F$, $\mathcal{C}_{\bm{X}}$}
	\Comment Algorithm \ref{algo:bound-constrained-contractor} or \ref{algo:constrained-contractor}
	\If{$\bm{X}$ cannot be discarded}
		\State \Call{UpperBoundUpdateTest}{$m(\bm{X})$}
		\Comment Algorithm \ref{alg:upper-bound-update-test}
		\State Partition $\bm{X}$ into $\{\bm{X}_1, \bm{X}_2\}$
		\State $lb \gets \max(lb_{\bm{X}}, lb_C)$
		\Comment best lower bound
		\State Insert $(\bm{X}_1, lb, \mathcal{C}_{\bm{X}})$ and $(\bm{X}_2, lb, \mathcal{C}_{\bm{X}})$ into $\mathcal{Q}$
		\Comment Section \ref{sec:maxdist}
	\EndIf
\EndWhile
\State \Return $(\tilde{\bm{x}}, \tilde{f})$
\EndFunction
\end{algorithmic}
\end{algorithm}

For each newly extracted box, Charibde invokes a contraction procedure to reduce the domains of the variables or to discard the box ; the contractor for bound constrained problems is described in Section \ref{sec:bound constrained-contractor} and the contractor for nonlinearly constrained problems is described in Section \ref{sec:constrained-contractor}. If the constraints become inconsistent (that is, if the current box does not contain feasible points), the processing of the box is interrupted by an exception mechanism. Lower bounds of the objective function are computed using various techniques: natural extension (Definition \ref{def:natural-extension}), second-order extensions (Section \ref{sec:lower-bounds}), convexification (Section \ref{sec:convexification-contraction}) and monotonicity-based extension (Section \ref{sec:monotonicity-contraction}).


The priority queue $\mathcal{Q}$ is implemented using a binary heap, that is a complete binary tree in which the priority of a node is greater than the priority of its children. Inserting an element or extracting the element with the highest priority is carried out in logarithmic time.

We describe in Algorithm \ref{alg:test-global-optimality} the elementary refutation test that compares the best known upper bound $\tilde{f}$ of the global minimum and a lower bound $lb$ of the objective function on a given box $\bm{X}$. If no point in $\bm{X}$ can improve $\tilde{f}$ by at least the precision $\varepsilon$, then $\bm{X}$ can be safely discarded.

\begin{algorithm}[htbp!]
\caption{Refutation test based on global optimality}
\label{alg:test-global-optimality}
\begin{algorithmic}
\Function{OptimalityRefutationTest}{$\bm{X}$: box, $lb$: lower bound}
\If{$\tilde{f} - \varepsilon < lb$}
	\State Discard $\bm{X}$
\EndIf
\EndFunction
\end{algorithmic}
\end{algorithm}

Algorithm \ref{alg:upper-bound-update-test} illustrates the update of the best known upper bound $\tilde{f}$ of the global minimum when a punctual solution found by the \gls{IBC} improves $\tilde{f}$. If $\tilde{f}$ is indeed improved, the corresponding point is injected into the \gls{DE} population.

\begin{algorithm}[htbp!]
\caption{Test to update the best known upper bound of the global minimum}
\label{alg:upper-bound-update-test}
\begin{algorithmic}
\Function{UpperBoundUpdateTest}{$\bm{x}$: point}
\State $\overline{F_{\bm{x}}} \gets \overline{F(\bm{x})}$
\If{$\overline{F_{\bm{x}}} < \tilde{f}$ and $\bm{x}$ feasible}
	\State $(\tilde{\bm{x}}, \tilde{f}) \gets (\bm{x}, \overline{F_{\bm{x}}})$
	\Comment update of the best known upper bound
	\State \textcolor{amblu}{\texttt{MPI\_SendToDE}$(\tilde{\bm{x}}, \tilde{f})$}
\EndIf
\EndFunction
\end{algorithmic}
\end{algorithm}

For the sake of simplicity, the best known upper bound $\tilde{f}$, the corresponding point $\tilde{\bm{x}}$ and the precision $\varepsilon$ are handled as global variables in the algorithms.

\subsection{Second-order lower bounding}
\label{sec:lower-bounds}
Computing accurate lower bounds is crucial in \gls{BB} algorithms in order to prune suboptimal or infeasible subspaces. Naive interval extensions (such as the natural extension) usually provide a crude enclosure of the range of a function over a box due to dependency. Second-order extensions are a relatively cheap alternative to the natural extension for computing lower bounds. Bear in mind that, in the presence of dependency, second-order extensions become more precise than the natural extension as the width of the box approaches zero (the overapproximation tends to zero quadratically).

Algorithm \ref{algo:second-order} describes the computation of the Baumann extension by choosing the optimal Baumann center $\bm{c}_B^-$ ; it is the optimal mean value extension with respect to the lower bound (see Equation \ref{eq:centre-baumann}).
In our implementation, we use $\bm{c}_B^-$ as a trial point for the update of the best known upper bound $\tilde{f}$: if it is feasible, the value $\overline{F(\bm{c}_B^-)}$ is compared with the current value of $\tilde{f}$.

\begin{algorithm}[htbp!]
\caption{Second-order lower bounding}
\label{algo:second-order}
\begin{algorithmic}
\Function{SecondOrderLowerBounding}{$\bm{X}$: box, $F$: objective function, $\bm{G}$: gradient of $F$ on $\bm{X}$}
\State $(lb_B, \bm{c}_B^-) \gets$ \Call{BaumannLowerBound}{$\bm{X}, F, \bm{G}$}
\Comment Equation \ref{eq:centre-baumann}
\State \Call{OptimalityRefutationTest}{$\bm{X}, lb_B$}
\Comment Algorithm \ref{alg:test-global-optimality}
\State \Call{UpperBoundUpdateTest}{$\bm{c}_B^-$}
\Comment Algorithm \ref{alg:upper-bound-update-test}
\State \Return $lb_B$
\EndFunction
\end{algorithmic}
\end{algorithm}

\subsection{Contractor for bound constrained optimization}
\label{sec:bound constrained-contractor}

This section addresses a filtering procedure for bound constrained optimization problems:
\begin{equation}
\begin{aligned}
(\mathcal{P}) \quad & \min_{\bm{x} \in \bm{D} \subset \mathbb{R}^n} 	& f(\bm{x})
\end{aligned}
\end{equation}

Our implementation in Charibde is given in Algorithm \ref{algo:bound-constrained-contractor}.
\begin{algorithm}[h!]
\caption{Contractor for bound constrained optimization}
\label{algo:bound-constrained-contractor}
\begin{algorithmic}
\Function{Contraction}{\inout{} $\bm{X}$: box, $F$: objective function}
\State $F_{\bm{X}} \gets$ \Call{HC4Revise}{$F(\bm{X}) \le \tilde{f} - \varepsilon$}
\Comment contraction
\State $lb \gets \underline{F_{\bm{X}}}$
\Comment natural lower bounding
\State \Call{OptimalityRefutationTest}{$\bm{X}, lb$}
\Comment Algorithm \ref{alg:test-global-optimality}
\State $\bm{G} \gets \nabla F(\bm{X})$
\Comment gradient computation
\State $lb_B \gets$ \Call{SecondOrderLowerBounding}{$\bm{X}, F, \bm{G}$}
\Comment Algorithm \ref{algo:second-order}
\State \Call{StationarityContraction}{$\bm{X}, F, \bm{G}$}
\Comment Algorithm \ref{algo:stationarity-contraction}
\State \Return $\max(lb, lb_B)$
\EndFunction
\end{algorithmic}
\end{algorithm}

Algorithm \ref{algo:stationarity-contraction} describes a monotonicity-based subcontractor for bound constrained problems: when a component $X_i$ of the current box $\bm{X}$ is on the frontier of the domain $\bm{D}$ (that is, when $\underline{X_i} = \underline{D_i}$ or $\overline{X_i} = \overline{D_i}$), it is sometimes possible to reduce $X_i$ to one of its bound. Conversely, a box in the interior of the domain that does not contain a stationary point can be discarded.

A component $x_i$ of a local minimizer $\bm{x}$ is either stationary ($\frac{\partial f}{\partial x_i}(\bm{x}) = 0$) or activates a bound constraint ($x_i = \underline{D_i}$ or $x_i = \overline{D_i}$). This property can be exploited to devise a stationarity-based subcontractor that can be invoked as the revise procedure of a fixed-point algorithm. For a given box $\bm{X}$:
\begin{itemize}
\item if $\underline{X_i} = \underline{D_i}$ and
	\begin{itemize}
	\item $f$ is increasing with respect to $x_i$ on $\bm{X}$, $X_i$ can be contracted to its left bound ;
	\item $f$ is strictly decreasing with respect to $x_i$ on $\bm{X}$, $\bm{X}$ can be discarded ;
	\end{itemize}
\item if $\overline{X_i} = \overline{D_i}$ and
	\begin{itemize}
	\item $f$ is strictly increasing with respect to $x_i$ on $\bm{X}$, $\bm{X}$ can be discarded ;
	\item $f$ is decreasing with respect to $x_i$ on $\bm{X}$, $X_i$ can be contracted to its right bound.
	\end{itemize}
\item if $X_i$ is interior to the domain ($X_i \subset int(D_i)$), the minimum of $f$ on $\bm{X}$ corresponds to a point $\bm{x}$ where $\frac{\partial f}{\partial x_i}(\bm{x}) = 0$. Consequently, if $\frac{\partial F}{\partial x_i}(\bm{X})$ does not contain zero, $\bm{X}$ can be discarded.
\end{itemize}

\begin{algorithm}[htbp!]
\caption{Stationarity-based contractor for bound constrained problems}
\label{algo:stationarity-contraction}
\begin{algorithmic}
\Function{StationarityContraction}{\inout{} $\bm{X}$: box, $F$: objective function, $\bm{G}$: gradient of $F$ on $\bm{X}$}
\For{$i \in \{1, \ldots, n\}$}
	\If{$X_i$ on the frontier of $\bm{D}$}
		\If{$G_i$ has constant sign}
			\State reduce $X_i$ to one of its bounds or discard $\bm{X}$
		\EndIf
	\ElsIf{$0 \notin G_i$}
	\Comment interior component that contains no stationary point
		\State discard $\bm{X}$
	\Else
		\State \Call{HC4Revise}{$\frac{\partial F}{\partial x_i}(\bm{X}) = 0$}
		\Comment discard non stationary values
	\EndIf
\EndFor
\EndFunction
\end{algorithmic}
\end{algorithm}

When $0 \in \frac{\partial F}{\partial x_i}(\bm{X})$, a cheap 2B contractor is invoked on the stationarity condition $\frac{\partial F}{\partial x_i}(\bm{X}) = 0$ in order to discard non-stationary points.
Note that filtering with respect to partial derivatives seems to be an uncommon practice in the global optimization community, although it may provide a substantial contraction of the search space.


\subsection{Contractor for nonlinearly constrained optimization}
\label{sec:constrained-contractor}

This section addresses a filtering procedure for nonlinearly constrained optimization problems:
\begin{equation}
\begin{aligned}
(\mathcal{P}) \quad & \min_{\bm{x} \in \mathbb{R}^n} 	& f(\bm{x}) & \\
					& s.t.	& g_j(\bm{x}) \le 0, & \quad j \in \{1, \ldots, m\} \\
					& 		& h_j(\bm{x}) = 0, & \quad j \in \{1, \ldots, p\} \\
\end{aligned}
\end{equation}

Interval-based solvers may address equality constraints $h_j(\bm{x}) = 0$ $(j \in \{1, \ldots, p\})$ in two different ways:
\begin{itemize}
\item GlobSol \cite{Kearfott1996Rigorous} and Icos \cite{Lebbah2005} produce a small box $\bm{X}$ guaranteed to contain a point $\bm{x}$ that minimizes $f$ and satisfies the constraints:
\begin{equation}
\begin{cases}
g_j(\bm{x}) \le 0, & \quad j \in \{1, \ldots, m\} \\
h_j(\bm{x}) = 0, & \quad j \in \{1, \ldots, p\}
\end{cases}
\end{equation}
The existence of $\bm{x}$ in $\bm{X}$ is numerically proven by a multivariate interval Newton method ;

\item IBBA \cite{Ninin2010Reliable}, Ibex \cite{Trombettoni2011Inner} and Charibde handle a relaxed problem in which each equality constraint $h_j(\bm{x}) = 0$ $(j \in \{1, \ldots, p\})$ is replaced with a pair of inequality constraints:
\begin{equation}
{-}\varepsilon_= \le h_j(\bm{x}) \le \varepsilon_=
\end{equation}
where $\varepsilon_=$ can be chosen arbitrarily small.
\end{itemize}

Our implementation in Charibde is given in Algorithm \ref{algo:constrained-contractor}. Various advanced lower bounding strategies, enabled by default in Charibde and described in the following subsections, can be disabled by the user.

\begin{algorithm}[htbp!]
\caption{Contractor for constrained optimization}
\label{algo:constrained-contractor}
\begin{algorithmic}
\Function{Contraction}{\inout{} $\bm{X}$: box, $F$: objective function, \inout{} $\mathcal{C}$: set of constraints}
\State $lb \gets -\infty$
\Comment lower bound
\Repeat
	\State $\bm{X}' \gets \bm{X}$
	\State $F_{\bm{X}} \gets$ \Call{HC4Revise}{$F(\bm{X}) \le \tilde{f} - \varepsilon$}
	\State $lb \gets \underline{F_{\bm{X}}}$
	\Comment natural lower bounding
	\State \Call{OptimalityRefutationTest}{$\bm{X}, lb$}
	\Comment Algorithm \ref{alg:test-global-optimality}
	\State $\bm{G} \gets \nabla F(\bm{X})$
	\Comment gradient computation
	\State $lb_B \gets$ \Call{SecondOrderLowerBounding}{$\bm{X}, F, \bm{G}$}
	\Comment Algorithm \ref{algo:second-order}
	\State $lb_M \gets$ \Call{MonotonicityBasedContraction}{$\bm{X}, F, \mathcal{C}$}
	\Comment Algorithm \ref{algo:monotonicity-based-contractor}
	\State $\mathcal{C} \gets$ \Call{HC4}{$\bm{X}, \mathcal{C}$} or \Call{Mohc}{$\bm{X}, \mathcal{C}$}
	\State $lb_{lp} \gets$ \Call{Convexification}{$\bm{X}, F, \bm{G}, \mathcal{C}$}
	\Comment Algorithm \ref{algo:convexification}
\Until{$\bm{X} = \varnothing$ or $gain(\bm{X}, \bm{X'}) < \eta$}
\State \Return $\max(lb, lb_B, lb_M, lb_{lp})$
\Comment best lower bound
\EndFunction
\end{algorithmic}
\end{algorithm}

\subsubsection{Contraction based on convexification}
\label{sec:convexification-contraction}
Similarly to X-Newton~\cite{Araya2012Contractor}, we implement a subcontractor based on a convexification of the problem: the objective function and the constraints are approximated by convex linear relaxations (Taylor extensions in which the expansion point is a corner of the box, see Section \ref{sec:convexification}), which results in a linear problem.
By default in Charibde, Dantzig's simplex algorithm is invoked once to compute a lower bound of the original problem ; $2n$ additional simplex calls may also contract the domains of the variables on demand (Algorithm \ref{algo:convexification}).

\begin{algorithm}[htbp!]
\caption{Lower bounding and contraction using convexification}
\label{algo:convexification}
\begin{algorithmic}
\Function{Convexification}{\inout{} $\bm{X}$: box, $F$: objective function, $\bm{G}$: gradient of $F$ on $\bm{X}$, $\mathcal{C}$: set of constraints}
\State $polytope \gets$ \Call{CornerTaylor}{$\bm{X}, \mathcal{C}$}
\Comment linear convex relaxation of the constraints
\State $lb_{lp} \gets$ \Call{Simplex}{$\bm{X}, G, polytope$}
\Comment lower bounding
\State \Call{OptimalityRefutationTest}{$\bm{X}, lb_{lp}$}
\Comment Algorithm \ref{alg:test-global-optimality}
\State \Call{X-Newton}{$\bm{X}, polytope$}
\Comment contraction of $2n$ variables
\State \Return $lb_{lp}$
\EndFunction
\end{algorithmic}
\end{algorithm}

Charibde invokes the binding ocaml-glpk~\cite{OCamlGlpk2004} to GLPK (GNU Linear Programming Kit), a library for solving large linear problems and mixed problems in \gls{FPA}. We implement a cheap postprocessing step~\cite{Neumaier2004Safe} that determines a rigorous lower bound using \gls{IA} based on the optimal solution of the linear problem.

\cite{Araya2012Contractor} suggest to compute the linear lower bound $lb_{lp}$ of the original problem by solving the linear problem from scratch, then to solve the additional $2n$ calls by taking the solution of the previous call as an initial feasible basis. Unfortunately, ocaml-glpk does not implement this feature. The contraction of variables in Charibde's X-Newton implementation thus comes with a substantial cost. However, we show in Section \ref{sec:coconut-comparison} that most numerical results are achieved without the $2n$ calls of X-Newton.

\subsubsection{Contraction based on monotonicity} 
\label{sec:monotonicity-contraction}

Charibde's main contractor for constrained problems (Algorithm \ref{algo:constrained-contractor}) invokes Mohc on the set of constraints $\mathcal{C}$. An additional call to MohcRevise is performed to contract the box $\bm{X}$ with respect to the inequality constraint $f(\bm{X}) \le \tilde{f} - \varepsilon$ if the expression of $f$ contains several occurrences of the variables. In this case, only a call to MinRevise and a call to Left/RightNarrowFmin are required.

Charibde implements a specific procedure for global optimization (Algorithm \ref{algo:monotonicity-based-contractor}): it may either contract $\bm{X}$ into a feasible subbox without losing the minimum of $f$ over $\bm{X}$, or extract an upper bound of the global minimum.
Remember that $\bm{X}^-$ (see Definition \ref{def:monotonicity-extension}) is a subbox of $\bm{X}$ obtained by replacing variables with multiples occurrences and with respect to which $f$ is monotonic with one of their bounds (the left bound if $f$ is increasing, the right bound if $f$ is decreasing). By construction, $\bm{X}^-$ contains the unconstrained minimizer $\bm{x}_{\bm{X}}^*$ of $f$ over $\bm{X}$, that is the solution to:
\begin{equation}
\min_{\bm{x} \in \bm{X}} f(\bm{x})
\end{equation}
If $\bm{X}^-$ is feasible, $\bm{x}_{\bm{X}}^*$ is also the constrained minimizer of $f$ over $\bm{X}$ ; in this case, $\bm{X}$ can therefore be contracted to $\bm{X}^-$. Otherwise, we attempt to update the best known upper bound $\tilde{f}$ by picking a point within $\bm{X}^-$ (for example, its midpoint).

\begin{algorithm}[htbp!]
\caption{Contraction and lower bounding using monotonicity}
\label{algo:monotonicity-based-contractor}
\begin{algorithmic}
\Function{MonotonicityBasedContraction}{\inout{} $\bm{X}$: box, $F$: objective function, \inout{} $\mathcal{C}$: set of constraints}
\State $(lb_M, \bm{X}^-) \gets$ \Call{MohcRevise}{$F(\bm{X}) \le \tilde{f} - \varepsilon$}
\Comment monotonicity-based contraction
\State \Call{OptimalityRefutationTest}{$\bm{X}, lb_M$}
\Comment Algorithm \ref{alg:test-global-optimality}
\If{$\bm{X}^-$ is feasible}
	\State $\bm{X} \gets \bm{X}^-$
	\Comment Definition \ref{def:monotonicity-extension}
	\State $\mathcal{C} \gets \varnothing$
\Else
	\State \Call{UpperBoundUpdateTest}{$m(\bm{X}^-)$}
	\Comment Algorithm \ref{alg:upper-bound-update-test}
\EndIf
\State \Return $lb_M$
\EndFunction
\end{algorithmic}
\end{algorithm}

\subsection{Contraction and automatic differentiation}
\label{sec:contraction-ad}

The similarities between the double traversal of the syntax tree of the HC4Revise algorithm (Section \ref{sec:hc4revise}) and the \gls{AD} in adjoint mode (Section \ref{sec:ad-adjoint}) were mentioned by~\cite{Schichl2005Interval}. Since the bottom-up evaluation phase is shared by both algorithms, it can be carried out only once. \cite{Schichl2005Interval} even suggest to carry out the top-down phase of the \gls{AD} \textit{after} the top-down propagation phase of HC4Revise in order to exploit the contracted \textit{intermediary} nodes. The resulting \textit{constrained derivatives} (the infeasible values that have been discarded do not contribute to the computation of the partial derivatives) are tighter than standard derivatives and may be used in derivative-based refutation techniques (Taylor form, Mohc).

Example \ref{ex:hc4-ad} illustrates the benefits of this approach on the constraint $x+(x+y)^2 - 1 = 0$. Computing the derivatives in the standard way on the box $X \times Y = [0, 1] \times [0, 1]$ yields $\frac{\partial G}{\partial x}(X, Y) = [1, 5]$ and $\frac{\partial G}{\partial y}(X, Y) = [0, 4]$. If however the top-down phase of HC4Revise is exploited, the constrained derivatives are $\frac{\partial G}{\partial x}(X, Y) = [1, 3]$ and $\frac{\partial G}{\partial y}(X, Y) = [0, 2]$. Here, the intermediary node $x + y \in [0, 1]$ obtained after the top-down phase of HC4Revise can be exploited to compute a tighter enclosure of the derivatives.

\begin{example}
\label{ex:hc4-ad}
Let:
\begin{itemize}
\item $g(x, y) = x+(x+y)^2 - 1 = 0$ be an equality constraint ;
\item $X = [0, 5]$ and $Y = [0, 5]$.
\end{itemize}

The elementary operations (the intermediary nodes) are:
\begin{equation}
\begin{aligned}
n_1 & \eqdef x + y 	& \quad n_3 & \eqdef n_2 + x \\
n_2 & \eqdef n_1^2 	& \quad n_4 & = n_3 - 1 \\
\end{aligned}
\end{equation}

The interval evaluation during the bottom-up phase yields:
\begin{equation}
\begin{aligned}
N_1 & = X + Y \in [0, 10] \\
N_2 & = N_1^2 \in [0, 100] \\
N_3 & = N_2 + X \in [0, 105] \\
N_4 & = N_3 - 1 \in [-1, 104]
\end{aligned}
\end{equation}

We now intersect $N_4$ with $[0, 0]$ and compute the top-down phase:
\begin{equation}
\begin{aligned}
N_3' & = N_3 \cap (1 + N_4') = [0, 105] \cap [1, 1] = [1, 1] \\
N_2' & = N_2 \cap (N_3' - X) = [0, 100] \cap ([1, 1] - [0, 5]) = [0, 100] \cap [-4, 1] = [0, 1] \\
N_1' & = N_1 \cap \square \left(-\sqrt{N_2'}, \sqrt{N_2'} \right) = [0, 10] \cap [-1, 1] = [0, 1] \\
X' & = X \cap (N_3' - N_2') = [0, 5] \cap (1 - [0, 1]) = [0, 5] \cap [0, 1] = [0, 1] \\
X' & = X' \cap (N_1' - Y) = [0, 1] \cap ([0, 1] - [0, 5]) = [0, 1] \cap [-5, 1] = [0, 1] \\
Y' & = Y \cap (N_1' - X') = [0, 5] \cap ([0, 1] - [0, 1]) = [0, 5] \cap [-1, 1] = [0, 1]
\end{aligned}
\end{equation}

The domains of $x$ and $y$ have both been contracted to $X = Y = [0, 1]$. The node $n_1 = x + y$ has been contracted to $N_1' = [0, 1]$, while the direct evaluation $X' + Y'$ yields $[0, 2] \supset N_1'$.
The relation $x+y \in [0, 1]$ over the box $[0, 1] \times [0, 1]$ is represented in Figure \ref{fig:implicit}. The exact domain (colored) is a polytope whose edges are not parallel to the axes, and cannot be represented exactly by a box.
\begin{figure}[htbp]
\centering
\def\svgwidth{0.7\columnwidth}
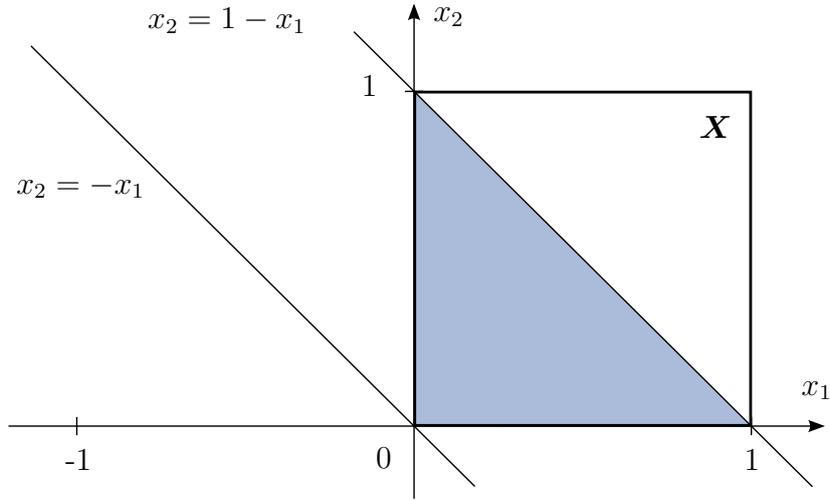
\caption{Polytope of the relation $x+y \in [0, 1]$}
\label{fig:implicit}
\end{figure}

The approach of~\cite{Schichl2005Interval} exploits the information of the polytope $\{x + y \in [0, 1], x \in [0, 1], y \in [0, 1]\}$
and computes the constrained derivatives:
\begin{align}
\frac{\partial G}{\partial x}(X, Y) & = 1 + 2(X + Y) = 1 + 2\bm{N_1'} = 1 + 2[0, 1] = [1, 3] \subset [1, 5]\\
\frac{\partial G}{\partial y}(X, Y) & = 2(X + Y) = 2\bm{N_1'} = 2[0, 1] = [0, 2] \subset [0, 4]
\end{align}

\end{example}

\subsection{MaxDist: a geometrical search strategy}
\label{sec:maxdist}
The order in which the remaining boxes are inserted into the priority queue $\mathcal{Q}$ determines the exploration strategy of the search space (see Section \ref{sec:search-heuristics}). Numerical tests suggest that the "best-first search" strategy is not consistent when the objective function is highly subject to dependency: the lower bound of $f$ over a large box usually provides little information on the actual range of $f$. The "largest first search" strategy performs a "breadth-first" search that does not lay emphasis on promising subspaces. The "depth-first search" strategy tends to quickly explore the neighborhood of local minima, however it struggles to escape from them efficiently. We therefore investigated a new strategy specific to the hybrid algorithm between the \gls{IBC} and the \gls{DE}.

We observed that the local convergence of the \gls{DE} in the neighborhood of a local minimum is excellent. However, when the population of the \gls{DE} agglutinates and gets stuck at a local minimum, it may be difficult to avoid premature convergence and to explore other areas of the search space.
We therefore drew the following conclusions about our hybridization:
\begin{itemize}
\item if the global minimizer $\bm{x}^*$ is located close to the best known solution $\tilde{\bm{x}}$, we hope that the \gls{DE} will quickly converge towards $\bm{x}^*$ ;
\item if $\bm{x}^*$ is distant from $\tilde{\bm{x}}$, the \gls{IBC} should not explore further the neighborhood of $\tilde{\bm{x}}$, and should focus on areas of the search space that are possibly out of reach for the \gls{DE}.
\end{itemize}

We now introduce a new search strategy called MaxDist. The idea is to process and hopefully discard the neighborhood of the global minimizer $\bm{x}^*$ (a priori not known) only when the best possible upper bound $\tilde{f}$ of the global minimum is available. This is usually an arduous task, on account of the similarity of the values $f(\bm{x})$ when $\bm{x}$ lives in the neighborhood of $\bm{x}^*$. In order to delay the costly processing of this neighborhood, we propose to extract from $\mathcal{Q}$ the box that is the farthest in the search space from the current solution $\tilde{\bm{x}}$. MaxDist thus explores (hopefully promising) areas of the search space that are hardly accessible to the \gls{DE}.

The notion of distance between a point $\bm{x}$ and a box $\bm{X}$ is detailed in Algorithm \ref{algo:distance}. Note that MaxDist is an adaptive strategy: whenever the best known solution $\tilde{\bm{x}}$ is updated, $\mathcal{Q}$ is reordered according to the new priorities of the boxes.

\begin{algorithm}[h!]
\caption{Distance between a point and a box}
\label{algo:distance}
\begin{algorithmic}
\Function{Distance}{$\bm{x}$: point, $\bm{X}$: box}
\State $d \gets 0$
\For{$j = 1$ to $n$}
	\If{$\overline{X_j} < x_j$}
		\State $d \gets d + (x_j - \overline{X_j})^2$
	\ElsIf{$x_j < \underline{X_j}$}
		\State $d \gets d + (\underline{X_j} - x_j)^2$
	\EndIf
\EndFor
\State \Return $d$
\EndFunction
\end{algorithmic}
\end{algorithm}

In Table \ref{tab:comparison-strategies}, we compare our strategy MaxDist against the standard strategies "best-first search" and "largest first" on eight optimization problems from the COCONUT benchmark\footnote{The COCONUT benchmark is available at \url{http://www.mat.univie.ac.at/~neum/glopt/coconut/Benchmark/Benchmark.html}}.
For each problem are given the CPU time (in seconds) and the maximal size $|\mathcal{Q}|_{max}$ of $\mathcal{Q}$.

\begin{table}[htbp!]
	\centering
	\caption{Comparison of exploration strategies}
	\begin{tabular}{|l|cc|cc|cc|}
	\hline
	& \multicolumn{2}{c|}{best-first search} & \multicolumn{2}{c|}{largest first} & \multicolumn{2}{c|}{MaxDist} \\
	Problem 	& CPU time (s) & $|\mathcal{Q}|_{max}$ & CPU time (s) & $|\mathcal{Q}|_{max}$ & CPU time (s) & $|\mathcal{Q}|_{max}$ \\
	\hline
	avgasa		& 6.26 	& 4,301 	& 6.12 	& 39,864 	& \bf{5.7} 	& 23 \\
	ex2\_1\_7	& 47.9 	& 29,530	& 28.3	& 78,784	& \bf{26} 	& 19 \\
	ex2\_1\_9	& 44	& 111,158	& 43.6 	& 54,587	& \bf{37.5}	& 134 \\
	ex7\_3\_5	& 12	& 34,430	& 8.9	& 14,064	& \bf{8.42}	& 87 \\
	ex6\_2\_6	& 2.13	& 3,800		& 2.15	& 2,236		& \bf{1.96}	& 45 \\
	ex6\_2\_8	& 3.2	& 6,377		& 3.13 	& 4,316		& \bf{3.03}	& 20 \\
	ex6\_2\_9	& 3.7	& 5,428		& 3.66	& 2,924		& \bf{3.47}	& 28 \\
	ex6\_2\_11	& 2.47	& 4,928 	& 2.41	& 2,556		& \bf{2.28}	& 38 \\
	\hline
	Sum			& 121.66&		& 98.27	&		& \bf{88.36}& \\
	\hline
	\end{tabular}
	\label{tab:comparison-strategies}
\end{table}

These preliminary results suggest that MaxDist is competitive with the standard strategies: the eight test problems are solved with MaxDist in 88.36s, that is 10.1\% faster than "largest first" (98.27s) and 27.4\% faster than "best-first search" (121.66s).
MaxDist also exhibits a remarkably low size of the priority queue (between 19 and 134 boxes at most) compared to "best-first search" (between 3,800 and 111,158 boxes) and "largest first" (between 2,236 and 78,784 boxes). Although this has a limited impact on the insertion time into the priority queue (a binary heap has a logarithmic worst-case time complexity), it may benefit the hybridization between the \gls{IBC} and the \gls{DE}: the priority queue can be sent to the \gls{DE} at low cost in order to update its search space. This technique is described in Section \ref{sec:interval-domain}.

\section{Differential evolution algorithm}
\label{sec:charibde-de}

The \gls{DE} algorithm is a credible candidate for computing an accurate upper bound of the global minimum, thus accelerating the pruning of the search space of the \gls{IBC} ; it proved very competitive on continuous optimization problems and has few hyperparameters, which makes the tuning of the algorithm less tedious.
Whenever the \gls{DE} improves the best known solution, the corresponding individual and its evaluation are sent to the \gls{IBC} (Algorithm \ref{alg:charibde-de}): the evaluation updates the best known upper bound $\tilde{f}$ of the global minimum, and the individual contributes to the computation of the MaxDist exploration strategy (see Section \ref{sec:maxdist}).
In return, whenever the \gls{IBC} evaluates a feasible point that improves the best known solution, the point is injected into the \gls{DE} population. In order to avoid replacing the whole population through successive injections of \gls{IBC} points, the same individual is systematically replaced.

\begin{algorithm}[h!]
\caption{Charibde: differential evolution algorithm}
\label{alg:charibde-de}
\begin{algorithmic}[]
\Function{DifferentialEvolution}{$f$: objective function, $\mathit{NP}$: population size, $W$: scaling factor, $\mathit{CR}$: crossover rate, $\bm{D}$: domain}
\State $P \gets$ initial population randomly generated in $\bm{D}$
\textcolor{amred}{\State $\tilde{f}_{nearest} \gets +\infty$}
\textcolor{amred}{\State $\tilde{f} \gets +\infty$}
\Repeat
	\State \textcolor{amblu}{$(\bm{x}, f_{\bm{x}}) \gets$ \texttt{MPI\_ReceiveFromIBC}()}
	\State \textcolor{amblu}{add $\bm{x}$ to $P$}
	\State \textcolor{amblu}{$\tilde{f} \gets f_{\bm{x}}$}
	\State generate temporary population $P'$ based on $P$
	\Comment Algorithm \ref{alg:de}
	\State $P \gets P'$
	\State $\bm{x}_{best} \gets$ \Call{BestIndividual}{$P$}
	\Comment Algorithm \ref{algo:de-reliable-evaluation}
	\textcolor{amred}{
	\State \Call{ReliableEvaluation}{$\bm{x}_{best}$, $\tilde{f}_{nearest}$, $\tilde{f}$}
	}%
\Until{termination criteria met}
\State \Return best individual of $P$
\EndFunction
\end{algorithmic}
\end{algorithm}

The cooperation between the \gls{DE} and the \gls{IBC} is of critical importance: since the \gls{IBC} is intrinsically reliable and guarantees the global optimality of the solution even in the presence of roundoff errors, the exchanges between the \gls{IBC} and the \gls{DE} must be equally reliable.
The robust handling of the objective function (respectively the constraints) is detailed in Section~\ref{sec:rigorous-objective} (respectively Section \ref{sec:rigorous-constraints}).


\subsection{Rigorous objective handling}
\label{sec:rigorous-objective}
Let $f_{nearest}$ be the machine implementation of $f$ using the round-to-nearest rounding method. For any $\bm{x}$ of the domain of $f$, $f_{nearest}(\bm{x}) \in F(\bm{x})$ holds. When $\bm{x}$ is feasible with respect to the constraints of the problem, the bound $\overline{F(\bm{x})}$ is a reliable upper bound of the global minimum $f^*$, while nothing can be said about $f_{nearest}(\bm{x})$. Two strategies may thus be considered:
\begin{itemize}
\item the objective function is systematically evaluated using \gls{IA}. This amounts to minimizing the function $\bm{x} \mapsto \overline{F(\bm{x})}$ ;
\item the objective function is evaluated using \gls{IA} only when the best known evaluation is improved.
\end{itemize}

\begin{figure}[htbp!]
\centering
\def\svgwidth{0.65\columnwidth}
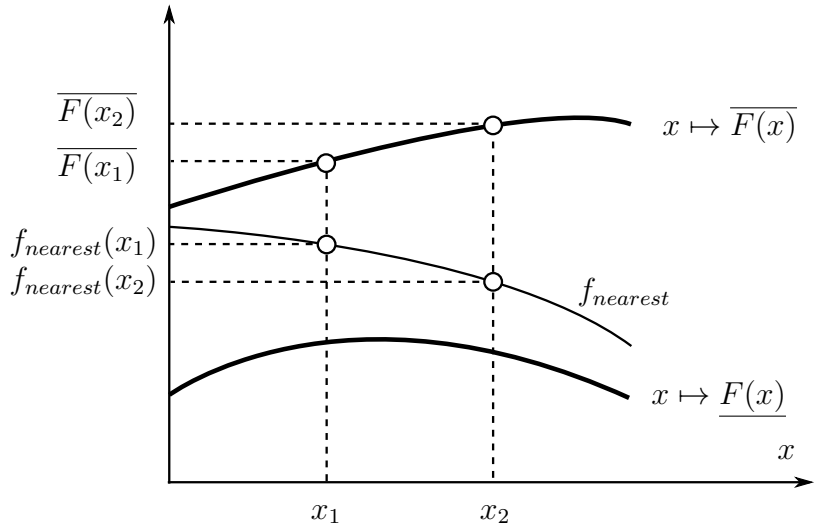
\caption{Opposite monotonicities of $f_{nearest}$ and $x \mapsto \overline{F(x)}$}
\label{fig:precision}
\end{figure}

The first strategy is the more rigorous, since both the \gls{DE} and the \gls{IBC} handle the same objective function $\bm{x} \mapsto \overline{F(\bm{x})}$ based on \gls{IA}. However, it induces a higher evaluation cost than an evaluation method based on \gls{FPA}. The second strategy assumes that $f_{nearest}$ and $\bm{x} \mapsto \overline{F(\bm{x})}$ have the same monotonicity.
Numerical tests on some test problems show that, for given $\bm{x}_1$ and $\bm{x}_2$, $f_{nearest}(\bm{x}_1) >
f_{nearest}(\bm{x}_2)$ and $\overline{F(\bm{x}_1)} < \overline{F(\bm{x}_2)}$ hold (Figure \ref{fig:precision}). This proves that the assumption on the monotonicity of $f_{nearest}$ and $\bm{x} \mapsto \overline{F(\bm{x})}$ does not hold. However, the gap between the round-to-nearest evaluation and the \gls{IA} evaluation is generally much smaller (around $10^{-15}$) than the user-defined tolerance $\varepsilon$.

This observation motivated the implementation of the second strategy in Charibde. The objective values of the \gls{DE} individuals are systematically evaluated using \gls{FPA} ; whenever the best round-to-nearest evaluation is improved, the objective function is rigorously bounded using \gls{IA} and the right bound of the result if compared with the best known reliable upper bound. If it is improved, the new reliable upper bound is sent to the \gls{IBC} (Algorithm \ref{algo:de-reliable-evaluation}). This choice substantially reduced the evaluation cost of the \gls{DE}, while guaranteeing that all the communications between the \gls{DE} and the \gls{IBC} are reliable.

\begin{algorithm}[htbp!]
\caption{Reliable evaluation of the best individual}
\label{algo:de-reliable-evaluation}
\begin{algorithmic}
\Function{ReliableEvaluation}{$\bm{x}$: best individual, \inout{} $\tilde{f}_{nearest}$: best known evaluation using floating-point arithmetic, \inout{} $\tilde{f}$: best known reliable evaluation using interval arithmetic}
\If{$\bm{x}$ is feasible}
	\State{$f_{\bm{x}} \gets f_{nearest}(\bm{x})$}
	\Comment round-to-nearest evaluation
	\If{$f_{\bm{x}} < \tilde{f}_{nearest}$}
		\State $\tilde{f}_{nearest} \gets f_{\bm{x}}$
		\Comment update of the best known evaluation
		\If{$\overline{F(\bm{x})} < \tilde{f}$}
			\State $\tilde{f} \gets \overline{F(\bm{x})}$
			\Comment update of the best known reliable evaluation
			\State \texttt{MPI\_SendToIBC}$(\bm{x}, \tilde{f})$
		\EndIf
	\EndIf
\EndIf
\EndFunction
\end{algorithmic}
\end{algorithm}

\subsection{Constraint handling}
\label{sec:rigorous-constraints}
The extension of \glspl{EA} to constrained optimization was the subject of numerous works. The most common approaches are penalty methods and direct constraint handling (see Section \ref{sec:constraint-handling}). The latter approach computes the violation of the constraints (the number of violated constraints and their magnitudes) ; we have adapted this strategy within the comparison operator between an individual and its parent.
For an individual $\bm{x}$:
\begin{itemize}
\item $f_{\bm{x}}$ is the objective value of $\bm{x}$ ;
\item $n_{\bm{x}}$ is the number of constraints violated by $\bm{x}$ ;
\item $\displaystyle s_{\bm{x}} \eqdef \sum_{j=1}^m \max(0, g_j(\bm{x}))$ is the constraint violation (the sum of magnitudes of the constraints violated by $\bm{x}$).
\end{itemize}
The evaluation of an individual $\bm{x}$ can be written as an enumerative type:
\begin{itemize}
\item $\mathit{Infeasible}(n_{\bm{x}}, s_{\bm{x}})$ when $\bm{x}$ is infeasible. The objective function is not evaluated ;
\item $\mathit{Feasible}(f_{\bm{x}})$ when $\bm{x}$ is a feasible individual.
\end{itemize}

\subsection{Selection by comparison}
\label{sec:evaluation-by-comparison}
\Gls{DE} algorithms require a comparison operator (a binary relation) that determines whether an individual $\bm{y}$ is better than its parent $\bm{x}$, and should replace it in the population. We propose the following rules:
\begin{enumerate}
\item a feasible individual is better than an infeasible individual ;
\item among two feasible individuals, the better one is the one with the lowest objective value ;
\item among two infeasible individuals, the better one is the one with the lowest number of violated constraints, then (in case of a tie) the one with the lowest constraint violation.
\end{enumerate}
The function \textsc{SelectionByComparison} (Algorithm \ref{algo:selection-by-comparison}) implements the evaluation of a newly generated individual $\bm{y}$ by comparison with its parent $\bm{x}$. The underlying strategy is to exploit the known evaluation of the parent in a refutation test to limit the number of evaluations of the constraints for $\bm{y}$.

\begin{algorithm}[h!]
\caption{Comparison operator between an individual and its parent}
\label{algo:selection-by-comparison}
\begin{algorithmic}
\Function{SelectionByComparison}{$\bm{x}$: parent individual, $\bm{y}$: individual, $f$: objective function, $\mathcal{C}$: set of constraints}
\If{$\bm{x}$ has evaluation $\mathit{Infeasible}(n_{\bm{x}}, s_{\bm{x}})$}
	\State $(\mathit{feasible}_{\bm{y}}, n_{\bm{y}}, s_{\bm{y}}) \gets$ \Call{RigorousFeasibilityMeasure}{$\bm{y}$, $\mathcal{C}$}
	\Comment Algorithm \ref{algo:rigorous-constraint-evaluation}
	\If{$\mathit{feasible}_{\bm{y}}$}
		\State \Return $(\bm{y}, \mathit{Feasible}(f(\bm{y})))$
		\Comment evaluation of the objective function
	\ElsIf{$(n_{\bm{y}}, s_{\bm{y}})$ better than $(n_{\bm{x}}, s_{\bm{x}})$} 
		\State \Return $(\bm{y}, \mathit{Infeasible}(n_{\bm{y}}, s_{\bm{y}}))$
	\Else
		\State \Return $(\bm{x}, \mathit{Infeasible}(n_{\bm{x}}, s_{\bm{x}}))$
	\EndIf
\Else
	\Comment $\bm{x}$ has evaluation $\mathit{Feasible}(f_{\bm{x}})$
	\State $f_{\bm{y}} \gets f(\bm{y})$
	\Comment evaluation of the objective function
	\If{$f_{\bm{y}} < f_{\bm{x}}$}
		\State $\mathit{feasible}_{\bm{y}}$ $\gets$ \Call{RigorousFeasibilityTest}{$\bm{y}$,
$\mathcal{C}$}
		\Comment Algorithm \ref{algo:is-feasible} or \ref{algo:is-feasible-floating}
		\If{$\mathit{feasible}_{\bm{y}}$}
			\State \Return $(\bm{y}, \mathit{Feasible}(f_{\bm{y}}))$
			
		\Else
			\State \Return $(\bm{x}, \mathit{Feasible}(f_{\bm{x}}))$
		\EndIf
	\Else
	\Comment $\bm{y}$ cannot improve the objective function
		\State \Return $(\bm{x}, \mathit{Feasible}(f_{\bm{x}}))$
	\EndIf
\EndIf
\EndFunction
\end{algorithmic}
\end{algorithm}

\subsection{Rigorous feasibility}
\label{sec:rigorous-feasibility}

Numerous nonlinear optimization solvers tolerate a slight numerical violation of the constraints (for example, $g(\bm{x}) \le 10^{-6}$ instead of $g(\bm{x}) \le 0$). Such a "pseudo-feasible" point $\bm{x}$ (that satisfies these relaxed constraints) brings no reliable information because of numerical errors ; in particular, its objective value cannot be deemed as a reliable upper bound of the global minimum. In practice, its objective value can be lower than the global minimum, or $\bm{x}$ can be rather distant from feasible solutions in the search space.

We opted for constraint evaluation using \gls{IA} ; an individual $\bm{x}$ is labeled as "feasible" only when it is numerically guaranteed that is satisfies the constraints of the problem $g_j \le 0$ ($i \in \{1, \ldots, m\})$:
\begin{equation}
\forall i \in \{1, \ldots, m\}, \quad \overline{G_j(\bm{x})} \le 0
\label{eq:individual-feasible}
\end{equation}

The function \textsc{RigorousFeasibilityMeasure} (Algorithm \ref{algo:rigorous-constraint-evaluation}) implements a feasibility evaluation: it evaluates all constraints from $\mathcal{C}$ using \gls{IA} and measures the number of violated constraints and their magnitudes.

\begin{algorithm}[h!]
\caption{Rigorous measure of feasibility}
\label{algo:rigorous-constraint-evaluation}
\begin{algorithmic}
\Function{RigorousFeasibilityMeasure}{$\bm{x}$: individual, $\mathcal{C}$: set of constraints}
\State $(\mathit{feasible}, violated, violation) \gets (\mathit{true}, 0, 0.)$
\For{$g_j \in \mathcal{C}$}
	\If{$0 < \overline{G_j(\bm{x})}$}
	\Comment interval constraint not satisfied
		\State{$(\mathit{feasible}, violated, violation) \gets (\mathit{false}, violated+1, violation + \overline{G_j(\bm{x})})$}
	\EndIf
\EndFor
\State \Return $(\mathit{feasible}, violated, violation)$
\EndFunction
\end{algorithmic}
\end{algorithm}

The function \textsc{RigorousFeasibilityTest} (Algorithm \ref{algo:is-feasible}) is a mere feasibility test: it verifies whether a point is feasible and stops at the first violated constraint.

\begin{algorithm}[htbp!]
\caption{Rigorous feasibility test}
\label{algo:is-feasible}
\begin{algorithmic}
\Function{RigorousFeasibilityTest}{$\bm{x}$: individual, $\mathcal{C}$: set of constraints}
\For{$g_j \in \mathcal{C}$}
	\If{$0 < \overline{G_j(\bm{x})}$}
	\Comment violated interval constraint
		\State \Return $\mathit{false}$
	\EndIf
\EndFor
\State \Return $\mathit{true}$
\EndFunction
\end{algorithmic}
\end{algorithm}

\subsection{Refutation using floating-point arithmetic}

\cite{Neumaier1990Interval} estimated that the evaluation of a function using \gls{IA} is 2 to 4 times costlier than using \gls{FPA}. Albeit not rigorous, Charibde exploits the refutation potential and the low cost of floating-point evaluations whenever possible. Since the round-to-nearest evaluation $g_{nearest}$ is enclosed by the \gls{IA} evaluation $G$, the following relation:
\begin{equation}
g_{nearest}(\bm{x}) \le \overline{G(\bm{x})}
\end{equation}
holds for any $\bm{x}$. The following refutation test was implemented in the feasibility test of the \gls{DE}:
\begin{equation}
\label{eq:refutation-constraint}
0 < g_{nearest}(\bm{x}) \Rightarrow 0 < \overline{G(\bm{x})}
\end{equation}
If the refutation test succeeds ($0 < g_{nearest}(\bm{x})$), $\bm{x}$ is infeasible with respect to the interval-based feasibility condition (Equation \ref{eq:individual-feasible}). The exact value of the constraint is however not known.
Otherwise ($g_{nearest}(\bm{x}) \le 0$), the constraint is evaluated using \gls{IA} in order to determine the sign of $\overline{G(\bm{x})}$.


The function \textsc{RigorousFeasibilityTest} (Algorithm \ref{algo:is-feasible-floating}) exploits the floating-point refutation test (Equation \ref{eq:refutation-constraint}) to determine the feasibility of a point $\bm{x}$: the constraints are first evaluated in the round-to-nearest mode using \gls{FPA}, in the hope that the refutation test succeeds at low cost. If \gls{FPA} evaluations cannot refute $\bm{x}$, the constraints are evaluated using \gls{IA}.

\begin{algorithm}[htbp!]
\caption{Rigorous feasibility test combining floating-point arithmetic and interval arithmetic}
\label{algo:is-feasible-floating}
\begin{algorithmic}
\Function{RigorousFeasibilityTest}{$\bm{x}$: individual, $\mathcal{C}$: set of constraints}
\For{$g_j \in \mathcal{C}$}
	\If{$0 < g_j(\bm{x})$}
	\Comment violated constraint
		\State \Return $\mathit{false}$
	\EndIf
\EndFor
\For{$g_j \in \mathcal{C}$}
	\If{$0 < \overline{G_j(\bm{x})}$}
	\Comment violated interval constraint
		\State \Return $\mathit{false}$
	\EndIf
\EndFor
\State \Return $\mathit{true}$
\EndFunction
\end{algorithmic}
\end{algorithm}

To demonstrate the validity of the approach, we compared experimentally the performance of the \gls{DE} using two different versions of the rigorous feasibility test:
\begin{itemize}
\item the "IA" version (Algorithm \ref{algo:is-feasible}) ;
\item the "FPA + IA" version (Algorithm \ref{algo:is-feasible-floating}) exploits the floating-point refutation test (Equation \ref{eq:refutation-constraint}) and the combined evaluation of \gls{FPA} and \gls{IA}.
\end{itemize}

\newpage

Table \ref{tab:convergence-versions-de} compares both approaches on 15 problems from the COCONUT benchmark. The hyperparameters of the \gls{DE} were set to $(\mathit{NP}, W, \mathit{CR}) = (40, 0.7, 0.9)$, except for ex2\_1\_7 ($\mathit{NP} = 20$). The second version ("FPA + IA") proves consistently faster than the first version ("IA"). The gain in CPU time ranges from $5\%$ to $60\%$, and the overall gain is $21\%$.

\begin{table}[htbp!]
	\centering
	\caption{Convergence time for both versions of the rigorous feasibility test}
	\begin{tabular}{|c|cc|c|}
	\hline
				& \multicolumn{2}{c|}{CPU time (s)} & \\
	Problem 	& "IA" version & "FPA + IA" version & Gain (\%) \\
	\hline
	avgasa		& 0.43	& 0.18	& 59 \\
	ex2\_1\_7	& 0.58	& 0.51	& 11.3 \\
	ex2\_1\_10	& 0.76	& 0.63	& 17.4 \\
	ex7\_2\_1	& 0.23	& 0.22	& 5.6 \\
	ex7\_2\_3	& 0.32	& 0.25	& 20.7 \\
	ex7\_2\_4	& 0.22	& 0.20	& 10.1 \\
	ex7\_2\_6	& 0.034	& 0.031	& 7.2 \\
	ex7\_2\_8	& 0.19	& 0.17	& 9.9 \\
	ex7\_2\_9	& 0.65	& 0.50	& 22.7 \\
	expfita		& 0.14	& 0.13	& 5.3 \\
	hexagon		& 0.30	& 0.20	& 33.2 \\
	hs100		& 0.15	& 0.14	& 5.3 \\
	hs118		& 2.23	& 1.74	& 21.7 \\
	keane		& 0.047	& 0.044	& 5.9 \\
	s365mod		& 0.14	& 0.13	& 8.8 \\
	\hline
	Sum		& 6.421	& 5.075 & 21 \\
	\hline
	\end{tabular}
	\label{tab:convergence-versions-de}
\end{table}

\begin{figure}[htbp!]
\centering
\includegraphics[width=0.9\columnwidth]{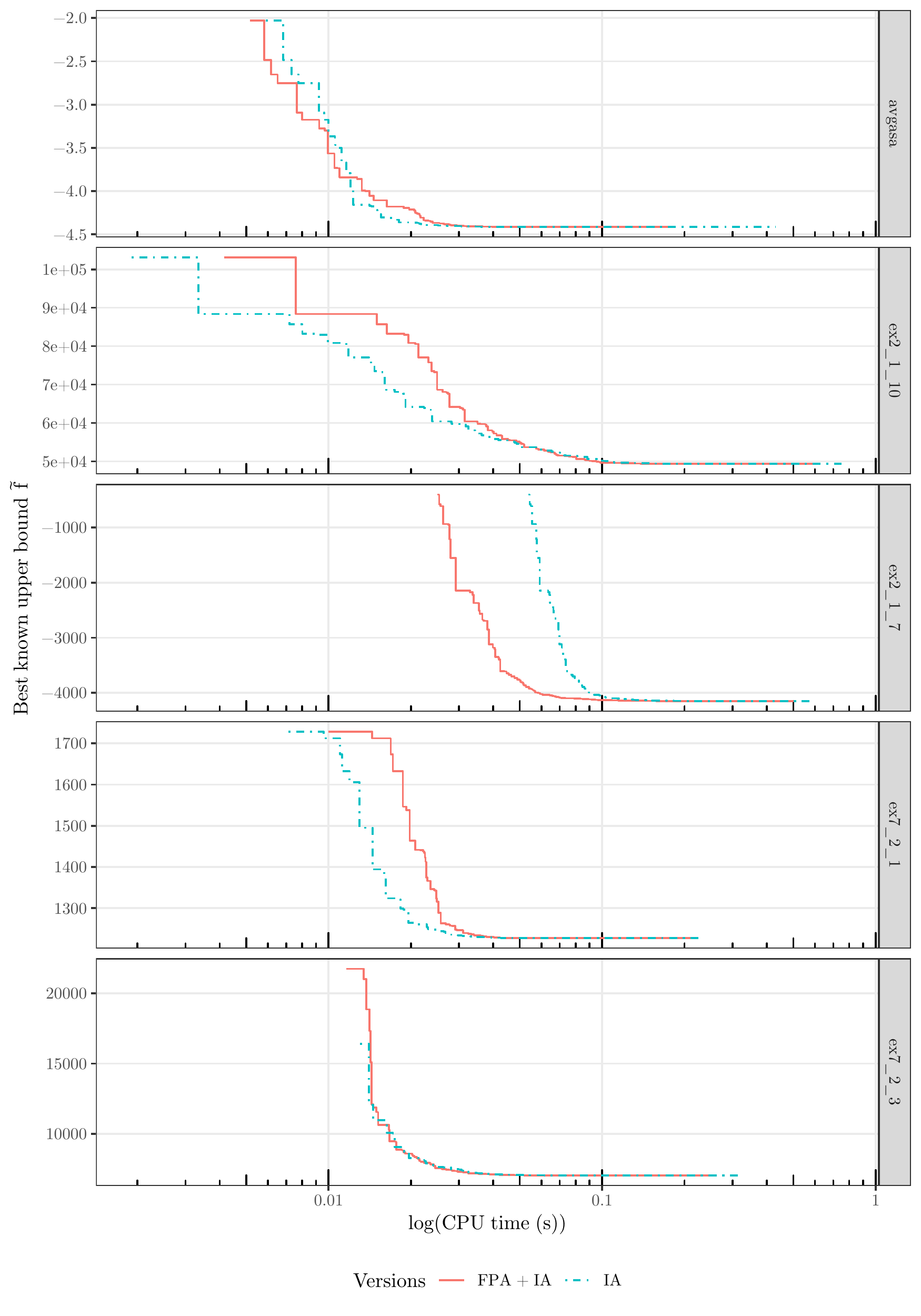}
\caption{Convergence of the differential evolution algorithm with both versions of the rigorous feasibility test (1/3)}
\label{fig:versions-de-1}
\end{figure}

\begin{figure}[htbp!]
\centering
\includegraphics[width=0.9\columnwidth]{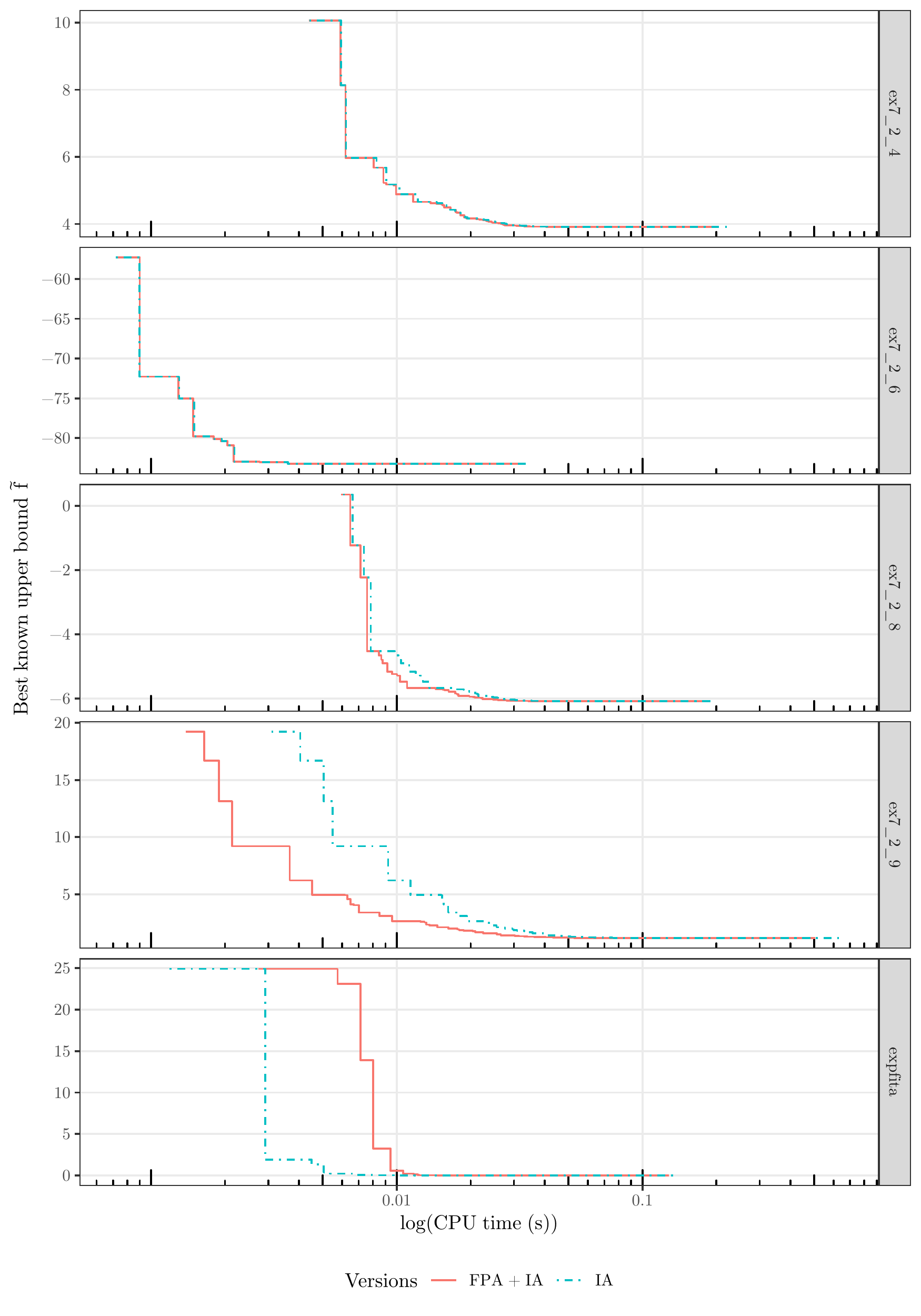}
\caption{Convergence of the differential evolution algorithm with both versions of the rigorous feasibility test (2/3)}
\label{fig:versions-de-2}
\end{figure}

\begin{figure}[htbp!]
\centering
\includegraphics[width=0.9\columnwidth]{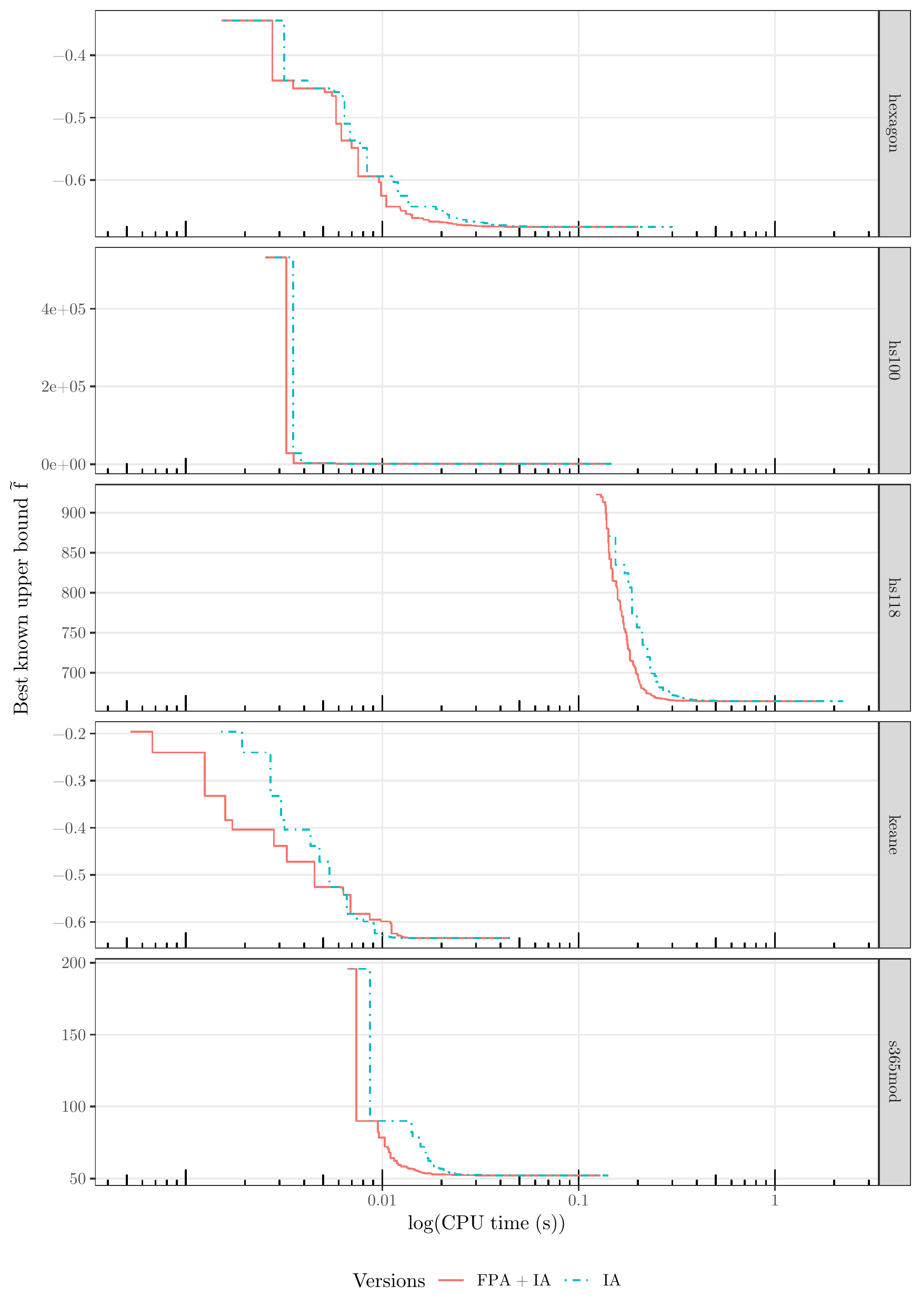}
\caption{Convergence of the differential evolution algorithm with both versions of the rigorous feasibility test (3/3)}
\label{fig:versions-de-3}
\end{figure}

Figures \ref{fig:versions-de-1}, \ref{fig:versions-de-2} and \ref{fig:versions-de-3} illustrate the evolution of the best known upper bound of the \gls{DE} (the evaluation of the best individual) for both versions of the rigorous feasibility test. The results show that the "IA" version starts off faster than the "FPA + IA" version on a few problems (ex2\_1\_10, ex7\_2\_1, expfita), while the "FPA + IA" version has the edge on the other half of the problems (ex2\_1\_7, ex7\_2\_9, keane, s365mod).
This can be explained by the fact that an individual $\bm{y}$ that improves its parent $\bm{x}$ requires $2m$ constraint evaluations ($m$ in \gls{FPA} + $m$ in \gls{IA}): consequently, the number of constraints $m$ and the frequence at which an individual is improved have a direct influence on the time of convergence.

\newpage

\subsection{Exploiting the interval domain}
\label{sec:interval-domain}

\subsubsection{Reduction of the initial domain}
The domain of the \gls{DE} (which corresponds to the initial box of the \gls{IBC}) may be initially contracted in order to discard infeasible values and avoid exploring infeasible or suboptimal areas of the search space.
The initial population of the \gls{DE} is then initialized within the contracted domain.

\begin{figure}[htbp!]
\centering
\includegraphics[width=\textwidth]{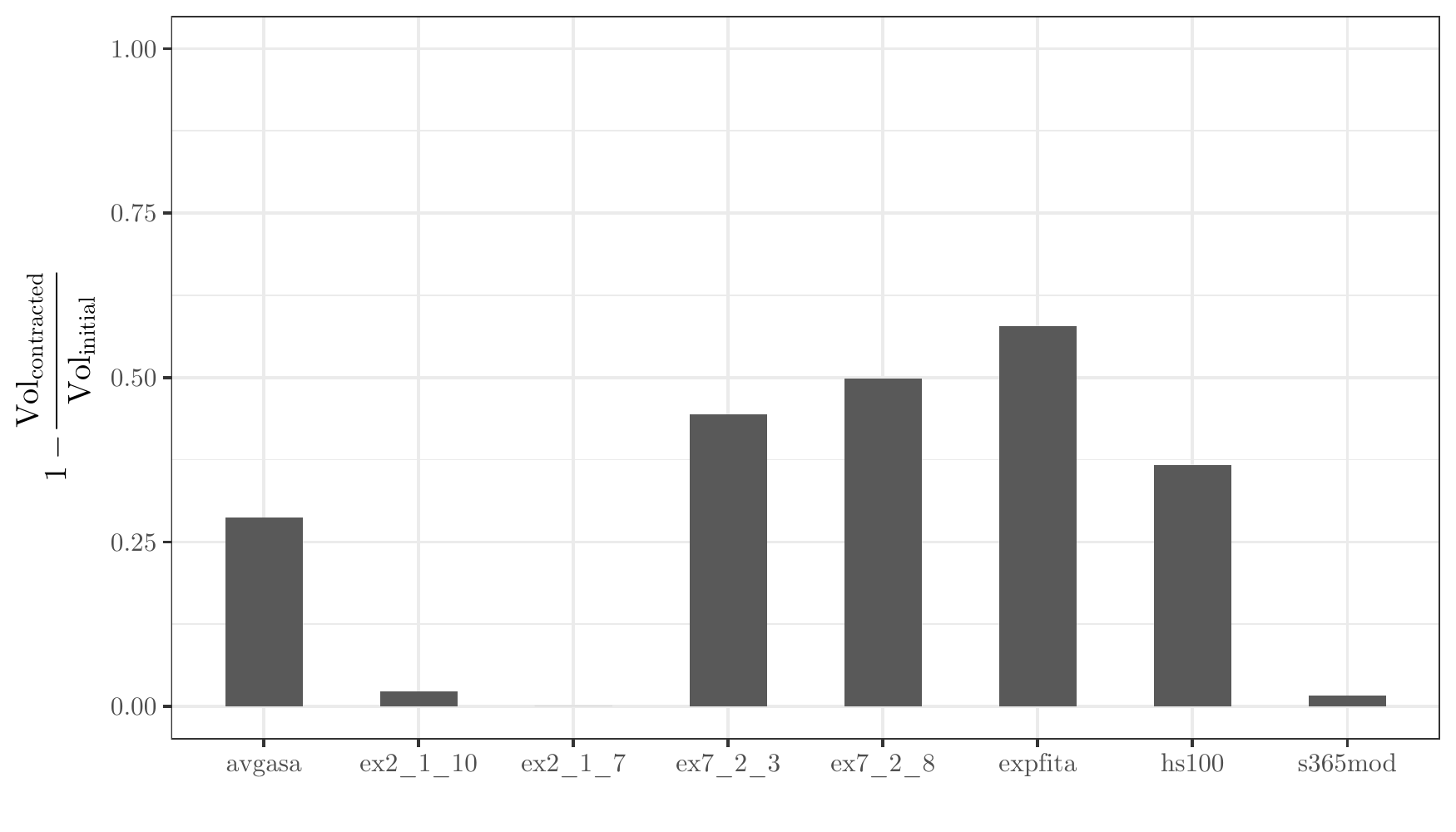}
\caption{Relative reduction of the volume of the initial domain on a subset of test problems}
\label{fig:reduction-initial-domain}
\end{figure}

Figure \ref{fig:reduction-initial-domain} shows the relative reduction of the volume of the initial domain on a subset of COCONUT problems, where the volume of a box $\bm{X} = X_1 \times \ldots \times X_n$ is defined as:
\begin{equation}
Vol(\bm{X}) \eqdef \prod_{i=1}^n w(X_i)
\end{equation}
The initial domains are contracted by a sequence of HC4, 3BCID(HC4) and X-Newton operators. The relative reductions range from 1.7\% (s365mod) to 57.8\% (expfita). Only the initial domain of the ex\_2\_1\_7 problem could not be contracted.

Our approach is similar to that of~\cite{Focacci2003Local} who suggested the use of a preprocessing phase using constraint programming techniques to reduce the initial search space of a local search method.

\subsubsection{Periodic reduction of the domain}
\label{sec:de-domain-reduction}

The priority queue $\mathcal{Q}$ maintained by the \gls{IBC} contains feasible or undecidable subspaces of the search space. This information may be exploited by the \gls{DE} in order to avoid evaluating individuals that lie in infeasible or suboptimal areas.
We propose a strategy that exploits the remaining boxes in the \gls{IBC} in order to regularly update the domain of the \gls{DE} and progressively discard infeasible or suboptimal subspaces. It consists in periodically sending the priority queue $\mathcal{Q}$ of the \gls{IBC} to the \gls{DE} ; the latter computes the convex hull $\square \mathcal{Q}$ of all boxes of $\mathcal{Q}$. The \gls{DE} population is then (possibly randomly) reinitialized within the new domain $\square \mathcal{Q}$.

In this context, the MaxDist search strategy (Section \ref{sec:maxdist}) exhibits two main advantages:
\begin{enumerate}
\item the \gls{IBC} usually maintains a small priority queue $\mathcal{Q}$, which limits the cost of sending $\mathcal{Q}$ to the \gls{DE} as well as the operation of convex hull (both have a linear complexity in the number of boxes) ;
\item by construction, MaxDist processes the boxes at the rim of the domain, which favors the quick reduction of the convex hull of the remaining boxes.
\end{enumerate}

Example \ref{ex:domain-contraction} shows how the size of the \gls{DE} domain is reduced with the generations, while retaining the global minimizer.
\begin{example}
\label{ex:domain-contraction}
Let:
\begin{equation}
\begin{aligned}
\min_{(x, y) \in X \times Y} \quad & -\frac{(x+y-10)^2}{30} - \frac{(x-y+10)^2}{120}  \\
\text{s.t.}	 \quad 	& \frac{20}{x^2} - y \le 0 \\
					& x^2 + 8y - 75 \le 0
\end{aligned}
\end{equation}
be a constrained optimization problem defined on the box $X \times Y = [0, 10] \times [0, 10]$.

\begin{figure}[htbp!]
	\centering
	\subfloat[Initial domain]{\label{fig:initial-domain}
		\includegraphics[width=0.5\columnwidth]{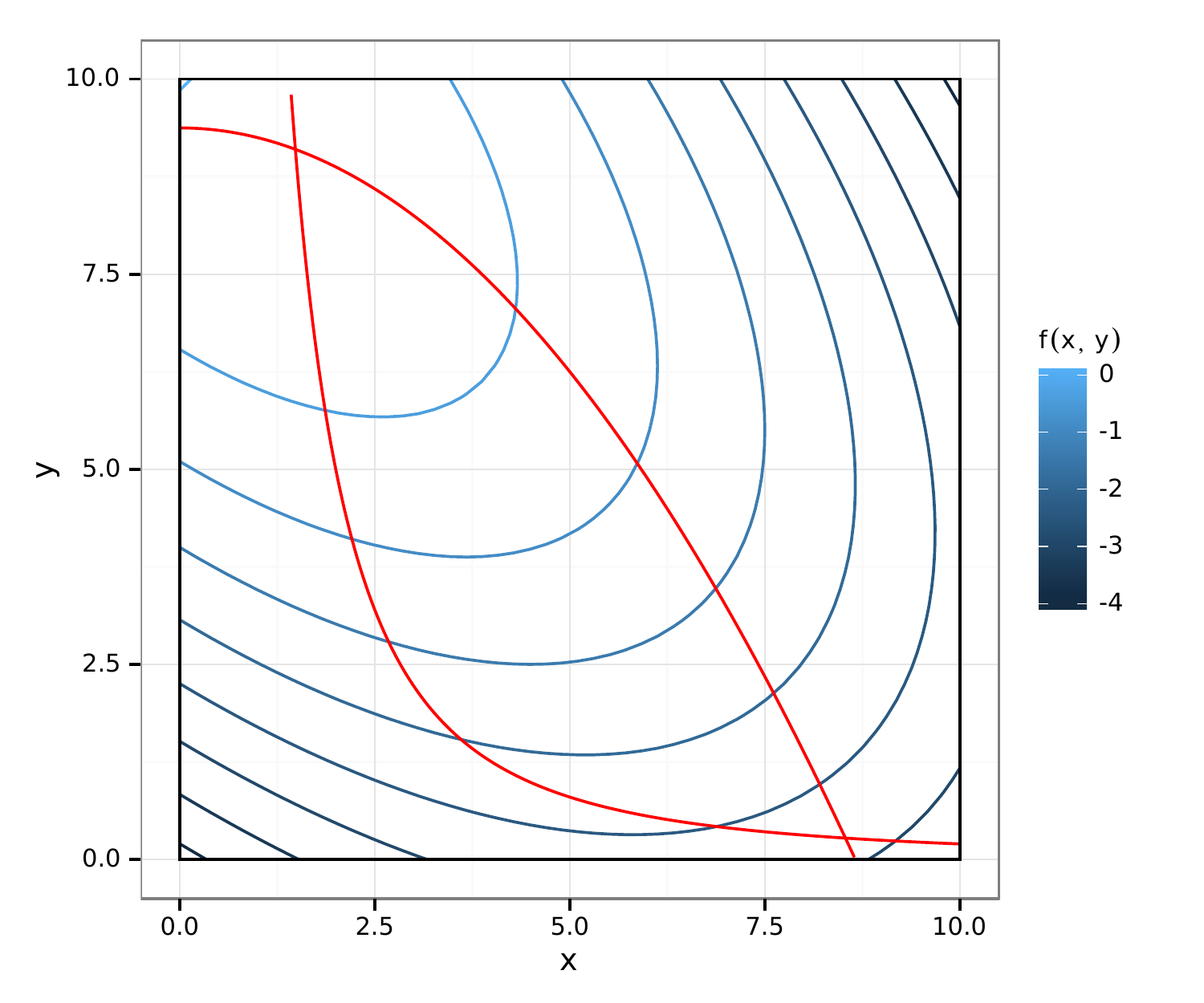}}
	\subfloat[Contracted initial domain]{\label{fig:contracted-domain}
		\includegraphics[width=0.5\columnwidth]{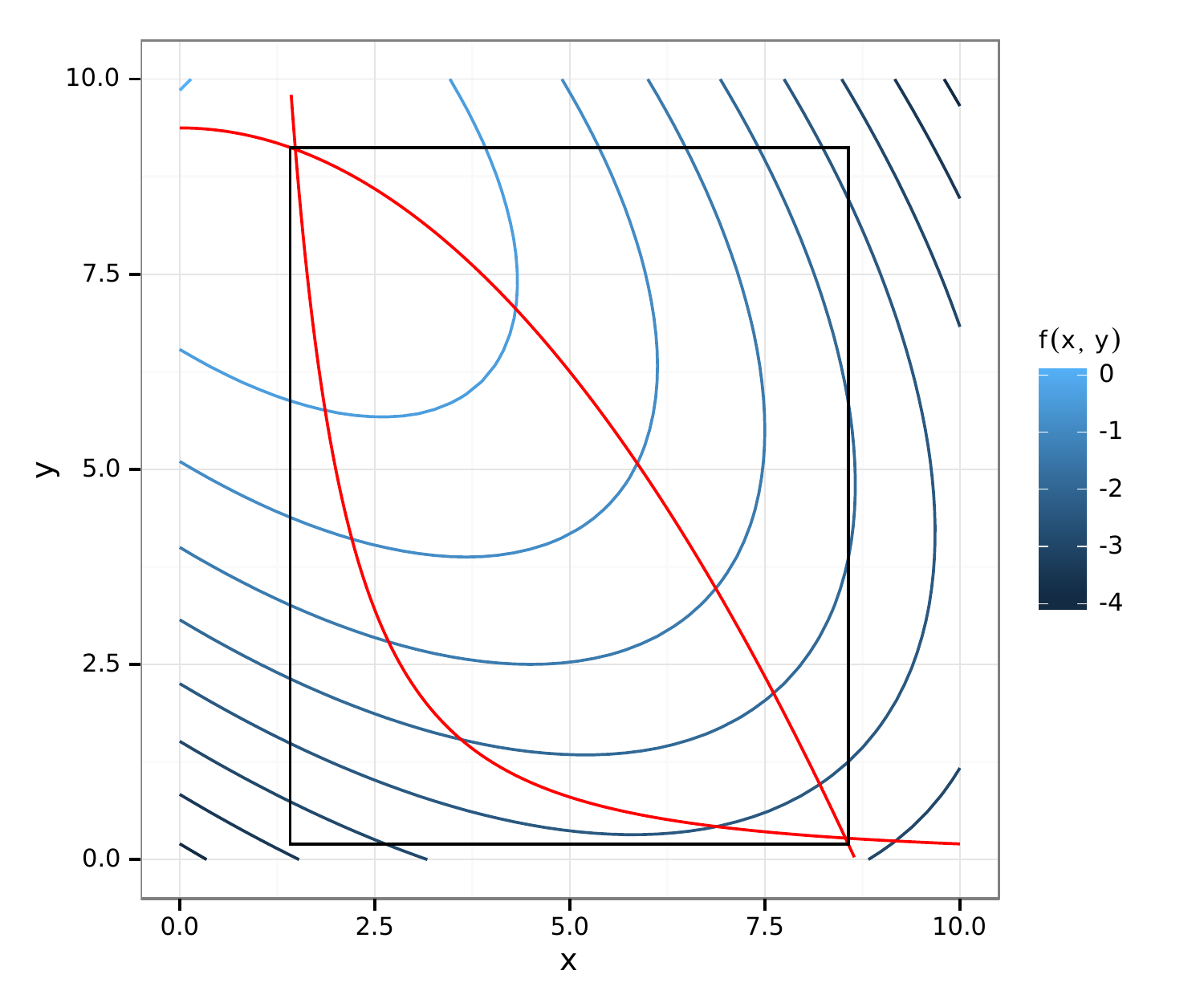}} \quad
	\subfloat[Generation 10]{\label{fig:domain-generation-10}
		\includegraphics[width=0.5\columnwidth]{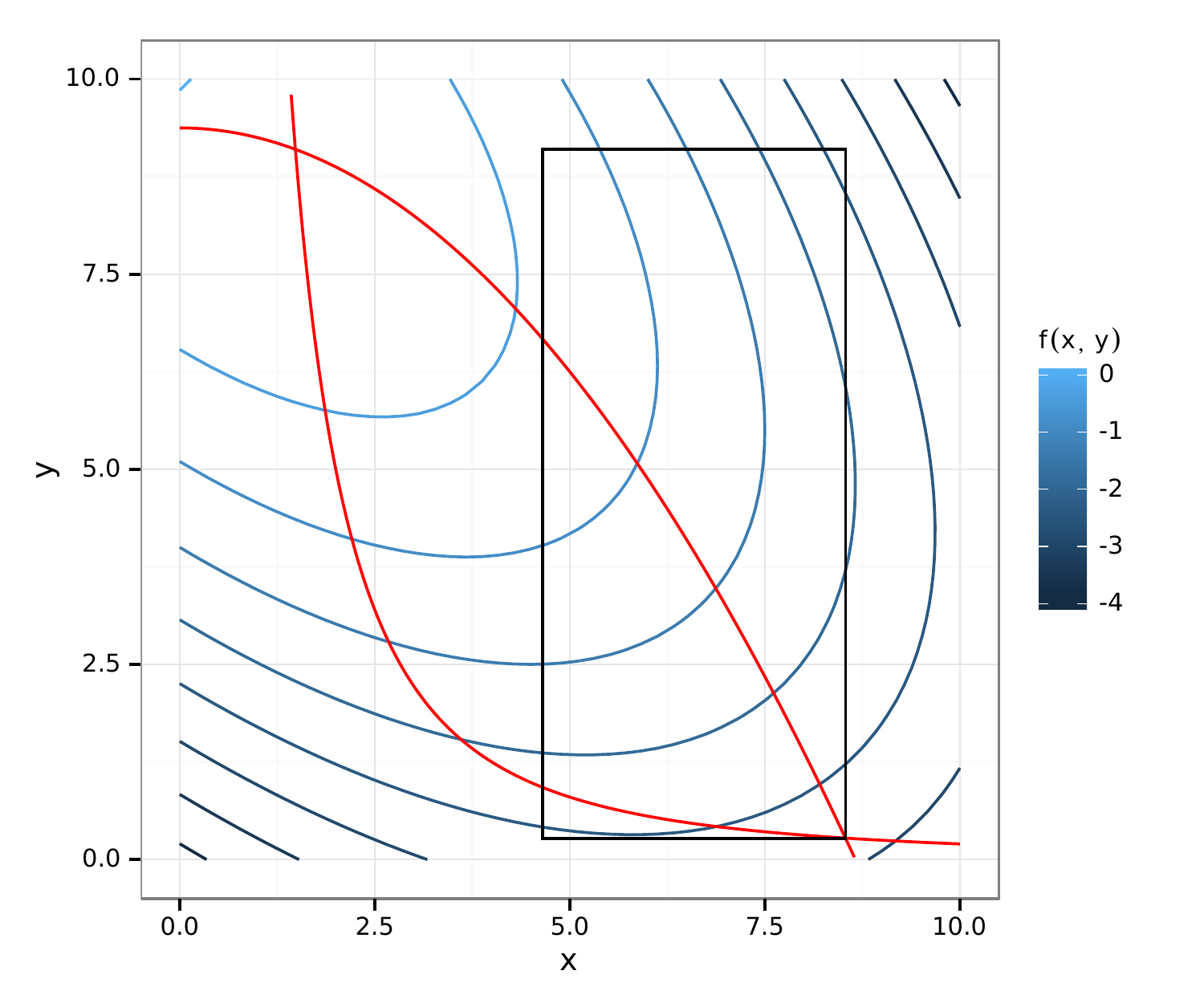}}
	\subfloat[Generation 20]{\label{fig:domain-generation-20}
		\includegraphics[width=0.5\columnwidth]{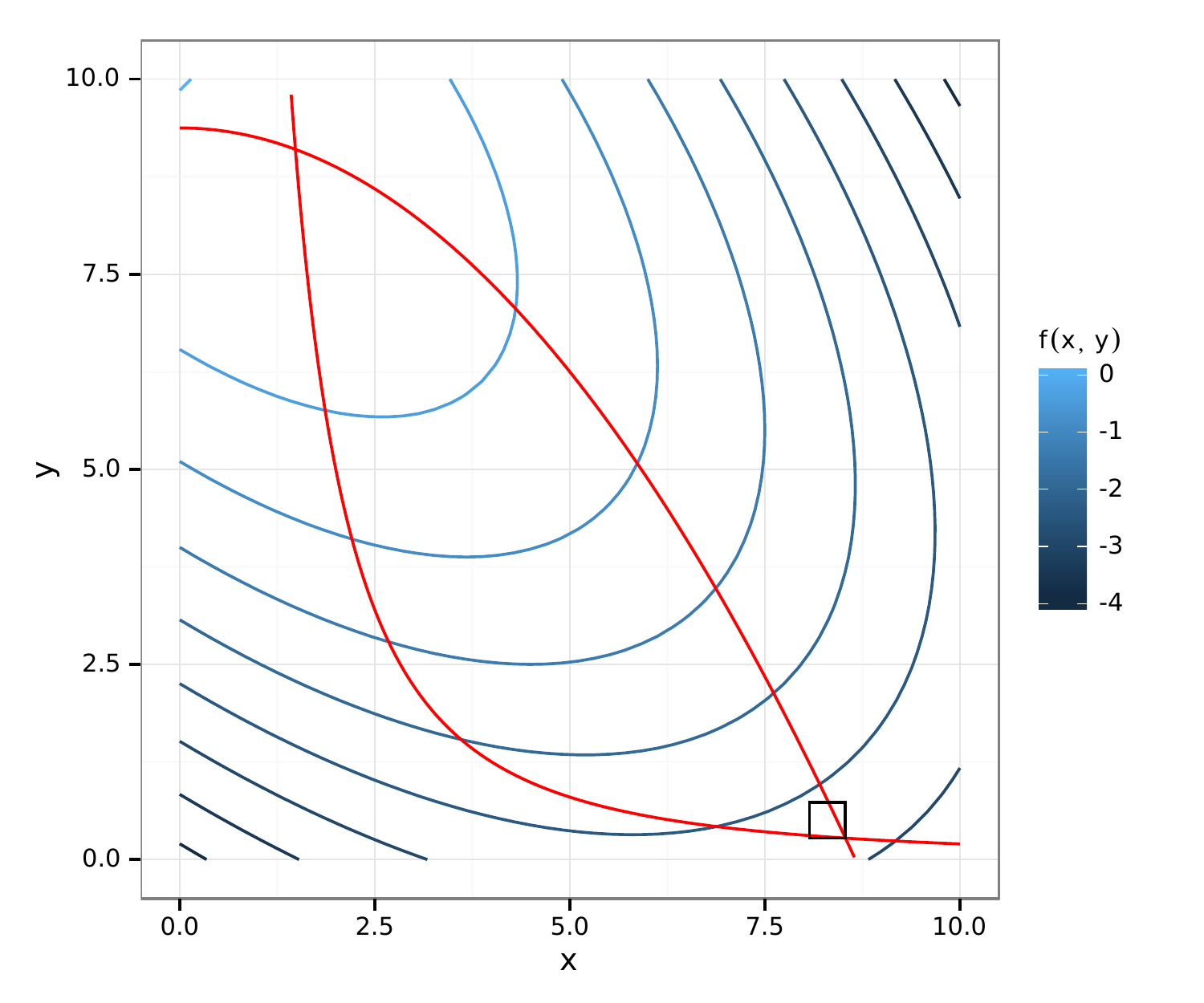}}
	\caption{Contraction of the domain of the differential evolution algorithm with the generations}
\end{figure}

The problem is represented in Figure \ref{fig:initial-domain}: the frontiers of the two inequality constraints are shown in red and the contours of the objective function are shown in blue. The feasible domain is the banana-shaped domain between the two red curves, and the global minimizer is located at its bottom right corner.
Figure \ref{fig:contracted-domain} illustrates the initial domain contracted with respect to the constraints: $(X, Y) = ([1.4142, 8.5674], [0.2, 9.125])$. The two constraints are merely handled one after the other by a 2B operator.
Figures \ref{fig:domain-generation-10} and \ref{fig:domain-generation-20} portray the convex hull of the remaining boxes maintained by the \gls{IBC}, after 10 and 20 \gls{DE} generations respectively. The global minimum found by Charibde with a precision $\varepsilon = 10^{-8}$ is $f(x^*, y^*) = f(8.532424, 0.274717) = -2.825296148$. Both constraints are active at the solution. The analytical expression of the solution is $(x^*, y^*) = (\sqrt{\frac{\sqrt{4985}+75}{2}}, \frac{40}{\sqrt{4985}+75}) \simeq (8.532424404, 0.274716723)$.
\end{example}

\newpage

\begin{figure}[htbp!]
\centering
\includegraphics[width=0.95\columnwidth]{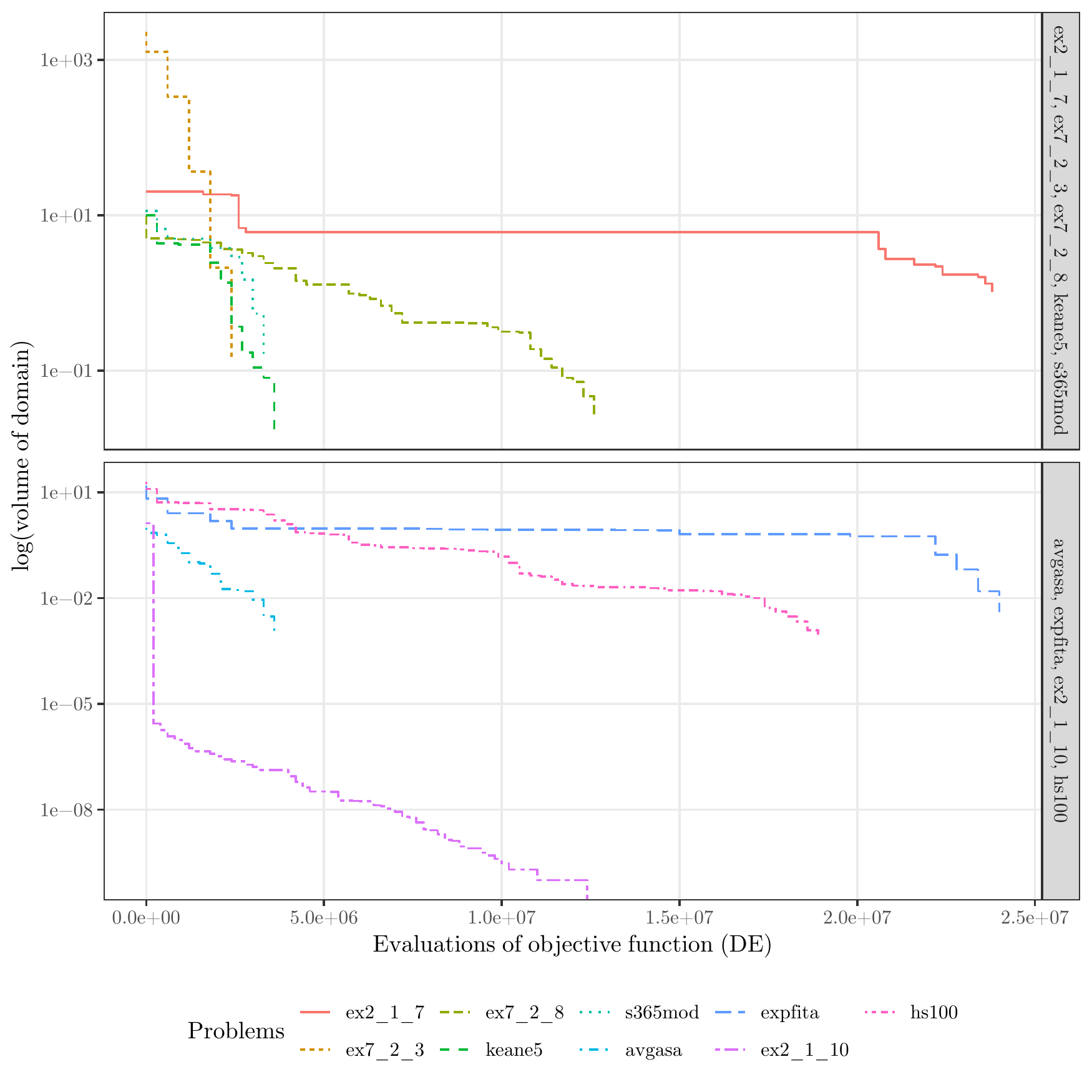}
\caption{Evolution of the volume of the domain (in logarithmic scale) with the number of objective evaluations in the differential evolution algorithm}
\label{fig:de-domain-reduction}
\end{figure}

Figure \ref{fig:de-domain-reduction} portrays the evolution of the size of the \gls{DE} domain (in logarithmic scale) with the number of objective evaluations for 9 COCONUT problems. For these preliminary results, the mechanism that sends the priority queue to the \gls{DE} and reduces its domain is arbitrarily triggered every 20,000 generations of the \gls{DE}.
The figure highlights the efficiency of the approach ; it quickly contracts the domain by several orders of magnitude and considerably reduces the risk that the \gls{DE} stays stuck in local minima.

\subsubsection{Lazy evaluation}

Since the priority queue $\mathcal{Q}$ is now available to the \gls{DE}, evaluating the constraints at a newly generated individual $\bm{x}$ may be carried out in a lazy fashion: if $\bm{x}$ does not belong to a box of $\mathcal{Q}$, it lies outside of the remaining domain of the \gls{IBC} ($\bm{x}$ is therefore either infeasible or suboptimal) and should not be evaluated. Otherwise, only the undetermined constraints (whose evaluation using \gls{IA} contains 0) on the box that contains $\bm{x}$ are evaluated. 
The enumerative type composed of $\mathit{Feasible}$ and $\mathit{Infeasible}$ (see Section \ref{sec:rigorous-constraints}) may be extended: a third option $\mathit{Outside}$ characterizes an individual that lies outside of the remaining domain and therefore has not been evaluated.

Algorithm~\ref{algo:individual-evaluation} describes the lazy evaluation function of the \gls{DE}. If there exists a box $\bm{X}$ in $\mathcal{Q}$ that contains the individual $\bm{x}$, only the undetermined constraints on $\bm{X}$ are evaluated by the procedure \textsc{RigorousFeasibilityMeasure} (Algorithm \ref{algo:rigorous-constraint-evaluation}).
If $\bm{X}$ is a feasible box, $\bm{x}$ is a labeled as a feasible point without evaluating the constraints.
If $\bm{x}$ lies outside of the boxes of $\mathcal{Q}$, the individual is not evaluated and is labeled $\mathit{Outside}$.

\begin{algorithm}[htbp!]
\caption{Evaluation of an individual with respect to the remaining subspaces}
\label{algo:individual-evaluation}
\begin{algorithmic}
\Function{IndividualEvaluation}{$\bm{x}$: individual, $\mathcal{Q}$: priority queue}
\State $(\bm{X}, d_{min}) \gets$ \Call{NearestBox}{$\bm{x}$, $\mathcal{Q}$}
\Comment Algorithm \ref{algo:nearest-box}
\If{$d_{min} = 0$}
	\Comment individual within remaining domain
	\If{$\bm{X}$ subject to undetermined constraints $\mathcal{C}_{\bm{X}}$}
		\State $(\mathit{feasible}, n, s) \gets$ \Call{RigorousFeasibilityMeasure}{$\bm{x}$, $\mathcal{C}_{\bm{X}}$}
		\Comment Algorithm \ref{algo:rigorous-constraint-evaluation}
		\If{$\mathit{feasible}$}
			\State \Return $\mathit{Feasible}(f(\bm{x}))$
		\Else
			\State \Return $\mathit{Infeasible}(n, s)$
		\EndIf
	\Else
		\State \Return $\mathit{Feasible}(f(\bm{x}))$
	\EndIf
\Else
	\State \Return $\mathit{Outside}$
	\Comment individual outside of domain: no evaluation
\EndIf
\EndFunction
\end{algorithmic}
\end{algorithm}

\begin{algorithm}[htbp!]
\caption{Computation of the box nearest to a point}
\label{algo:nearest-box}
\begin{algorithmic}
\Function{NearestBox}{$\bm{x}$: point, $\mathcal{Q}$: priority queue}
\State $\bm{X}_{nearest} \gets \varnothing$
\Comment nearest box
\State $d_{min} \gets +\infty$
\Comment distance between individual and nearest box
\Repeat
	\State $\mathcal{Q} \gets \mathcal{Q} \setminus \{\bm{X}\}$
	\Comment box extraction
	\If{$\bm{x} \in \bm{X}$}
		\Comment the point belongs to the box
		\State $d_{min} = 0$
		\State $\bm{X}_{nearest} \gets \bm{X}$
	\Else
		\State $d \gets $ \Call{Distance}{$\bm{x}, \bm{X}$}
		\Comment Algorithm \ref{algo:distance}
		\If{$d < d_{min}$}
			\State $d_{min} \gets d$
			\State $\bm{X}_{nearest} \gets \bm{X}$
		\EndIf
	\EndIf
\Until{$\mathcal{Q} = \varnothing$ or $d_{min} = 0$}
\State \Return $(\bm{X}_{nearest}, d_{min})$
\EndFunction
\end{algorithmic}
\end{algorithm}

\newpage

\section{Comparison of solvers on the COCONUT benchmark}
\label{sec:coconut-comparison}

GlobSol, IBBA and Ibex are nowadays the most efficient rigorous solvers for constrained optimization. Although they share a common skeleton of \gls{IBB}, they embed distinct contraction techniques.
GlobSol~\cite{Kearfott1996Rigorous,Kearfott2005Validated} computes linear relaxations based on a reformulation-linearization technique.
IBBA~\cite{Ninin2010Reliable} invokes the evaluation-propagation algorithm, and computes a relaxation of the feasible set using the AF2 variant of affine arithmetic.
Ibex~\cite{ChabertJaulin2009} is a powerful solver dedicated to numerical constraint satisfaction and global optimization ; it embeds the majority of the most recent and efficient contractors (HC4, 3B, Mohc, CID, X-Newton).

\subsection{COCONUT benchmark}

\cite{Araya2012Contractor} selected a subset of 11 COCONUT constrained problems that are particularly challenging for rigorous solvers: ex2\_1\_7, ex2\_1\_9, ex6\_2\_6, ex6\_2\_8, ex6\_2\_9, ex6\_2\_11, ex6\_2\_12, ex7\_2\_3, ex7\_3\_5, ex14\_1\_7 and ex14\_2\_7. The smallest instance counts 3 variables and one constraint, while the largest counts 20 variables and 10 constraints.
Due to numerical instabilities ("assert failure"), the linear programming library ocaml-glpk interrupts the execution of Charibde on problems ex6\_1\_1, ex6\_1\_3 et ex\_6\_2\_10, for which no results will be shown.

The test problems are described in Table \ref{tab:coconut-problems}. The second and third columns give respectively the number of variables $n$ and  the  number  of  constraints $m$.  The  fourth  (resp.  fifth)  column specifies the type of the objective function (resp. the constraints): L is linear, Q is quadratic and NL is nonlinear. The logsize of the domain $\bm{D}$ (sixth column) is $\displaystyle \log \left( \prod_{i=1}^n (\overline{D_i} - \underline{D_i}) \right)$.

\begin{table}[htbp!]
	\centering
	\caption{Description of difficult COCONUT problems}
	\begin{tabular}{|l|cc|cc|c|}
	\hline
				& 		&		& \multicolumn{2}{c|}{Type}	& \\
	Problem 	& $n$ 	& $m$	& $f$ 	& $g_j, h_j$ 	& Domain logsize \\
	\hline
	ex2\_1\_7	& 20	& 10	& Q		& L				& $+\infty$	\\
	ex2\_1\_9	& 10	& 1		& Q		& L 			& $+\infty$	\\
	ex6\_2\_6	& 3		& 1		& NL	& L 			& $-3 \cdot 10^{-6}$ \\
	ex6\_2\_8	& 3		& 1		& NL	& L 			& $-3 \cdot 10^{-6}$ \\
	ex6\_2\_9	& 4		& 2		& NL	& L 			& $-2.77$ \\
	ex6\_2\_11	& 3		& 1		& NL	& L			 	& $-3 \cdot 10^{-6}$ \\
	ex6\_2\_12	& 4		& 2		& NL	& L			 	& $-2.77$ \\
	ex7\_2\_3	& 8		& 6		& L		& NL			& $61.90$ \\
	ex7\_3\_5	& 13 	& 15	& L		& NL 			& $+\infty$ \\
	ex14\_1\_7	& 10	& 17	& L		& NL		 	& 23.03 \\
	ex14\_2\_7	& 6		& 9		& L		& NL			& $+\infty$	\\
	\hline
	\end{tabular}
	\label{tab:coconut-problems}
\end{table}

\subsection{Numerical results}

In Table \ref{tab:comparison-solvers}, we compare the rigorous solvers GlobSol, IBBA (results available in \cite{Ninin2010Reliable}), Ibex (results available in \cite{Araya2012Contractor}) and Charibde on the benchmark of 11 optimization problems, along with the non rigorous solvers Couenne and BARON.
For each problem, the first line shows the CPU time (in seconds) and the second line shows the number of bisections.
The value of the precision is identical for all six solvers ($\varepsilon = 10^{-8}$), as well as the relaxation factor for the equality constraints ($\varepsilon_= = 10^{-8}$). TO (timeout) indicates that convergence was not achieved within an hour.

\begin{table}[htbp!]
	\centering
	\caption{Comparison of solvers GlobSol, IBBA, Ibex, Charibde, Couenne and BARON on 11 COCONUT problems: CPU time (in s) and number of bisections}
	\begin{tabular}{|l|cccc|cc|}
	\hline
				& \multicolumn{4}{c|}{Rigorous}									& \multicolumn{2}{c|}{Non rigorous} \\
	Problem 	& GlobSol 	& IBBA	 		& Ibex 		& Charibde 				& Couenne	& BARON \\
	\hline
	ex2\_1\_7	& 			& 16.7			& \bf{7.74}	& 26 					& 476 		& 16.23 \\
				&			& 1574			& 1344		& 67249 				& 			& \\
	ex2\_1\_9	&  			& 154			& 9.07		& 36 					& \bf{3.01}	& 3.58 \\
				&			& 60007			& 5760		& 328056 				& 			& \\
	ex6\_2\_6	& 306		& 1575			& 136		& \textbf{1.96/6.92} 	& TO 		& 5.7 \\
				&			& 922664		& 61969		& 15787 				& 			& \\
	ex6\_2\_8	& 204		& 458			& 59.3		& \textbf{2.77/10.66} 	& TO		& TO \\
				&			& 265276		& 25168		& 23047 				& 			& \\
	ex6\_2\_9	& 463		& 523			& 25.2		& \textbf{2.76/4.54} 	& TO		& TO \\
				&			& 203775		& 27892		& 34591 				& 			& \\
	ex6\_2\_11	& 273		& 140			& 7.51		& \textbf{1.97/2.76} 	& TO		& TO \\
				&			& 83487			& 8498		& 26952 				& 			& \\
	ex6\_2\_12	& 196		& 112			& 22.2		& \textbf{8.7/10.75} 	& TO		& TO \\
				& 			& 58231			& 7954		& 127198 				& 			& \\
	ex7\_2\_3	& 			& TO			& 544		& \bf{1.59} 			& TO		& TO \\
				&			& 				& 611438	& 743 					& 			& \\
	ex7\_3\_5	& 			& TO			& 28.91		& \bf{8}.8 				& TO		& 4.95 \\
				&			&				& 5519		& 36072 				& 			& \\
	ex14\_1\_7	& 			& TO			& 406		& 4		 				& 13.86		& \bf{0.56} \\
				&			&				& 156834	& 8065 					& 			& \\
	ex14\_2\_7	& 			& TO			& 66.39		& 0.3	 				& \bf{0.01}	& 0.02 \\
				&			&				& 12555		& 587 					& 			& \\	
	\hline
	Sum 		& > 1442	& TO			& 1312.32	& \textbf{94.85/112.32} & TO		& TO \\
	\hline
	\end{tabular}
	\label{tab:comparison-solvers}
\end{table}

Two different CPU times are given for Charibde on problems ex6\_2\_6, ex6\_2\_8, ex6\_2\_9, ex6\_2\_11 and ex6\_2\_12: the second value corresponds to the convergence time on the original problem, while the first value is the convergence time on an automatically reformulated version of the problem (see following parapraphs).

The results of GlobSol (commercial solver) are not available for all test problems ; are shown only the results given in \cite{Ninin2010Reliable}.
Only the CPU time of the best Ibex strategy (call to the simplex method, to X-NewIter or to X-Newton) for each problem is provided here ; the reader can refer to \cite{Araya2012Contractor} for more details.
Charibde was run on an Intel Xeon(R) CPU E31270 @ 3.40GHz x 8 with 7.8 GB of RAM.
Couenne and BARON (only the commercial version of the code is available) were run on the NEOS server~\cite{Gropp1997Optimization} on 2 Intel Xeon X5660 @ 2.8GHz x 12 with 64 GB of RAM.
IBBA and Ibex were run on similar processors (Intel x86, 3GHz). The difference in CPU time between computers is about 10\%~\cite{Araya2014Upper}, which makes the comparison quite fair.

The hyperparameters of Charibde for the test problems are shown in Table \ref{tab:parametres-charibde} ; $NP$ is the \gls{DE} population size, the bisection schemes are briefly presented in Section \ref{sec:box-partitioning}, $\eta$ is the fixed-point ratio, and the lower-bounding convexification method and the X-Newton contractor (Section \ref{sec:convexification}) can be toggled on or off.
The scaling factor $W = 0.7$, the crossover rate $CR = 0.9$ and the MaxDist heuristic are common for all test problems. 
Tuning the hyperparameters is generally problem-dependent, and requires structural knowledge about the problem: the population size $\mathit{NP}$ may be set according to the dimension and the number of local minima, the crossover rate $CR$ is related to the separability of the problem, and the techniques based on convexification have little influence for problems with few constraints, but are cheap when the constraints are linear.

\begin{table}[htbp!]
	\centering
	\caption{Hyperparameters of Charibde on 11 COCONUT problems}
	\begin{tabular}{|l|c|cccc|}
	\hline
	Problem 	& $\mathit{NP}$ & Bisection scheme	& $\eta$ 	& Convexification	& X-Newton \\
	\hline
	ex2\_1\_7	& 20			& RR				& 0.9		& \checkmark 		& \checkmark \\
	ex2\_1\_9	& 100			& RR				& 0.8		& \checkmark 		& \\
	ex6\_2\_6	& 30			& Smear				& 0			& \checkmark 		& \\
	ex6\_2\_8	& 30 			& Smear				& 0			& \checkmark 		& \\
	ex6\_2\_9	& 70			& Smear				& 0			& 			 		& \\
	ex6\_2\_11	& 35			& Smear				& 0			& 			 		& \\
	ex6\_2\_12	& 35			& RR				& 0			& \checkmark 		& \\
	ex7\_2\_3	& 40			& Largest			& 0			& \checkmark 		& \checkmark \\
	ex7\_3\_5	& 30			& RR				& 0			& \checkmark 		& \\
	ex14\_1\_7	& 40			& RR				& 0			& \checkmark 		& \\
	ex14\_2\_7	& 40			& RR				& 0			& \checkmark 		& \\
	\hline
	\end{tabular}
	\label{tab:parametres-charibde}
\end{table}

\begin{figure}[htbp!]
\centering
\includegraphics[width=0.9\columnwidth]{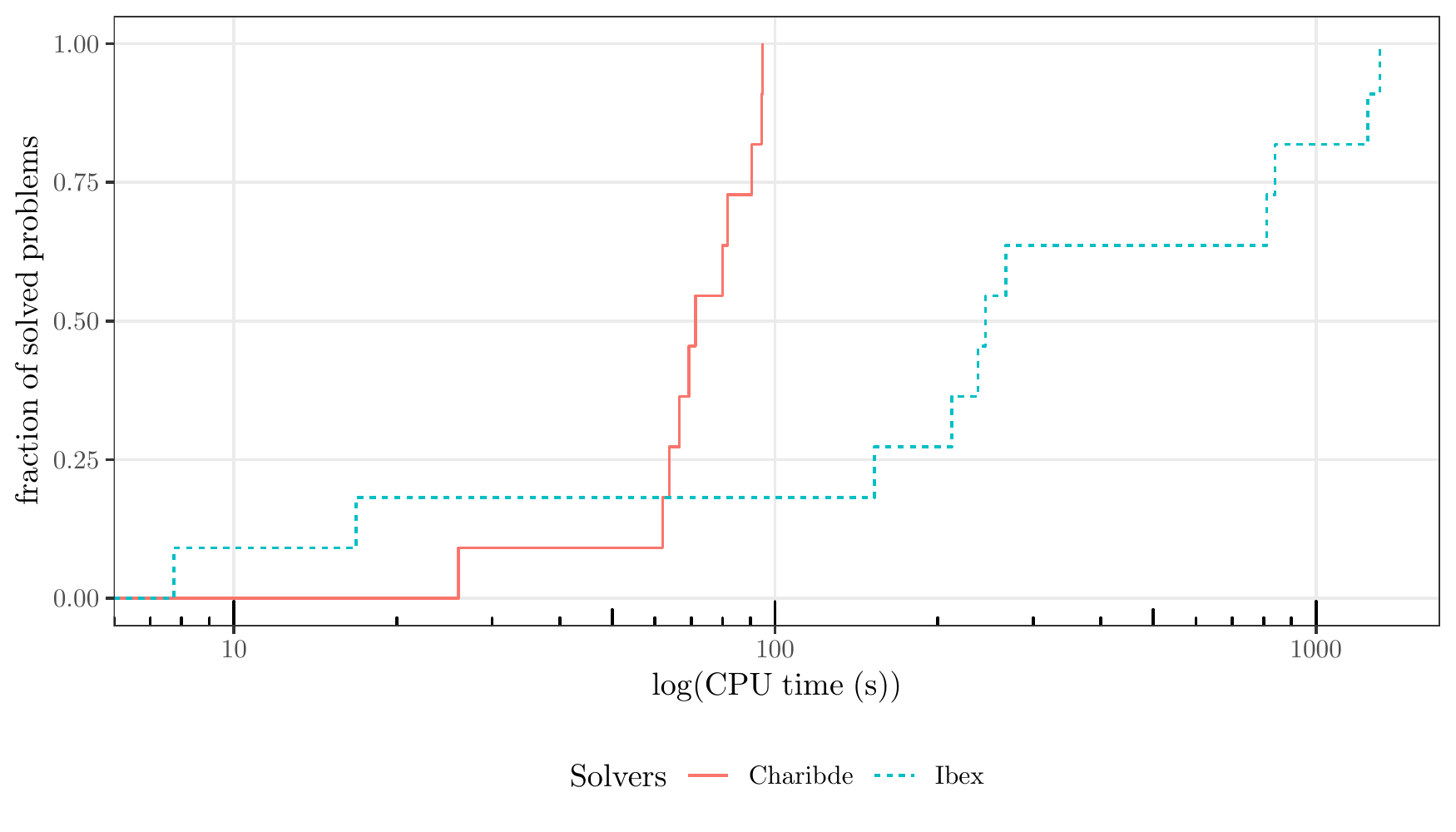}
\caption{Fraction of solved COCONUT problems with respect to time for Ibex and Charibde}
\label{fig:profile-charibde-ibex}
\end{figure}

Charibde surpasses Ibex on 9 out of the 11 test problems, IBBA on 10 out of 11 problems and GlobSol on all the available problems. The cumulated CPU time over the 11 test problems shows that Charibde (112.32s, 94.85s on the reformulated instances) improves the performances of Ibex by an order of magnitude (1312.32s) on this benchmark.
The low number of bisections, in particular on the first two problems, suggests that Ibex adopts a resolution strategy based on strong filtering. On the contrary, Charibde performs more bisections, although the convergence time is often lower (problems ex6\_2\_9, ex6\_2\_11, ex6\_2\_12, ex7\_3\_5).
Figure \ref{fig:profile-charibde-ibex} illustrates how the fraction of solved problems evolves with time (in logarithmic scale) for both Ibex and Charibde.

Charibde also proves highly competitive against the non rigorous solvers Couenne and BARON. The latter are faster or have similar CPU times on some of the 11 problems, however they both timeout on at least five problems (seven for Couenne, five for BARON). Overall, Charibde seems more robust and solves all the problems of the benchmark, while providing a numerical proof of optimality. Surprisingly, the convergence times do not seem directly related to the dimensions of the instances.  They may be explained by the nature of the objective function and constraints (in particular, Charibde seems to struggle  when the objective function is quadratic) and the dependency induced by the multiple occurrences of the variables.

\subsubsection{Automatic reformulation of problems ex6\_2\_6 and ex6\_2\_8}
The objective function was automatically reformulated following a simple rule: factors of variables $x_1, x_2, x_3$ were collected to obtain an expression similar to:
\begin{equation}
f(\bm{x}) = x_1 f_1(x_1, x_2, x_3) + x_2 f_2(x_1, x_2, x_3) + x_3 f_3(x_1, x_2, x_3) 
\end{equation}
where $f_1$, $f_2$ and $f_3$ contain logarithmic terms and constants. This symbolic rewriting reduces the dependency problem in the expression of $f$ and improves the quality of enclosure by \gls{IA}.

\subsubsection{Automatic reformulation of problems ex6\_2\_9, ex6\_2\_11 and ex6\_2\_12}
A reformulation similar to that of problems ex6\_2\_6 and ex6\_2\_8 was carried out. In addition, a basic rewriting of logarithmic terms eliminates the dependency effect within each term:
\begin{equation}
\log(\frac{x_i}{ax_i+bx_j+cx_k}) = \log(\frac{1}{a+\frac{bx_j+cx_k}{x_i}}) =
-\log(a+\frac{bx_j+cx_k}{x_i})
\end{equation}

\subsubsection{Note on problem ex7\_2\_3}
Problem ex7\_2\_3 is an arduous challenge. Even the best solvers struggle to solve it (544s for Ibex) ; \cite{Araya2012Contractor} explains that Ibex converges within reasonable time only when the gradients are computed using Hansen's recursive variant.
However, Charibde reaches convergence in 1.6s: the upper bound $\tilde{f} = 7049.248020528667439$ quickly obtained by the \gls{DE} (4,663 generations, that is 0.31s) intensifies the pruning in the \gls{IBC} in a surprising manner.

It can be shown that, when the best known upper bound is located within the range $\tilde{f} \in [7049.2480205\mathbf{28667439}, 7049.2480205\mathbf{344641}]$, Charibde's \gls{IBC} converges in less than 2s.
This suggests that it reaches fast convergence when the global minimum is approximated with a precision $5.8 \cdot 10^{-9}$, that is yet lower that the asked precision ($\varepsilon = 10^{-8}$)!

\section{Conclusion}
Our cooperative solver Charibde combines the efficiency of \glspl{EA} and the reliability of \gls{IA} to boost the convergence of traditional methods and certify optimality for difficult optimization problems. The \gls{DE} and the \gls{IBC} run individually and exchange best known upper bound, best individual and remaining domain in order to intensify the pruning of the search space and to prevent convergence towards local minima.
Preliminary results show that Charibde improves the performances of Ibex by an order of magnitude on a subset of difficult COCONUT problems.

Much like the algorithm of \cite{Alliot2012Finding}, Charibde is a generic framework, likely to solve optimization problems stemming from various fields of application. Numerous numerical techniques can exploit the structure of the problem (decomposability, differentiability, monotonicity) and greatly accelerate convergence. In particular, combining constraint programming and automatic differentiation (see Section \ref{sec:contraction-ad}) exploits the contraction of the intermediary nodes and computes tighter partial derivatives for gradient-based refutation tests.

\chapter{New global minima of multimodal problems}
\label{sec:new-optima}

\minitoc

In this section, we present new optimality results for five $n$-dimensional test problems that are highly multimodal: Michalewicz~\cite{Michalewicz1996Genetic}, Sine Envelope Sine Wave, Eggholder~\cite{Whitley1996}, Keane~\cite{Keane1994} and Rana~\cite{Whitley1996}. We provide the global minima for various instances and the corresponding solutions. $f_n^*$ denotes the global minimum of $f$ for a given value of the dimension $n$.

Section \ref{sec:multimodal-functions} introduces the expressions, domains and best known solutions of the five test problems. The new optima certified by Charibde are provided in Section \ref{sec:multimodal-functions-results}. We compare Charibde with standalone \gls{DE} and \gls{IBC} algorithms on a particular instance of each problem in Section \ref{sec:multimodal-standalone}. In Section \ref{sec:multimodal-comparison}, we compare Charibde against NLP solvers (Ipopt, LOQO, Minos), metaheuristics (PGAPack, PSwarm) and the spatial branch and bound solver Couenne~\cite{Belotti2009Branching}.

\section{Multimodal test functions}
\label{sec:multimodal-functions}

\begin{enumerate}
\item The Michalewicz function is a separable function, and is highly multimodal (it has around $n!$ local optima). The best known solutions for up to 50 variables were achieved by~\cite{Mishra2006Some} using a repulsive particle swarm optimization algorithm: $f_{10}^* = -9.6602, f_{20}^* = -19.6370, f_{30}^* = -29.6309, f_{50}^* = -49.6248$.

Few results were obtained using global deterministic methods ; the rugged surface of the function is an arduous challenge for global optimization methods. \cite{Alliot2012Finding} claim to have found the global minimum for $n = 12$ with a precision $\varepsilon = 10^{-4}$ in 6,000s: $f_{12}^* = -11.64957$. However, their implementation is flawed: boxes whose size is smaller than $10^{-3}$ are discarded without being further explored. Since the search is not exhaustive, the optimality cannot be guaranteed. In comparison, Charibde achieves convergence on this instance in 0.03s ;

\item The best known solution for the Sine Envelope Sine Wave function is \\
$f_2^* = f_2(-1.1773, -1.6985) = -1.4915$~\cite{Pohl2010Stochastic} ;

\item The best known solution for the Eggholder function is \\
$f_2^* = -959.641$~\cite{Oplatkova2008Metaevolution} ;

\item The best known solutions for Keane's problem are $f_2^* = -0.36497975$, \\
$f_3^* = -0.51578550$, $f_4^* = -0.62228103$ and $f_5^* = -0.63444869$~\cite{Kang2002Level} ;

\item The best known solution for Rana's function is \\
$f_2^* = f_2(-488.63, 512) = -511.7329$~\cite{Tao2007Tao}.
\end{enumerate}

\def\arraystretch{1.6}
\begin{table}[htbp!]
	\centering
	\caption{Expressions and domains of the multimodal test problems}
	\begin{tabular}{|l|l|c|}
	\hline
	Problem 		& Expression														& Domain \\
	\hline
	Michalewicz 	& $-\sum_{i=1}^n \sin(x_i) \left[ \sin(\frac{i x_i^2}{\pi}) \right]^{20}$ & $[0, \pi]^n$ \\
	Sine Envelope 	& $-\sum_{i=1}^{n-1}\left( 0.5 + \frac{\sin^2(\sqrt{x_{i+1}^2 + x_{i}^2} -
						0.5)}{(0.001(x_{i+1}^2 + x_{i}^2) + 1)^2} \right)$ & $[-100, 100]^n$ \\
	Eggholder 		& $-\sum_{i=1}^{n-1} \left[ (x_{i+1}+47) \sin \left(\sqrt{|x_{i+1} + 47 + \frac{x_i}{2}|} \right) + \right.$ & $[-512, 512]^n$ \\
					& $\left. x_i \sin \left(\sqrt{|x_i-(x_{i+1}+47)|} \right) \right]$ & \\
	Keane 			& $-\frac{|\sum_{i=1}^n \cos^4(x_i) - 2\prod_{i=1}^n \cos^2 (x_i)|}{\sqrt{\sum_{i=1}^n i x_i^2}}$ & $[0, 10]^n$ \\
					& s.t. $0.75 \le \prod_{i=1}^n x_i$ and $\sum_{i=1}^n x_i \le 7.5n$ & \\
	Rana 			& $\sum_{i=1}^{n-1} \left( x_i \cos \sqrt{|x_{i+1}+x_i+1|} \sin \sqrt{|x_{i+1}-x_i+1|} + \right.$ & $[-512, 512]^n$ \\
					& $\left. (1+x_{i+1}) \sin \sqrt{|x_{i+1}+x_i+1|} \cos \sqrt{|x_{i+1}-x_i+1|} \right)$ & \\
	\hline
	\end{tabular}
	\label{tab:multimodal-functions}
\end{table}
\def\arraystretch{1}

\begin{figure}[htbp!]
\subfloat[Michalewicz]{\includegraphics[width=0.45\textwidth]{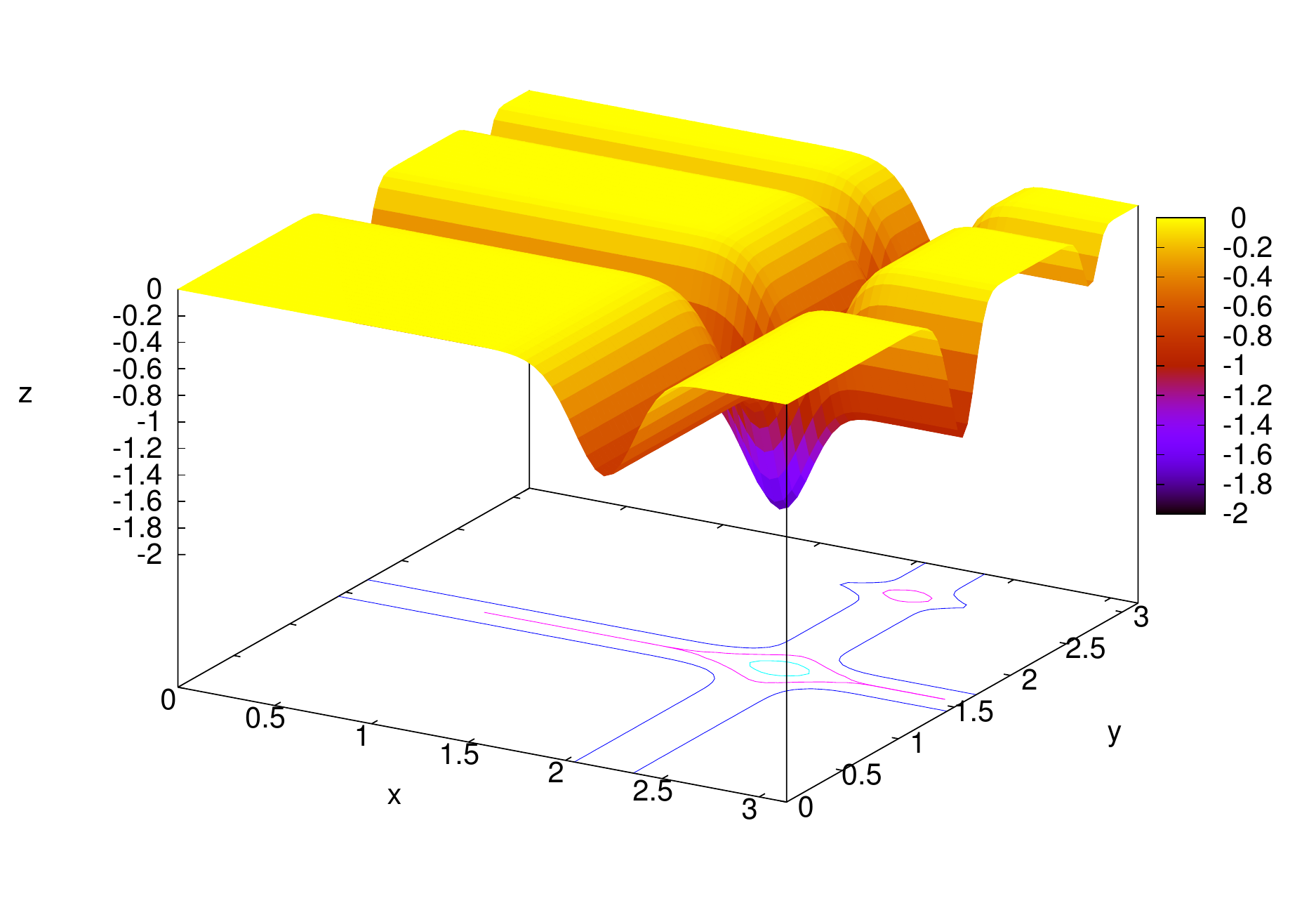}}
\subfloat[Sine Envelope Sine Wave]{\includegraphics[width=0.45\textwidth]{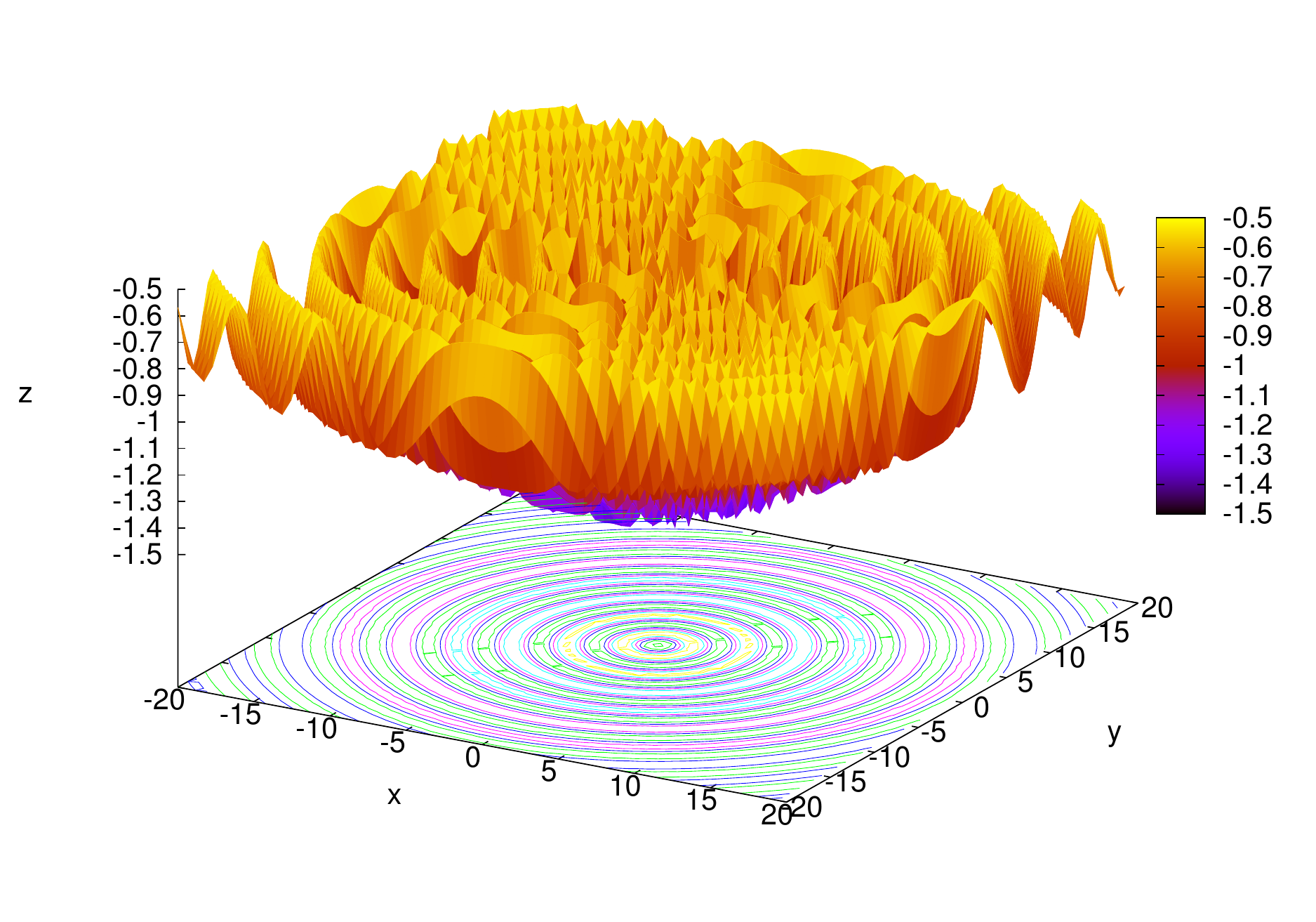}} \quad
\subfloat[Eggholder]{\includegraphics[width=0.45\textwidth]{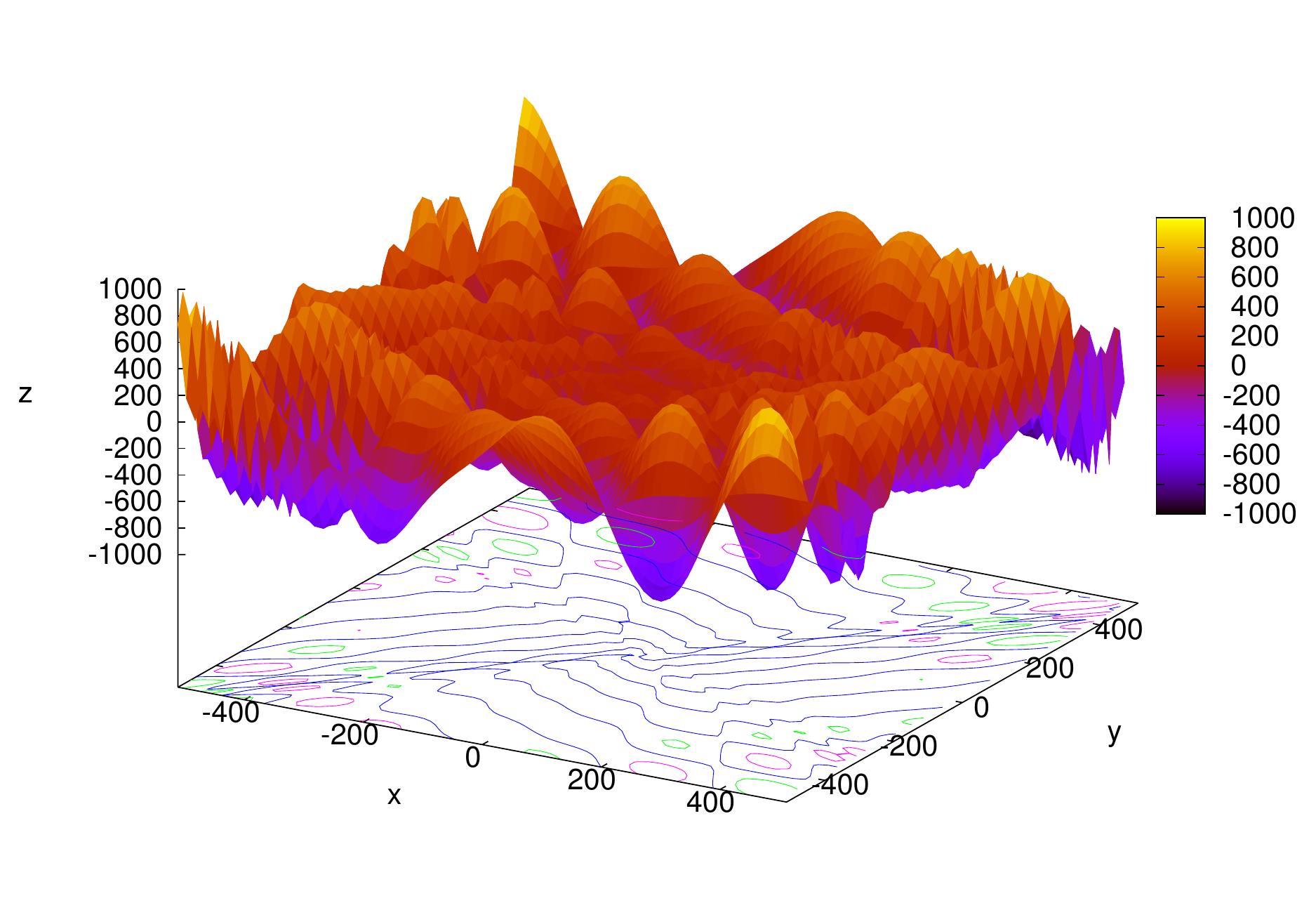}}
\subfloat[Keane]{\includegraphics[width=0.45\textwidth]{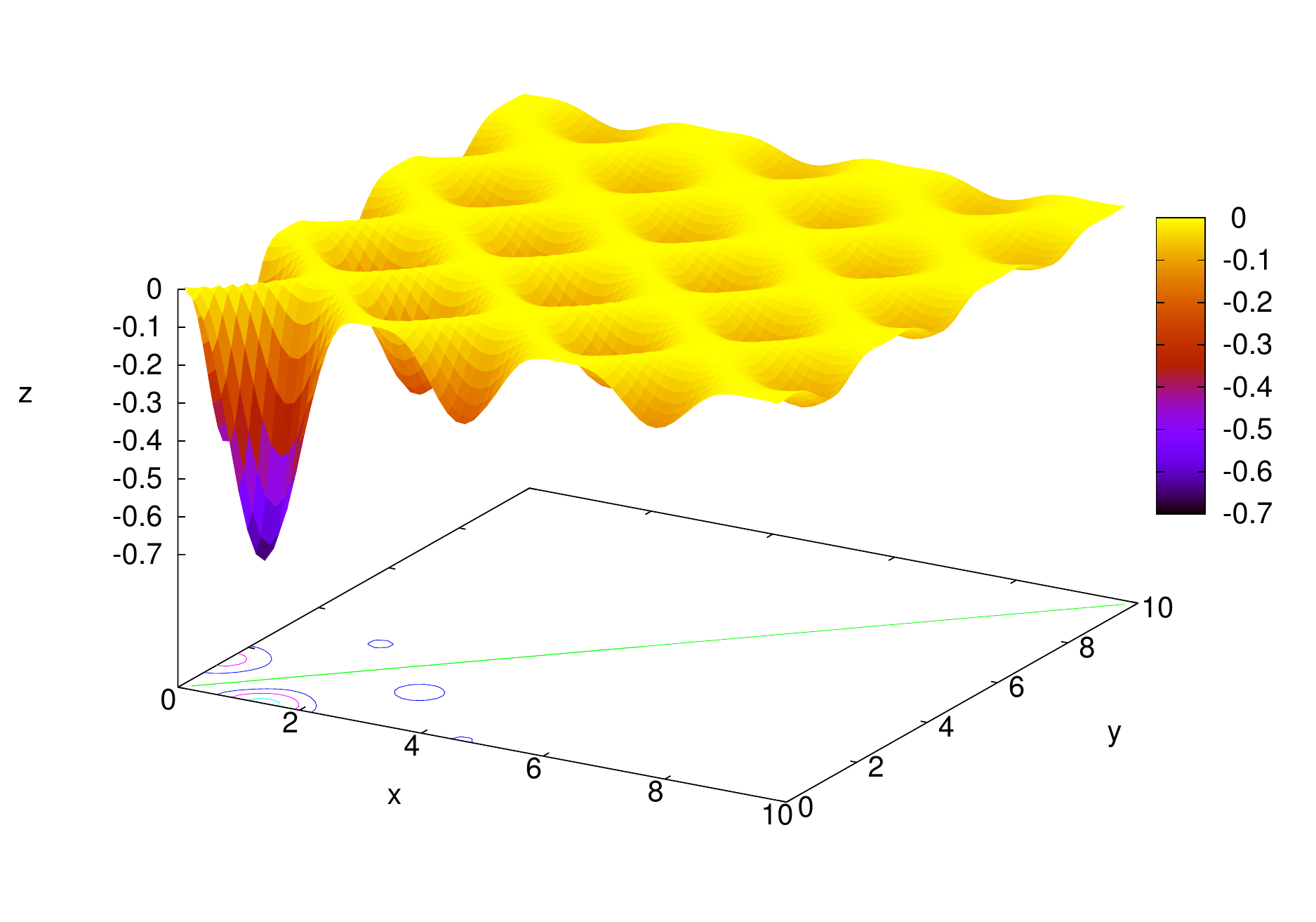}} \quad
\subfloat[Rana]{\includegraphics[width=0.45\textwidth]{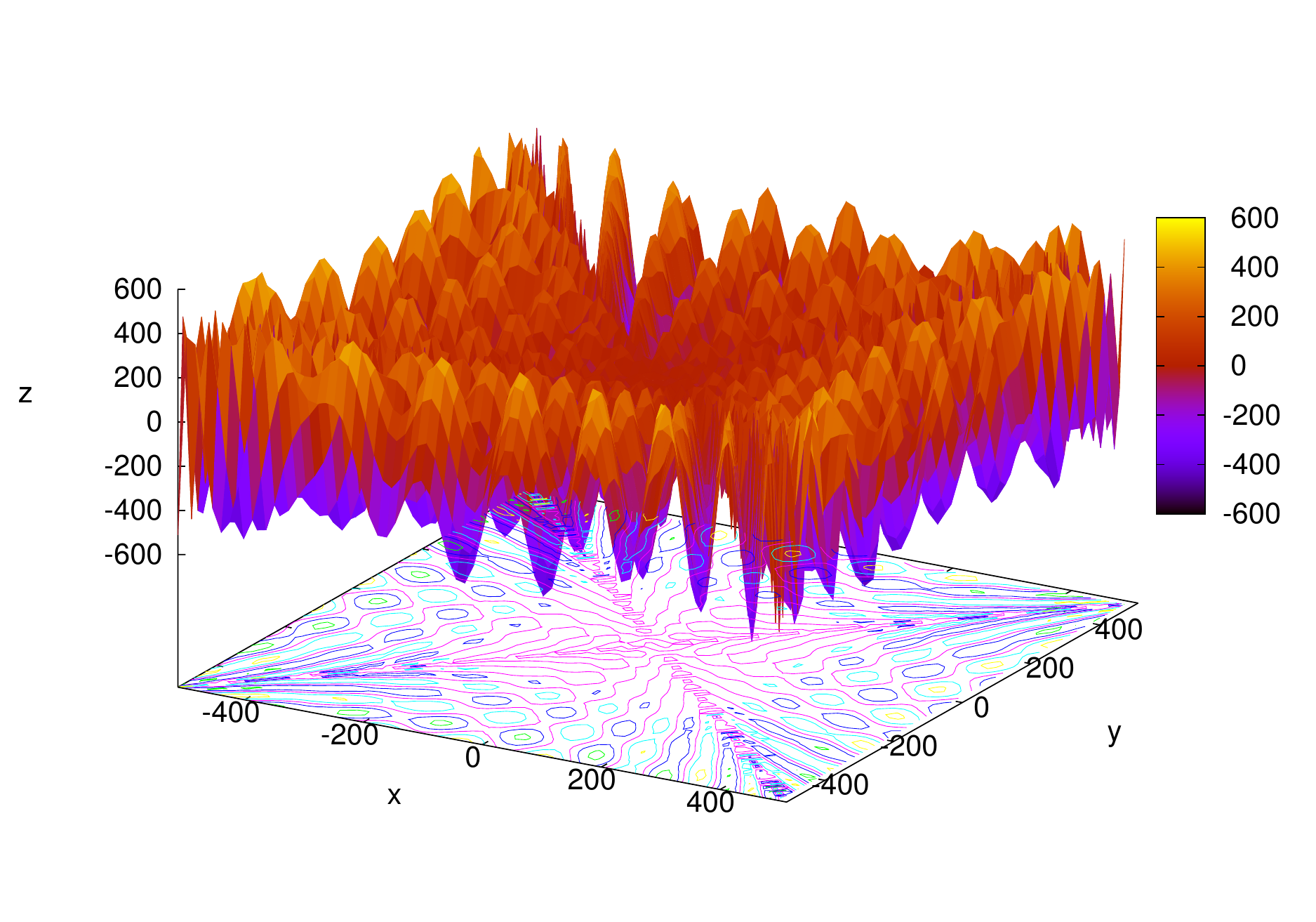}}
\caption{Multimodal test problems ($n = 2$)}
\label{fig:multimodal-problems}
\end{figure}

The expressions and domains of the five multimodal test problems are given in Table \ref{tab:multimodal-functions}. Keane's problem is the only nonlinearly constrained problem ; the others are bound constrained problems. Figure \ref{fig:multimodal-problems} portrays the surfaces of the test problems for $n = 2$.

\clearpage

\section{Numerical results}
\label{sec:multimodal-functions-results}
Charibde converges within a reasonable time on the following instances:
\begin{itemize}
\item up to 70 variables for the Michalewicz function ;
\item up to 5 variables for the Sine Envelope Sine Wave function ;
\item up to 10 variables for the Eggholder function ;
\item up to 5 variables for Keane's problem ;
\item up to 7 variables for Rana's function.
\end{itemize}

The problems Eggholder, Keane and Rana contain absolute values ; the absolute value function is not differentiable at $x = 0$, which hinders the computation of gradients. Based on the suggestion of \cite{Kearfott1996Interval} to use its subgradient, Charibde computes a straightforward enclosure of gradients using \gls{AD}:
\begin{equation}
| \cdot |'(X) =
	\begin{cases}
	[-1, -1] & \text{if } \overline{X} < 0 \\
	[1, 1] & \text{if } \underline{X} > 0 \\
	[-1, 1] & \text{otherwise}
	\end{cases}
\label{eq:subgradient}
\end{equation}

Table \ref{tab:multimodal-stats} details the average and maximum CPU times (in seconds) and the average number of evaluations (NE) of the objective function $f$, its interval extension $F$ and the gradient $\nabla F$ after 100 runs of Charibde.
The hyperparameters of Charibde used on the various instances are given in Table \ref{tab:multimodal-hyperparameters}.

\begin{table}[htbp]
	\centering
	\small
	\caption{Hyperparameters of Charibde for the multimodal test problems}
	\begin{tabular}{|l|ccccc|}
	\hline
	Hyperparameter 					& Michalewicz 	& Sine Envelope & Eggholder	& Keane 	& Rana \\
	\hline
	$\varepsilon$ (precision)		& $10^{-8}$ 	& $10^{-6}$ 	& $10^{-8}$ & $10^{-8}$ & $10^{-8}$ \\
	$\mathit{NP}$ (population size)	& 10 to 70 		& 50			& 50 		& 30 		& 50 \\
	$W$ (scaling factor)			& 0.7 			& 0.7 			& 0.7 		& 0.7 		& 0.7 \\
	$\mathit{CR}$ (crossover rate)	& 0 			& 0.9 			& 0.4 		& 0.9 		& 0.5 \\
	Bisection scheme				& RR 			& RR 			& Largest 	& Largest 	& RR \\
	Search strategy					& MaxDist 		& MaxDist 		& MaxDist 	& MaxDist 	& MaxDist \\
	$\eta$ (fixed-point ratio)		& 0 			& 0.8 			& 0.8 		& 0 		& 0.9\\
	\hline
	\end{tabular}
	\label{tab:multimodal-hyperparameters}
\end{table}

\def\arraystretch{1.05}
\begin{table}[htbp!]
	\centering
	\small
	\caption{Average and maximum CPU time (in seconds) and average number of evaluations (NE) after 100 runs of Charibde for the multimodal test problems}
	\begin{tabular}{|c|c|cccccc|}
	\hline
	& $n$ & Average time & Max. time & NE $f$ & NE $F$ (DE) & NE $F$ (IBC) & NE $\nabla
		F$ \\
	\hline
	\multirow{7}{*}{\begin{turn}{90}Michalewicz\end{turn}} &
 	  10	& 0.018 & 0.022 & 2,601 & 54.5 & 1,015.6 & 633 \\
 	& 20	& 0.092 & 0.13 & 12,987 & 182.1 & 8,552.9	& 6,623.6 \\
 	& 30	& 0.32 & 0.39 & 41,016 & 397	& 33,403.8 & 28,260.8 \\
	& 40	& 1.37 & 1.49 & 141,978 & 655.9 & 146,354.4 & 129,950.1 \\
	& 50	& 5.09 & 6.39 & 427,344 & 852.7 & 515,656.9 & 468,812.3 \\
	& 60	& 26.08 & 34.61 & 1,815,786.4 & 1,106.6 & 2,806,303.5 & 2,608,734 \\
	& 70	& 94.14 & 113.3 & 5,363,505.3 & 1,288.8 & 10,088,068.6 & 9,482,448.7 \\
	\hline
	\multirow{4}{*}{\begin{turn}{90}Sine\end{turn}} &
	  2 & 0.33 & 0.57 & 697,975 & 9.7 & 102,105.4 & 45,449.2 \\
 	& 3 & 2.10 & 2.15 & 2,657,588.1 & 16.6 & 473,732.4 & 264,384.9 \\
	& 4 & 24.29 & 25,43 & 21,673,266.5 & 31.3 & 4,144,215.5 & 2,543,524 \\
	& 5 & 194.25 & 195.7 & 142,043,888.5 & 51.4 & 29,603,730.4 & 19,362,164.1 \\
	\hline
	\multirow{9}{*}{\begin{turn}{90}Eggholder\end{turn}} &
	  2 & 0.0035 & 0.005 & 3,325 & 9.1 & 484.3 & 163.2 \\
	& 3 & 0.05 & 0.067 & 73,830 & 69.1 & 9,703.7 & 5,648.1 \\
	& 4 & 0.18 & 0.26 & 201,585 & 105.1 & 28,070.2 & 16,266.1 \\
	& 5 & 1.72 & 1.87 & 1,523,307.6 & 146 & 270,582.5 & 181,441.3 \\
	& 6 & 4.45 & 5.17 & 2,920,999 & 582.7 & 628,228.7 & 413,157.5 \\
	& 7 & 8.37 & 8.49 & 5,115,327.3 & 146.7 & 1,107,179.1 & 716,550.8 \\
	& 8 & 28.52 & 28.7 & 13,250,573.8 & 158.2 & 3,478,698.3 & 2,208,158 \\
	& 9 & 185.47 & 187.42 & 76,566,411.1 & 550 & 20,526,451.9 & 15,084,418.4 \\
	& 10 & 606.44 & 621.11 & 229,408,972.7 & 504.3 & 64,625,870.1 & 46,654,258.9 \\
	\hline
	\multirow{4}{*}{\begin{turn}{90}Keane\end{turn}} &
	  2 & 0.012 & 0.023 & 9,824.9 & 99.8	& 653.9 & 154.2 \\
 	& 3 & 0.047 & 0.061 & 43,334.3	& 70.8 & 3,140.8 & 698.6 \\
	& 4 & 0.41 & 0.74 & 386,745.8 & 157.3 & 31,236.5 & 6,863.2 \\
	& 5 & 2.72 & 3.11 & 1,804,713 & 158.8 & 188,803.6 & 42,587.4 \\
	\hline
	\multirow{5}{*}{\begin{turn}{90}Rana\end{turn}} &
	  2 & 0.011 & 0.014 & 12,855 & 39.8 & 1,513.9 & 439.3 \\
	& 3 & 0.13 & 0.13 & 199,770 & 41 & 21,645 & 8,068.1 \\
	& 4 & 1.44 & 1.56 & 1,613,409.2 & 64.7 & 187,906.9 & 78,977.9 \\
	& 5 & 18.85 & 19.13 & 15,819,144.5 & 89 & 2,055,106.5 & 957,616.6 \\
	& 6 & 247.29 & 248.98 & 161,239,558.6 & 78.5 & 23,158,792.9 & 11,753,033 \\
	\hline
	\end{tabular}
	\label{tab:multimodal-stats}
\end{table}
\def\arraystretch{1}

\clearpage

\subsection{Note on the Michalewicz function}
The separability of the function allows to find global minima on quite large instances. Choosing a value $\mathit{CR} = 0$ for the crossover rate has a decisive influence on the convergence of the \gls{DE} in Charibde: a newly generated individually will differ from its parent by only one component. It becomes straightforward to independently optimize terms of $f$ that depend on a single variable. 

\subsection{Note on Rana's function}
Charibde converges within a reasonable time for up to 4 variables on the original syntax (Table \ref{tab:multimodal-functions}). However, rewriting the expression de $f$ drastically reduces the computation time and allows to solve instances with up to 7 variables. Applying the trigonometric identity $\forall (u, v) \in \mathbb{R}^2, \cos u \sin v = \frac{1}{2}(\sin(u+v)-\sin(u-v))$ to $f_n$, we get:
\begin{equation}
\begin{split}
f_n(\bm{x}) = 	& \quad \frac{1}{2} \sum_{i=1}^{n-1} \left( (x_{i+1}+1+x_i)
\sin(\sqrt{|x_{i+1}+x_i+1|} + 
\sqrt{|x_{i+1}-x_i+1|}) - \right. \\
				& \left. (x_{i+1}+1-x_i) \sin(\sqrt{|x_{i+1}-x_i+1|}-\sqrt{|x_{i+1}+x_i+1|}) \right)
\label{eq:rana2}
\end{split}
\end{equation}

Table \ref{tab:rana-comparison} compares the convergence time of Charibde between the original syntax and the proposed reformulation (Equation \ref{eq:rana2}).
TO (timeout) indicates that convergence was not achieved within an hour.

\begin{table}[htbp]
	\centering
	\caption{Convergence time of Charibde for two syntaxes of Rana's function}
	\begin{tabular}{|c|cc|}
	\hline
	$n$ & \multicolumn{2}{c|}{CPU time (s)} \\
		& Original syntax & Reformulation \\
	\hline
	2 & 0.25 	& 0.009 \\
	3 & 6.5	 	& 0.12 \\
	4 & 254		& 1.45 \\
	5 & TO		& 18.5 \\
	6 & TO 		& 244 \\
	7 & TO 		& 3300 \\
	\hline
	\end{tabular}
	\label{tab:rana-comparison}
\end{table}

\subsection{Global minima and corresponding solutions}
The global minima and the corresponding solutions of the Sine Envelope Sine Wave, Eggholder, Keane and Rana problems are provided in Table \ref{tab:multimodal-minima}. Charibde provided a numerical proof of optimality for a precision $\varepsilon = 10^{-8}$.

\def\arraystretch{1.05}
\begin{table}[htbp!]
	\centering
	\caption{Global minima and corresponding solutions of the Sine Envelope Sine Wave, Eggholder, Keane and Rana problems $(\varepsilon = 10^{-8})$}
	\begin{tabular}{|c|c|l|l|}
	\hline
	& $n$ & $f_n^*$ & $x_n^*$ \\
	\hline
	\multirow{4}{*}{\begin{turn}{90}Sine\end{turn}} &
	  2 & -1.4914953 & (-0.086537, 2.064868) \\
	& 3 & -2.9829906 & (1.845281, -0.930648, 1.845281) \\
	& 4 & -4.4744859 & (2.066680, 0.001365, 2.066680, 0.001422) \\
	& 5 & -5.9659811 & (-1.906893, -0.796823, 1.906893, 0.796823, -1.906893) \\
	\hline
	\multirow{14}{*}{\begin{turn}{90}Eggholder\end{turn}} &
	  2 & -959.6406627 	& (512, 404.231805) \\
	& 3 & -1888.3213909 & (481.462894, 436.929541, 451.769713) \\
	& 4 & -2808.1847922 & (482.427433, 432.953312, 446.959624, 460.488762) \\
	& 5 & -3719.7248363 & (485.589834, 436.123707, 451.083199, 466.431218, 421.958519) \\
	& 6 & -4625.1447737 & (480.343729, 430.864212, 444.246857, 456.599885, 470.538525, \\
	& 	&				& 426.043891) \\
	& 7 & -5548.9775483 & (483.116792, 438.587598, 453.927920, 470.278609, 425.874994, \\
	& 	&				& 441.797326, 455.987180) \\
	& 8 & -6467.0193267 & (481.138627, 431.661180, 445.281208, 458.080834, 472.765498, \\
	& 	&				& 428.316909, 443.566304, 457.526007) \\
	& 9 & -7376.2797668 & (482.785353, 438.255330, 453.495379, 469.651208, 425.235102, \\
	& 	& 				& 440.658933, 454.142063, 468.699867, 424.215061) \\
	& 10 & -8291.2400675 & (480.852413, 431.374221, 444.908694, 457.547223, 471.962527, \\
	& 	&				& 427.497291, 442.091345, 455.119420, 469.429312, 424.940608) \\
	\hline
	\multirow{4}{*}{\begin{turn}{90}Keane\end{turn}} &
	  2 & -0.3649797 & (1.600860, 0.468498) \\
	& 3 & -0.5157855 & (3.042963, 1.482875, 0.166211) \\
	& 4 & -0.6222810 & (3.065318, 1.531047, 0.405617, 0.393987) \\
	& 5 & -0.6344487 & (3.075819, 2.991995, 1.475794, 0.236691, 0.233309) \\
	\hline
	\multirow{6}{*}{\begin{turn}{90}Rana\end{turn}} &
	  2 & -511.7328819  & (-488.632577, 512) \\
	& 3 & -1023.4166105 & (-512, -512, -511.995602) \\
	& 4 & -1535.1243381 & (-512, -512, -512, -511.995602) \\
	& 5 & -2046.8320657 & (-512, -512, -512, -512, -511.995602) \\
	& 6 & -2558.5397934 & (-512, -512, -512, -512, -512, -511.995602) \\
	& 7 & -3070.2475210 & (-512, -512, -512, -512, -512, -512, -511.995602) \\
	\hline
	\end{tabular}
	\label{tab:multimodal-minima}
\end{table}
\def\arraystretch{1}

\clearpage

The global minima and the corresponding solutions of the Michalewicz function for $n \in \{10, 20, 30, 40, 50, 60, 70\}$ are provided in Table \ref{tab:micha-minima} and Table \ref{tab:micha-solutions}, respectively. Since the function is separable, only the solution for $n = 70$ is given ; the solution for a smaller instance of size $k < 70$ can be inferred by taking the $k$ first components of the given vector.

\begin{table}[htb]
	\centering
	\caption{Global minima of the Michalewicz function $(\varepsilon = 10^{-8})$}
	\begin{tabular}{|cc|}
	\hline
	$n$ & $f_n^*$ \\
	\hline
	10 & -9.66015171564 \\
	20 & -19.63701359935 \\
	30 & -29.63088385032 \\
	40 & -39.62674886468 \\
	50 & -49.62483231828 \\
	60 & -59.62314622857 \\
	70 & -69.62222020764 \\
	\hline
	\end{tabular}
	\label{tab:micha-minima}
\end{table}

\begin{table}[htb]
	\centering
	\caption{Global minimizer of the Michalewicz function for $n = 1$ to 70}
	\begin{tabular}{|ccccc|}
	\hline
	(2.202905, & 1.5707963, & 1.2849915, & 1.9230584, & 1.7204697, \\
	1.5707963, & 1.4544139, & 1.7560865, & 1.6557174, & 1.5707963, \\
	1.4977288, & 1.6966163, & 1.6300760, & 1.5707963, & 1.5175461, \\
	1.6660645, & 1.6163286, & 1.5707963, & 1.5289070, & 1.6474563, \\
	1.6077572, & 1.5707963, & 1.5362725, & 1.6349315, & 1.6019018, \\
	1.5707963, & 1.5414351, & 1.6259253, & 1.5976479, & 1.5707963, \\
	1.5452545, & 1.6191375, & 1.5944175, & 1.5707963, & 1.5481947, \\
	1.6138382, & 1.5918810, & 1.5707963, & 1.5505278, & 1.6095861, \\
	1.5898364, & 1.5707963, & 1.5524243, & 1.6060986, & 1.5881533, \\
	1.5707963, & 1.5539962, & 1.6031866, & 1.5867435, & 1.5707963, \\
	1.5553204, & 1.6007184, & 1.5855456, & 1.5707963, & 1.5564510, \\
	1.5985997, & 1.5845151, & 1.5707963, & 1.5574277, & 1.5967613, \\
	1.5836191, & 1.5707963, & 1.5582799, & 1.5951509, & 1.5828331, \\
	1.5707963, & 1.5590300, & 1.5937286, & 1.5821378, & 1.5707963) \\
	\hline
	\end{tabular}
	\label{tab:micha-solutions}
\end{table}

\begin{figure}[htbp!]
\centering
\includegraphics[width=0.9\columnwidth]{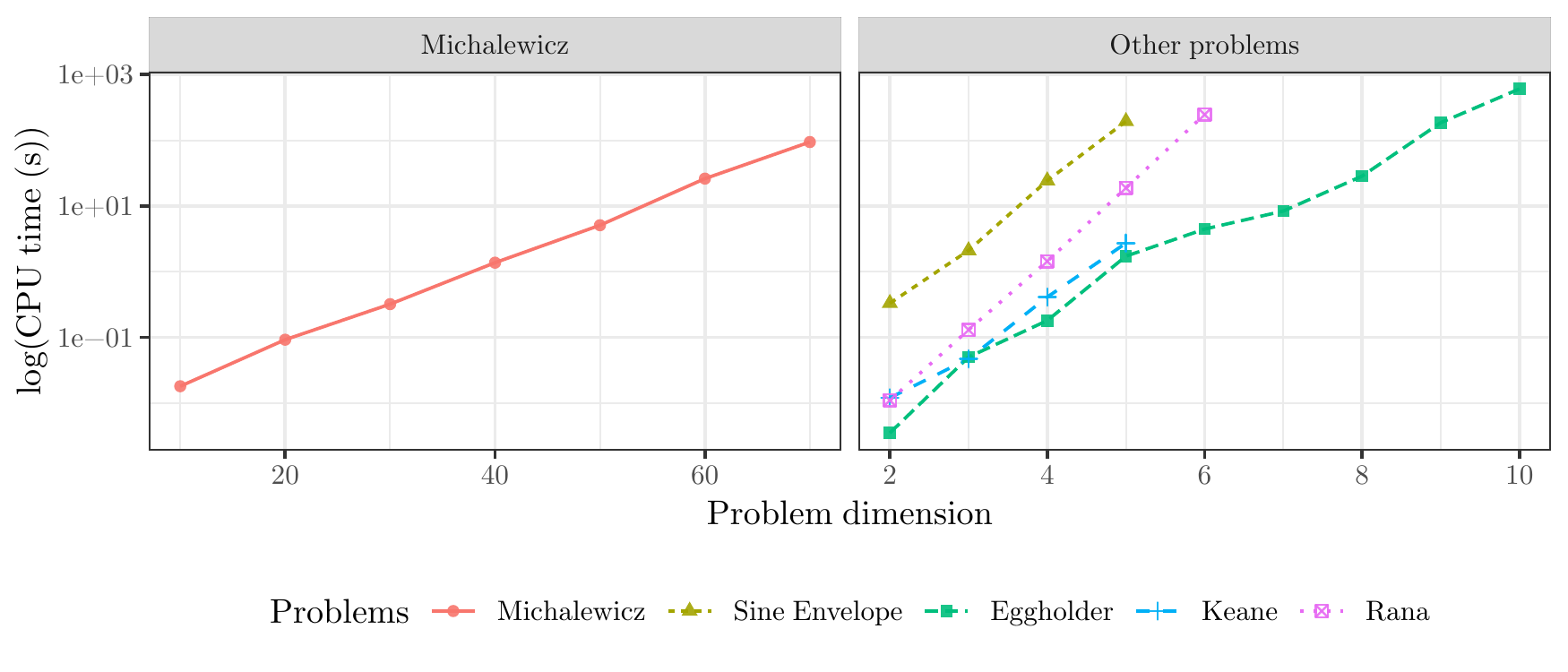}
\caption{Convergence time against the dimension on the multimodal test problems}
\label{fig:time-against-size}
\end{figure}

Figure \ref{fig:time-against-size} illustrates the average convergence time of Charibde (in logarithmic scale) plotted against the size of the instances of the test problems ; it corroborates the exponential complexity of \gls{IBB} algorithms with the number of variables.

\clearpage

\subsection{Interpolation}

Table \ref{tab:multimodal-minima} and Table \ref{tab:micha-minima} suggest that the values of the global minimum $f_n^*$ of four out of five problems satisfy affine equations in $n$. We carried out a linear regression for each of the Michalewicz, Sine Envelope Sine Wave, Eggholder and Rana problems ; the affine expressions and the corresponding coefficients of determination $R^2$ can be found in Table \ref{tab:regressions}.

\begin{table}[htbp]
	\centering
	\caption{Affine relationship between the global minima of multimodal problems and the number of variables}
	\begin{tabular}{|lcc|}
	\hline
	Problem 		& $f_n^*$ 							& $R^2$ \\
	\hline
	Michalewicz		& $-0.9994729257n + 0.3467746311$	& 0.9999998949 \\
	Sine Envelope	& $-1.49150n + 1.49150$				& 1 \\
	Eggholder		& $-915.61991n + 862.10466$			& 0.9999950 \\
	Rana			& $-511.70430n + 511.68714$			& 1 \\
	\hline
	\end{tabular}
	\label{tab:regressions}
\end{table}

\section{Benefits of the hybridization}
\label{sec:multimodal-standalone}

Table \ref{tab:multimodal-standalone} compares Charibde against its standalone components (a \gls{DE} and an \gls{IBC}) on a particular instance of each of the test problems. The "best-first search" and random (the next box to explore is randomly chosen) search strategies performed best for the standalone \gls{IBC}, while Charibde consistently performed best with the MaxDist strategy (Section \ref{sec:maxdist}). TO (timeout) indicates that convergence was not achieved within an hour.

\begin{landscape}
\def\arraystretch{1.05}
\begin{table}[htbp!]
	\centering
 	\footnotesize
	\caption{Comparison of Charibde against its standalone components on the multimodal test problems}
	\begin{tabular}{|c|l|c|c|c|c|c|}
	\hline
	 &				& Michalewicz 	& Sine Envelope & Eggholder  	& Keane  		& Rana \\
	&				& ($n = 50$)	& ($n = 5$)		& ($n = 7$)		&($n = 4$)		& ($n = 5$) \\
	\hline
	\multirow{7}{*}{\begin{turn}{90}DE\end{turn}} &
	$\tilde{f}$ 	& -49.6237668 	& -5.9659811 	& -5548.9775483 & -0.6222810	& -2021.9520121 \\
	& CPU time (s)	& 4.22			& 0.060			& 0.51			& 0.028			& 0.076 \\
	& NE $f$		& 335,580 		& 49,900		& 350,350  		& 13,500		& 51,650 \\
	& NE $F$		& 878 			& 74 			& 142  			& 74	  		& 93 \\
	& $\mathit{NP}$ & 60			& 50			& 50 			& 30			& 50 \\
	& $W$ 			& 0.7			& 0.7 			& 0.7 			& 0.7 			& 0.7 \\
	& $\mathit{CR}$ & 0 			& 0.9			& 0.4 			& 0.9 			& 0.5 \\

	\hline
	\multirow{9}{*}{\begin{turn}{90}IBC\end{turn}} &
	$\tilde{f}$		& -16.09385		& \bf{-5.9659811} & \bf{-5548.9775483} & \bf{-0.6222810} & \bf{-2046.8320657} \\
	& CPU time (s) 	& TO 			& 269.9 		& 12.7 			& 0.56			& 75.76\\
	& Bisections 	& - 			& 1,934,639 	& 105,659 		& 6,059  		& 756,812 \\
	& NE $F$	 	& - 			& 29,766,788 	& 1,475,165 	& 42,755 		& 6,240,504 \\
	& NE $\nabla F$ & - 			& 19,462,956 	& 983,730 		& 9,180			& 3,269,639 \\
	& $|\mathcal{Q}|_{max}$ & - 	& 870,564		& 50,330		& 4,605			& 481,099 \\
	& Bisection scheme 	& RR 		& RR 			& Largest 		& Largest 		& RR \\
	& Search strategy & random 		& best-first 	& best-first 	& best-first 	& best-first \\
	& $\eta$ (fixed-point ratio)	& 0.9 			& 0.8 			& 0.9 			& 0 			& 0.9 \\
	
	\hline
	\multirow{15}{*}{\begin{turn}{90}Charibde\end{turn}} &
	$\tilde{f}$		& \bf{-49.6248323} & \bf{-5.9659811} & \bf{-5548.9775483} & \bf{-0.6222810} & \bf{-2046.8320657} \\
	& CPU time (s) 	& 6.3 			& 194.9 		& 7.85 			& 0.36 			& 18.8 \\
	& Bisections 	& 11,823 		& 1,922,833 	& 55,177 		& 3,047 		& 240,945 \\
	& NE $f$		& 520,800 		& 108,944,542 	& 4,756,387 	& 342,780 		& 16,478,673 \\
	& NE $F$ (DE) 	& 809 			& 52 			& 141 			& 104 			& 64 \\
	& NE $F$ (IBC) 	& 673,529 		& 29,603,444 	& 980,929 		& 26,523  		& 2,049,552 \\
	& NE $\nabla F$ & 614,451 		& 19,361,735 	& 689,304 		& 5,860  		& 954,871 \\
	& $|\mathcal{Q}|_{max}$ &	142	& 60			& 43			& 9				& 47 \\
	& $\mathit{NP}$ & 60			& 50			& 50 			& 30			& 50 \\
	& $W$ 			& 0.7			& 0.7 			& 0.7 			& 0.7 			& 0.7 \\
	& $\mathit{CR}$ & 0 			& 0.9			& 0.4 			& 0.9 			& 0.5 \\
	& Bisection scheme 	& RR 		& RR 			& Largest 		& Largest 		& RR \\
	& Search strategy 		& MaxDist 		& MaxDist 		& MaxDist 		& MaxDist 		& MaxDist \\
	& $\eta$ (fixed-point ratio)	& 0 		& 0.8 			& 0.9 			& 0 		& 0.9 \\
	& Gain /IBC	(\%)& -				& 27.8			& 38.2			& 35.7  		& 75.2 \\
	\hline
	\end{tabular}
	\label{tab:multimodal-standalone}
\end{table}
\def\arraystretch{1}
\end{landscape}

Note that the \gls{DE} component of Charibde keeps running until the proof of optimality is achieved by the \gls{IBC}. The total number of evaluations of the objective function "NE $f$ (DE)" may thus be much larger than the number of evaluations required to reach the global minimum.

The standalone components of Charibde exhibit various behaviors on the test problems. The \gls{DE} suffers from premature convergence on the Michalewicz and Rana functions, and achieves the global minimum on the Sine Envelope Sine Wave, Eggholder and Keane problems (although with no proof of optimality). The \gls{IBC} converges within a reasonable time on four out of five test problems and times out on the Michalewicz function.
On the other hand, Charibde converges in 6.3s on the Michalewicz function. Although a suboptimal upper bound is produced by the \gls{DE}, the pruning is improved and a large portion of the search space is discarded. In return, the \gls{IBC} injects new individuals into the \gls{DE} population whenever it improves the best known solution. Overall, Charibde surpasses the standalone \gls{IBC}. The results on the other problems show that Charibde consistently improves upon the computing time of the standalone \gls{IBC}: the line "Gain /IBC (\%)" indicates the relative gain of Charibde with respect to the standalone \gls{IBC}. On the given instances (excepted Michalewicz), it varies between 27.8 and 75.2\%. Our search strategy MaxDist maintains a very low maximum size $|\mathcal{Q}|_{max}$ of the priority queue $\mathcal{Q}$ in Charibde (between 9 and 142 boxes), compared to the standalone \gls{IBC} (between 4,605 and 870,564 boxes). The domain reduction strategy described in Section \ref{sec:de-domain-reduction} can thus be applied at a low cost.

\begin{figure}[htbp!]
\centering
\includegraphics[width=\textwidth]{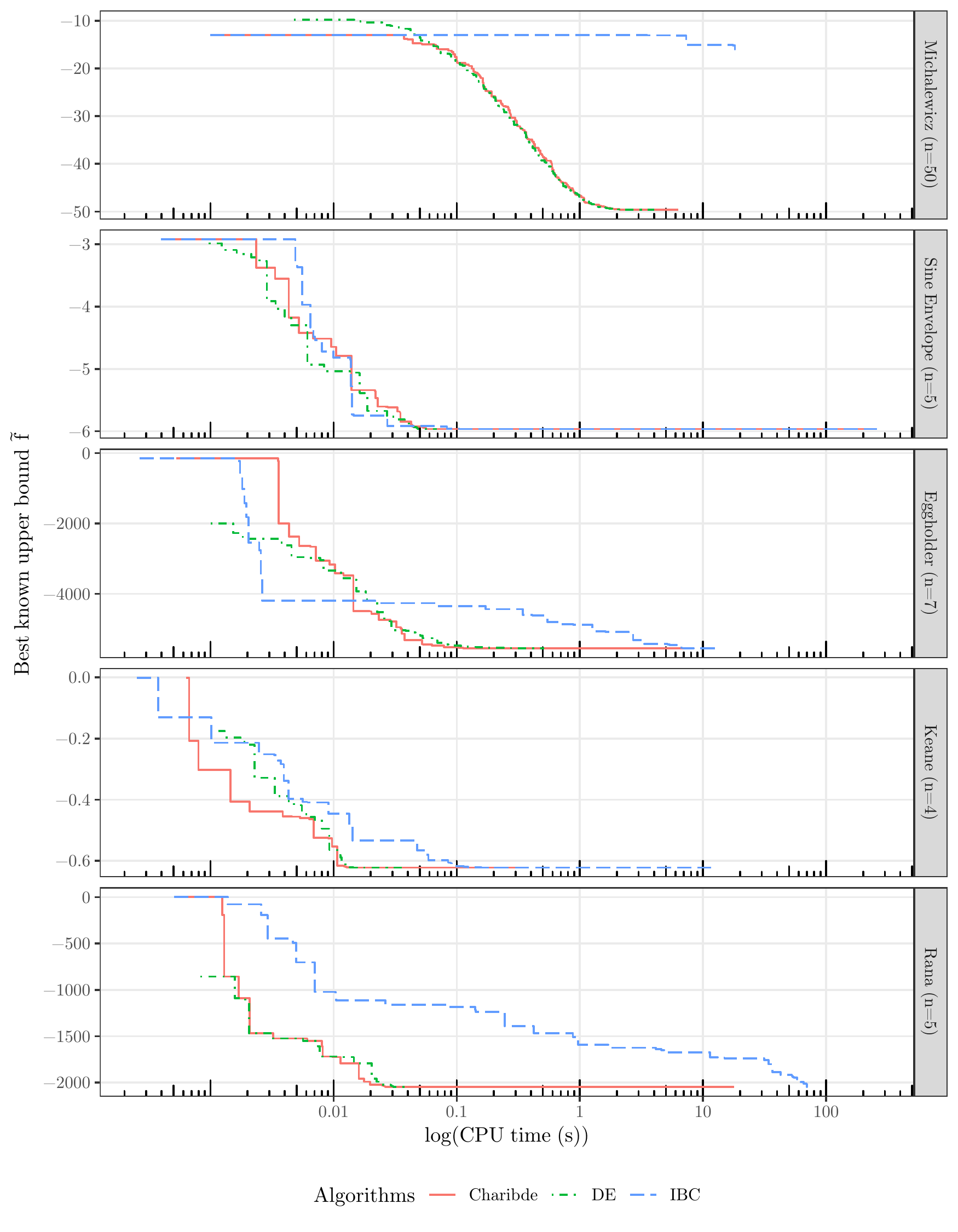}
\caption{Best known upper bound against computation time for Charibde and its standalone components on the multimodal test problems}
\label{fig:multimodal-standalone}
\end{figure}

Figure \ref{fig:multimodal-standalone} illustrates the evolution of the best known upper bound $\tilde{f}$ against the computation time of Charibde and its standalone components for particular instances of the multimodal test problems. It confirms that a standard \gls{IBC} algorithm, even equipped with the usual "best-first search" strategy, struggles to quickly find a good upper bound of the global minimum.

Even though Charibde outperforms the standalone \gls{IBC} (this trend generally intensifies on larger instances), the benefits of the hybridization seem less obvious than on a nonlinearly constrained benchmark (Section \ref{sec:coconut-comparison}). The test problems in the present section are highly multimodal (the functions involve trigonometric terms) and suffer from severe dependency ; a good upper bound of the global minimum is generally not sufficient to prune the search space efficiently without resorting to partitioning. On the contrary, nonlinearly constrained problems may be easier to solve, since they are often less subject to dependency ; the core challenge is often to find a feasible point, in which case the benefits of \gls{EA} are evident.

\clearpage

\section{Solver comparison}
\label{sec:multimodal-comparison}

Table \ref{tab:multimodal-solvers} presents a comparison of Charibde against six solvers (Ipopt, LOQO, MINOS, PGAPack, PSwarm and Couenne) -- available on the NEOS server~\cite{Gropp1997Optimization} -- on a particular instance of each of the test problems. When available, the number of evaluations of the objective function or the computation time is shown under the achieved objective value. Note that the solver BARON~\cite{Sahinidis1996BARON} does not support trigonometric functions.

\begin{table}[htbp]
	\centering
	\footnotesize
	\caption{Comparison of Charibde against six solvers on the multimodal test problems}
	\begin{tabular}{|l|c|c|c|c|c|}
	\hline
					& Michalewicz 	& Sine Envelope & Eggholder  	& Rana  		& Keane \\
					& ($n = 50$)	& ($n = 5$)		& ($n = 7$)		&($n = 5$)		& ($n = 4$) \\
	\hline
	Ipopt 			& -19.773742 	& -5.8351843 	& -5199.5967304 & -75.512076	& -0.2010427 \\
					& (167 eval) 	& (24 eval) 	& (7 eval) 		& (16 eval) 	& (8 eval) \\
	LOQO 			& -0.0048572 	& -5.8351843 	& -44.45892854 	& -69.5206		& -0.0983083 \\
					& (88 eval) 	& (17 eval) 	& (5406 eval)	& (138 eval)	& (50 eval) \\
	MINOS 			& 0 			& -5.87878 		& -5199.59673 	& -233.592		& -0.2347459 \\
					& (3 eval) 		& (38 eval) 	& (3 eval)		& (1 eval)		& (3 eval) \\
	PGAPack 		& -37.60465 	& -5.569544 	& -4369.204 	& -2091.068		& - \\
					& (9582 eval) 	& (9615 eval) 	& (9593 eval) 	& (9622 eval)	& - \\
	PSwarm 			& -24.38158 	& -5.835182 	& -3429.485 	& -1595.056		& - \\
					& (2035 eval) 	& (2049 eval) 	& (2054 eval) 	& (2046 eval)	& - \\
	Couenne & -49.6\underline{19042} & -5.96\underline{60007} & -55\underline{10.513933430007}
& \bf{-2046.8320657} & -0.6222\underline{999} \\
					& (265s) 		& (0.4s) 		& (44s) 		& (20.3s) 		& (2s) \\
	\hline
	\bf{Charibde} 	& \bf{-49.6248323} & \bf{-5.9659811} & \bf{-5548.9775483} & \bf{-2046.8320657} 
& \bf{-0.6222810} \\
					& (4.9s) 		& (194.9s)		& (7.85s)		& (18.8s) 		& (0.36s) \\
	\hline
	\end{tabular}
	\label{tab:multimodal-solvers}
\end{table}

To illustrate the multimodality of the test problems, three local solvers (Ipopt, LOQO, MINOS) were chosen. They usually require few iterations to converge towards a local minimum, starting from an initial point ; the quality of the local minimum depends on the initial point and the size of the basins of attraction (the sets of initial points that converge towards a given minimum). The three solvers produce poor local minima for the considered problems. Two metaheuristics for bound constrained optimization, PGAPack (\gls{GA}) and PSwarm (Particle Swarm Optimization), were also included in the benchmark. Overall, PGAPack produces better local minima than PSwarm, albeit at a higher cost. Their hyperparameters were set to default values ; the numerical results could probably be improved by picking more appropriate values of the hyperparameters.
Couenne is a deterministic global optimization solver that performs a comprehensive exploration of the search space ; it implements a spatial branch and bound algorithm based on reformulation techniques that constructs a linear programming relaxation in each subspace. Although considered as one of the most efficient solvers nowadays, Couenne is not rigorous: the underapproximations and overapproximations obtained by relaxing the objective function and the constraints are not conservative and are subject to numerical approximations. Therefore, optimality cannot be reliably proven ; Table \ref{tab:multimodal-solvers} demonstrates that the optimal results for the Michalewicz, Sine Envelope Sine Wave, Eggholder and Keane problems achieved by Couenne are erroneous (incorrect digits are underlined).

These results suggest that Charibde is highly competitive with respect to Couenne in terms of computation time (except for the Sine Envelope Sine Wave function), while providing a numerical proof of optimality with a precision $\varepsilon$. Charibde converges faster that Couenne on the Michalewicz (time ratio 54), Eggholder (5.6), Rana (1.1) and Keane (5.6) problems. Couenne converges faster than Charibde on the Sine Envelope Sine Wave function (time ratio 487), however it provides an upper bound that is too low and the third digit is incorrect (-5.9660007 instead of -5.9659811). Consequently, Couenne performs a non reliable pruning of the search space and quickly converges towards an incorrect solution.

\begin{figure}[htbp]
\centering
\includegraphics[width=\textwidth]{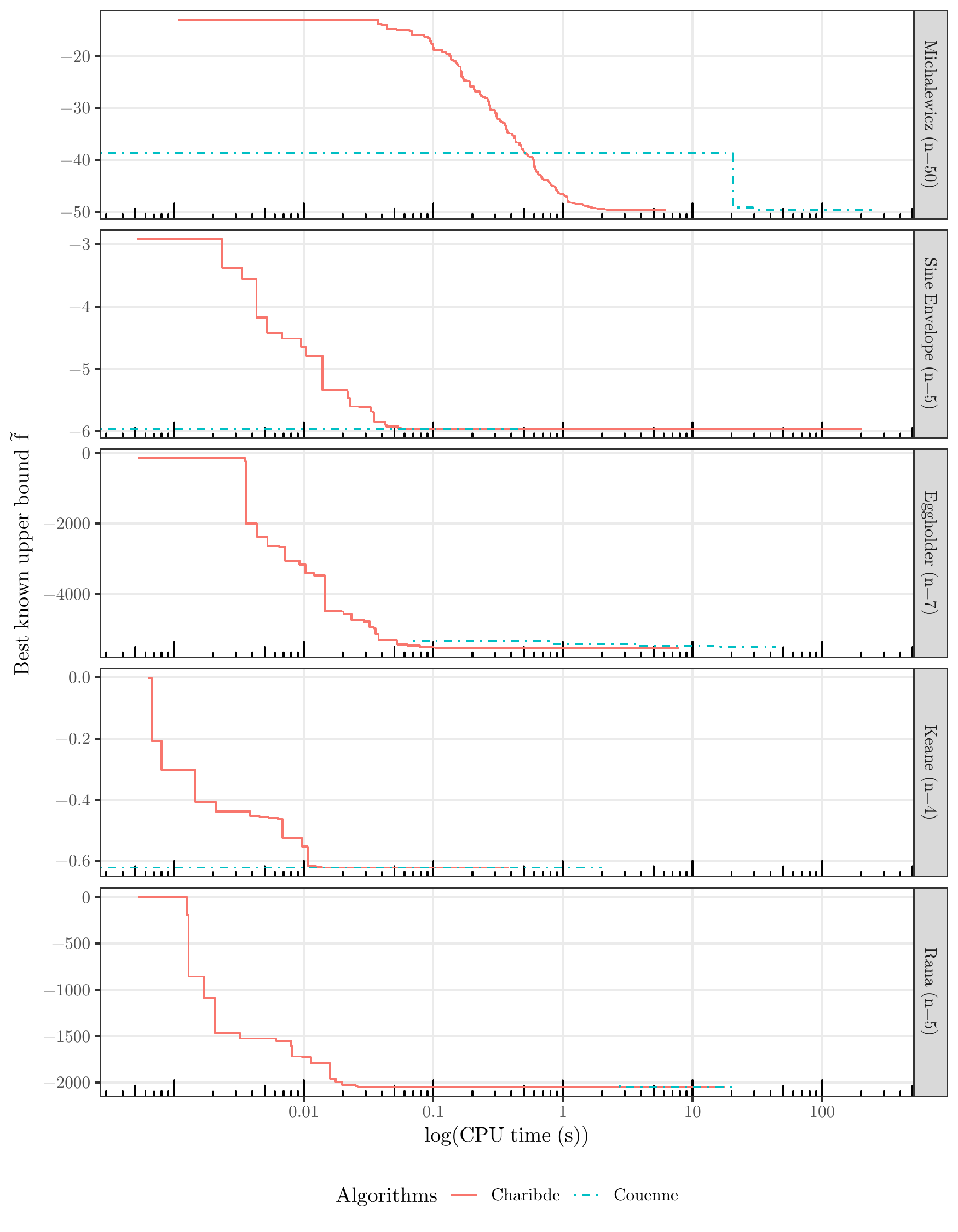}
\caption{Best known upper bound against computation time for Charibde and Couenne on the multimodal test problems}
\label{fig:multimodal-charibde-couenne}
\end{figure}

Figure \ref{fig:multimodal-charibde-couenne} presents the best known upper bound of Charibde and Couenne against the computation time (in logarithmic scale). Couenne benefits from an efficient cooperation between the spatial branch and bound algorithm and the nonlinear solver Ipopt, invoked in order to quickly compute a good upper bound of the global minimum.

\chapter{Rigorous proof of optimality in molecular dynamics}
\label{chap:lj}

\minitoc

\section{Lennard-Jones potential}
The Lennard-Jones potential is a simplified model proposed by~\cite{Jones1924} that
approximates the interaction between a pair of spherical atoms. This model is accurate
for noble gases in which atoms repel each other at short range and attract each other
at long range. The Lennard-Jones potential is given by:
\begin{equation}
V(d_{ij}) = 4\epsilon \left[ \left( \frac{\sigma}{d_{ij}} \right)^{12} -
	\left( \frac{\sigma}{d_{ij}} \right)^6 \right]
\label{eq:lj-potentiel-general}
\end{equation}
where $d_{ij}$ is the distance (in ångströms) between atoms $i$ and $j$,
$\epsilon$ is the depth of the potential well and $\sigma$ is the distance at
which the potential is zero. $d_{min} \eqdef 2^{1/6} \sigma$ is the distance at
which the potential reaches its minimum:
\begin{itemize}
\item when $d_{ij} > d_{min}$, the attractive Van der Waals forces, modeled by the term
$(\frac{\sigma}{d_{ij}})^6$, prevail over the repulsive forces ;
\item when $d_{ij} < d_{min}$, the repulsive forces, modeled by the approximate term
$(\frac{\sigma}{d_{ij}})^{12}$, prevail over the attractive forces.
\end{itemize}

\begin{figure}[htbp]
\centering
\def\svgwidth{0.7\columnwidth}
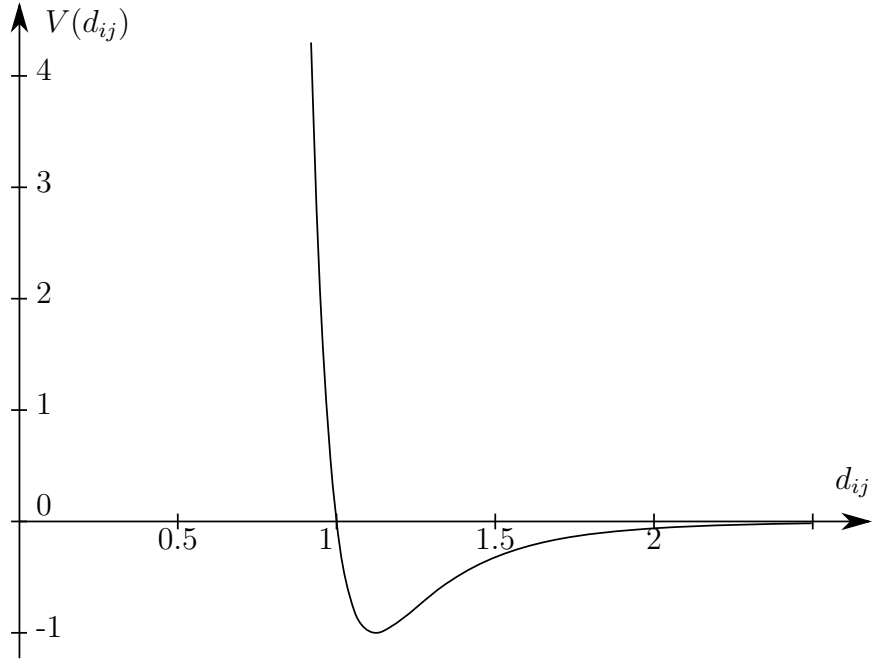
\caption{Reduced Lennard-Jones potential}
\label{fig:lj}
\end{figure}

The expression of the Lennard-Jones potential is generally given with
reduced units $\epsilon = 1eV$ and $\sigma = 1\AA$ (Figure \ref{fig:lj}):
\begin{equation}
V(d_{ij}) = 4 \left( \frac{1}{d_{ij}^{12}} - \frac{1}{d_{ij}^6} \right)
\label{eq:lj-potentiel}
\end{equation}

\section{Lennard-Jones cluster problem}
\subsection{Minimum-energy spatial configuration}
Finding the minimum-energy spatial configuration of a cluster of noble gas
with $N$ atoms amounts to minimizing the pairwise interactions within the cluster:
\begin{equation}
f_N(\bm{x}) = \sum_{i < j}^N V(d_{ij}) = 4 \sum_{i < j}^N \left( \frac{1}{d_{ij}(\bm{x})^{12}} - 
\frac{1}{d_{ij}(\bm{x})^6} \right)
\label{eq:lj}
\end{equation}
where $(x_i, y_i, z_i)$ are the Cartesian coordinates of atom $i$,
$\bm{x} = (x_1, y_1, z_1, \ldots, x_N, y_N, z_N)$ is the vector of $3N$ variables
and $d_{ij} > 0$ is the distance between atoms $i$ and $j$:
\begin{equation}
d_{ij}(\bm{x})^2 = (x_i - x_j)^2 + (y_i - y_j)^2 + (z_i - z_j)^2
\end{equation}

\subsection{An open problem}

Beneath the deceiving simplicity of $f_N$ (Equation \ref{eq:lj}) lies a nonconvex and
highly multimodal problem ; the numerical tests of~\cite{Locatelli2003} suggest that
the number of local minima of $f_N$ increases exponentially with $N$.
In optimal configurations with $N = 2$, $N = 3$ and $N = 4$, the atoms are vertices
of a regular tetrahedron. Configurations with $N \ge 5$ however have never been solved using
reliable global methods~\cite{Vavasis1994}.
Numerous putative solutions\footnote{\url{http://doye.chem.ox.ac.uk/jon/structures/LJ.html}} were
obtained using a variety of approximate methods~\cite{Northby1987,Hoare1971,Leary1997,Wales1997}.

\section{The first rigorous proof of optimality for 5 atoms}

In this section, we present the first rigorous proof of optimality for the Lennard-Jones cluster with 5 atoms achieved by Charibde \cite{Vanaret2015Premiere}. The best known solution in the literature is a triangular bipyramid~\cite{Sloane1995} whose objective value if $-9.103852415708$.

We demonstrate that the best known solution, which has never been numerically certified, is indeed optimal. We explain how to reformulate the problem in order to reduce the overestimation computed by \gls{IA} and decrease the size of the optimization problem. Lastly, we exhibit the first numerical proof of optimality for the open Lennard-Jones cluster problem with 5 atoms.

\subsection{Dependency reduction}
Equation \ref{eq:lj-potentiel} contains two occurrences of the distance $d_{ij}$. A well-known trick to reduce dependency is to complete the square:
\begin{equation}
V(d_{ij}) = 4 \left( \frac{1}{d_{ij}^{6}} - \frac{1}{2} \right)^2 - 1
\label{eq:lj-potentiel-factorise}
\end{equation}
The natural inclusion function of the potential $V$ becomes optimal with respect to the occurrences of $d_{ij}$. Note however that the objective function (Equation \ref{eq:lj}) still suffers from dependency, since the coordinates $(x_i, y_i, z_i)$ of the atoms have multiple occurrences in the distance terms $d_{ij}(\bm{x})$.

\subsection{Symmetry breaking}
Equation \ref{eq:lj} involves pairwise distance terms, therefore the optimal solutions of the
cluster problem are invariant under translation and rotation. The coordinates of some of the
atoms should be fixed in order to partially break symmetries and reduce the size of the search
space.

We fix the first atom at the origin of the coordinate system, the second atom on the
half-line $x \ge 0$, the third atom in the first quadrant of the plane $z = 0$ and the fourth
atom in the first octant:
\begin{equation}
\begin{cases}
x_1 = y_1 = z_1 = 0 \\
x_2 \ge 0, y_2 = z_2 = 0 \\
x_3 \ge 0, y_3 \ge 0, z_3 = 0 \\
x_4 \ge 0, y_4 \ge 0, z_4 \ge 0
\end{cases}
\end{equation}

Reducing the symmetries thus decreases the size of the problem with 5 atoms from 15 to 9
variables, and reduces the domains of the remaining variables.

\subsection{Proof of optimality}
\begin{figure}[h!]
	\centering
	\includegraphics[width=0.35\columnwidth]{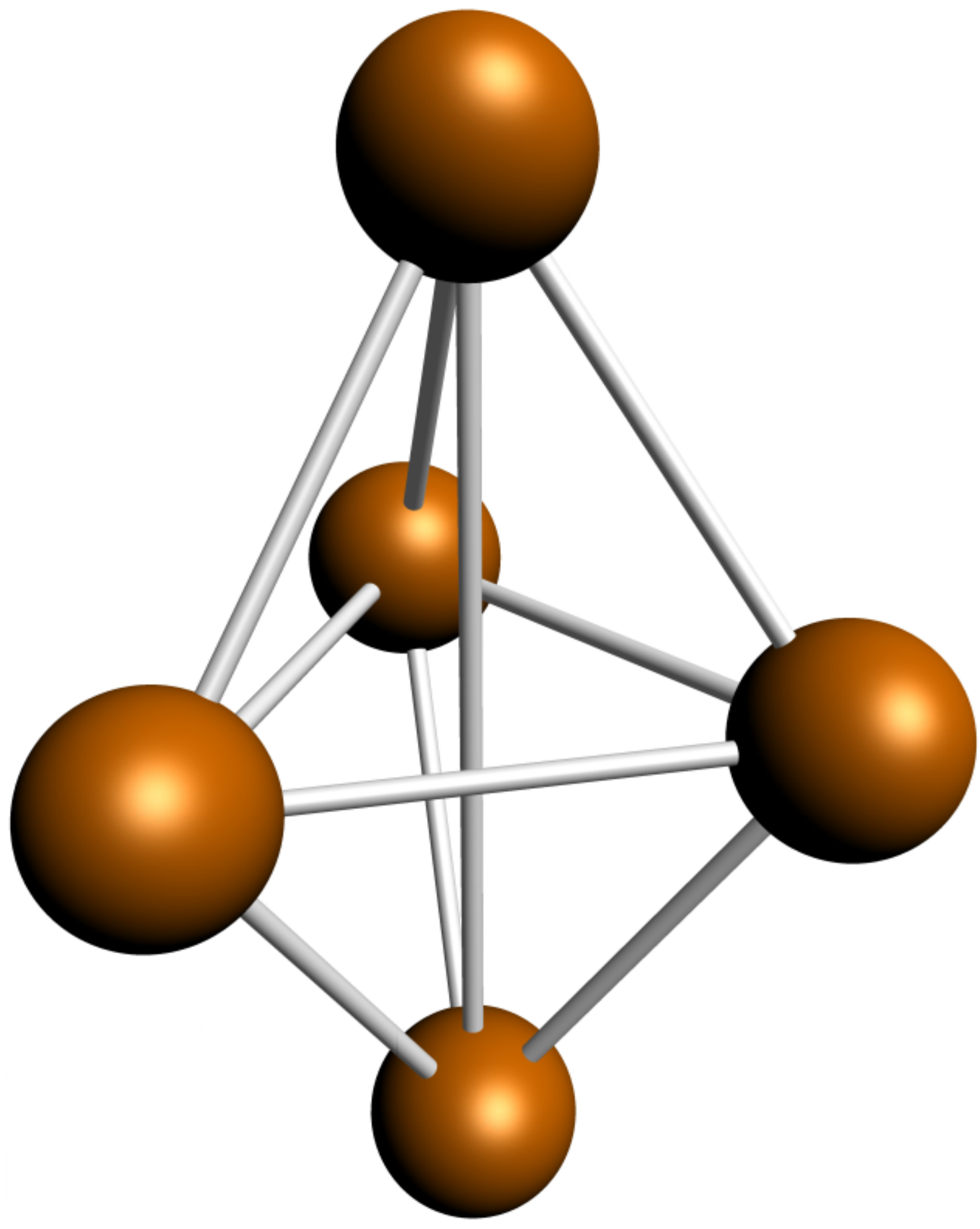}
	\caption{Lennard-Jones problem with 5 atoms: optimal configuration}
	\label{fig:lj5}
\end{figure}

Charibde proved the optimality of the best known solution over the initial domain
$(x_i, y_i, z_i) \in [-1.2, 1.2]$, and reached the global minimum $f_5^* = -9.103852415707552$
with a precision $\varepsilon = 10^{-9}$. The corresponding solution is given in Table \ref{tab:lj5}
and is represented in Figure \ref{fig:lj5}. The hyperparameters of Charibde
can be found in Table \ref{tab:lj5-param}.

\vspace{1cm}

\begin{table}[htbp]
	\centering
	\caption{Lennard-Jones problem with 5 atoms: optimal solution}
	\begin{tabular}{|cccc|}
	\hline
	Atom & $x$ & $y$ & $z$ \\
	\hline
	1 & 0 & 0 & 0 \\
	2 & 1.1240936 & 0 & 0 \\
	3 & 0.5620468 & 0.9734936 & 0 \\
	4 & 0.5620468 & 0.3244979 & 0.9129386 \\
	5 & 0.5620468 & 0.3244979 & -0.9129385 \\
	\hline
	\end{tabular}
	\label{tab:lj5}
\end{table}

\begin{table}[htbp]
	\centering
	\caption{Lennard-Jones problem with 5 atoms: hyperparameters of Charibde}
	\begin{tabular}{|lc|}
	\hline
	Hyperparameter	 						& Value \\
	\hline
	$\varepsilon$ (precision) 				& $10^{-9}$ \\
	$\mathit{NP}$ (\gls{DE} population size) & $40$ \\
	$W$ (amplitude factor)	 				& $0.7$ \\
	$\mathit{CR}$ (crossover rate)	 		& $0.4$ \\
	Bisection strategy						& Largest \\
	Priority								& MaxDist \\
	$\eta$ (fixed-point ratio)				& 0 \\
	\hline
	\end{tabular}
	\label{tab:lj5-param}
\end{table}

The average computation times, queue size and numbers of evaluations of the objective function
and its gradient over 100 runs are gathered in Table \ref{tab:lj5-results}. The \gls{DE}
reaches the global optimum $f_5^*$ after $764$ iterations ($0.11$s), whereas the proof of
optimality is achieved only after $1436$s.
\begin{table}[htbp]
	\centering
	\caption{Lennard-Jones problem with 5 atoms: average results of Charibde over 100 runs}
	\begin{tabular}{|lc|}
	\hline
	Metrics & Value \\
	\hline
	Average CPU time (s) & 1,436 \\
	Maximum CPU time (s) & 1,800 \\
	Maximum size of $\mathcal{Q}$ & 46 \\
	Number of evaluations of $F$ (\gls{IBC}) & 7,088,758 \\
	Number of evaluations of $\nabla F$ (\gls{IBC}) & 78,229,737 \\
	Number of evaluations of $f$ (\gls{DE}) & 483,642,320 \\
	Number of evaluations of $F$ (\gls{DE}) & 132 \\
	\hline
	\end{tabular}
	\label{tab:lj5-results}
\end{table}


\clearpage

\subsection{Comparison with state-of-the-art solvers}
Table \ref{tab:comparaison-lj} compares Charibde against state-of-the-art nonlinear solvers
BARON and Couenne on the reformulated problem.
The search time indicates the time required to find the global minimum, and the proof time
indicates the total time required to prove its optimality.

These results corroborate the fact that state-of-the-art nonlinear solvers, while exhibiting
excellent convergence times, cannot guarantee the optimality of the solution in the presence
of roundoff errors. Incorrect digits in Table \ref{tab:comparaison-lj} are underlined.
The solution found by BARON is incorrect after the fifth digit, while that of Couenne is
incorrect after the fourth digit, although the required tolerance
is $\varepsilon = 10^{-9}$.
Moreover, while Couenne finds a global minimum, BARON cannot prove the global optimality of
its solution.

\begin{table}[htbp]
	\centering
	\caption{Lennard-Jones problem with 5 atoms: comparison of BARON, Couenne and Charibde}
	\begin{tabular}{|l|c|c|c|}
	\hline
						& BARON  	& Couenne	& Charibde \\
	\hline
	Minimum 			& -9.10385\underline{346444055} & -9.1038\underline{70325603582} &
-9.103852415707552 \\
	Search time (s) 	& 0.23 		& 41.94		& 0.11 \\
	Proof time (s)		& 0.23		& 61.7		& 1436 \\
	Status				& locally optimal	& optimal	& certified ($\varepsilon = 10^{-9}$) \\
	\hline
	\end{tabular}
	\label{tab:comparaison-lj}
\end{table}

\chapter*{General conclusion}
\addstarredchapter{General conclusion} 

In this document, we introduced a new cooperative framework that combines interval methods and evolutionary algorithms in a parallel fashion. The interval methods alternates partitioning of the search space and filtering of inconsistent values, and converges towards the global minimum, even in the presence of roundoff errors. We devised operators that help the evolutionary algorithm escape local minima, and accelerate the interval method.

\section*{Contributions}
Inspired by~\cite{Alliot2012Finding}, our reliable hybrid solver Charibde combines interval methods and evolutionary algorithms, and certifies the global optimality of the solution with a given precision. The differential evolution algorithm quickly explores the search space in the search of a satisfactory feasible solution, then updates the best solution of the interval branch and contract algorithm, a framework that alternates between partitioning of the search space and filtering of inconsistent values. Maintaining the best known upper bound of the global minimum allows the pruning of infeasible or suboptimal subspaces through the use of powerful refutation operators. In return, new solutions are injected into the population of the evolutionary algorithm whenever the interval methods identify promising regions. This avoids premature convergence towards local minima. The remaining subspaces to be processed by the interval branch and contract algorithm may be exploited periodically by the differential evolution algorithm ; its population is reinitialized within the convex hull of the remaining subspaces, thus avoiding the exploration of subspaces that have already been discarded. A novel exploration strategy implemented within the branch and contract algorithm, based on the maximum distance between the current solution and the remaining subspaces, tends to increase the efficiency of the convex hull operation.

Although they are not used in our work, iterative local optimization methods are good candidates for computing good feasible points and to provide the interval methods with upper bounds of the global minimum. Since all optimization problems tackled in this document contain analytical functions whose first and second derivatives are available,  Newton-based methods (such as SQP or interior point methods) may be invoked.
Our numerical results suggest however that the differential evolution algorithm often reaches the global minimum of the considered optimization problems. When it is not the case, the advanced cooperation techniques presented in Chapter \ref{chap:charibde} -- the periodic reduction of the domain of the evolutionary algorithm and the injection of fresh individuals into the population -- usually address the problem of premature convergence and guide the population towards the optimal solution.

Charibde has proven competitive with state-of-the-art solvers: it outperforms interval-based solvers GlobSol, IBBA and Ibex on a subset of 11 challenging COCONUT problems selected by \cite{Araya2012Contractor} by an order of magnitude. These results suggest that Charibde's strategy is to partition the search space more efficiently, while Ibex opts for aggressive filtering of the subspaces. We provide new optimality results for five scalable multimodal problems (Michalewicz, Sine Wave Sine Envelope, Eggholder, Keane, Rana) for which few solutions are known in the literature. We observe however that the impact of an excellent approximation of the global minimum provided by the evolutionary algorithm remains limited for highly multimodal problems.

\begin{figure*}[h!]
	\centering
	\includegraphics[width=0.35\columnwidth]{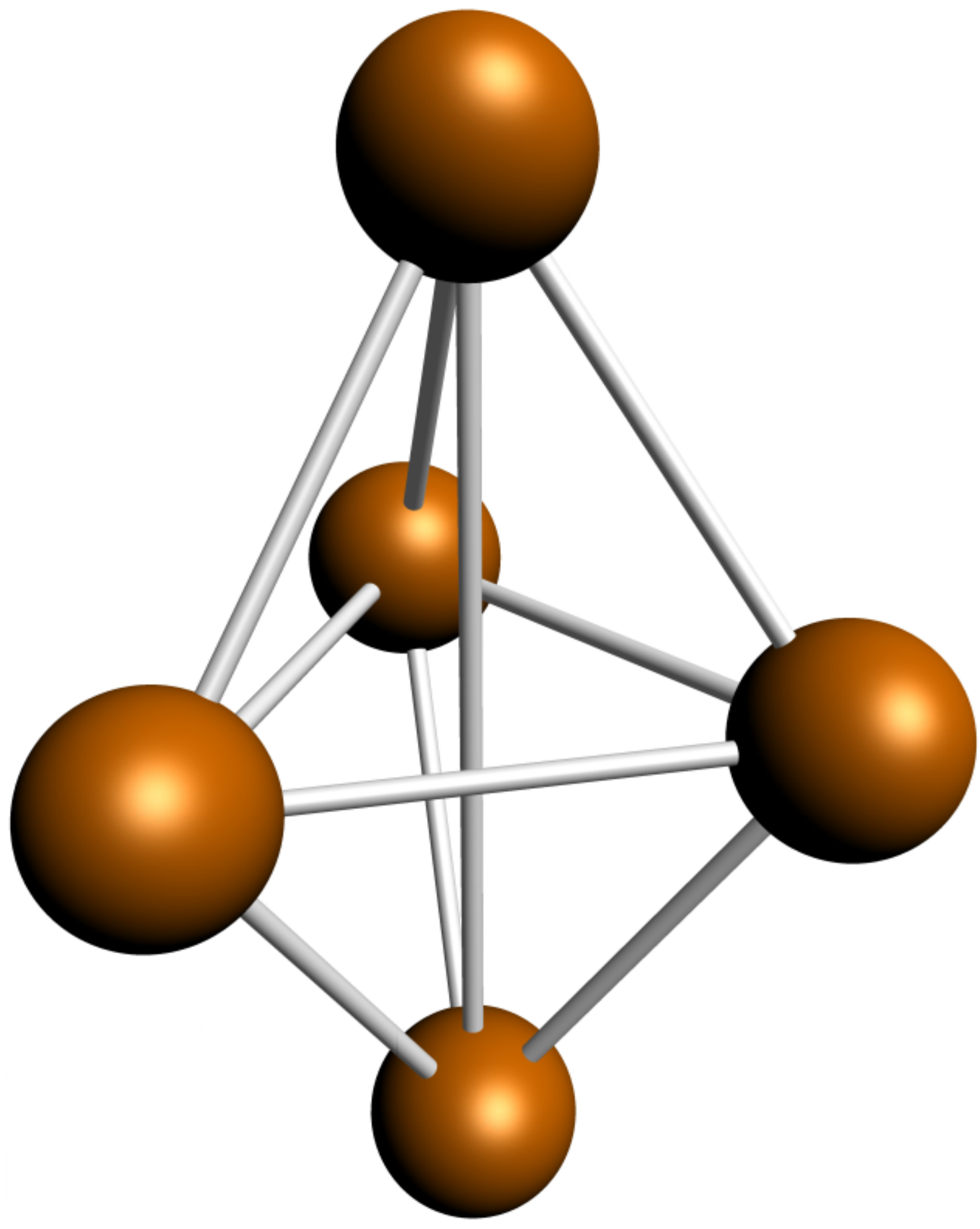}
	\captionsetup{labelformat=empty}
	\caption{Solution to the open Lennard-Jones cluster problem with 5 atoms}
	\label{fig:lj-conclusion}
\end{figure*}

Finally, we closed the open Lennard-Jones cluster problem with five atoms, a multimodal optimization problem stemming from molecular dynamics. Charibde provided the first numerical proof of global optimality of the solution with a precision of $10^{-9}$ (Figure \ref{fig:lj-conclusion}). We showed that the exhaustive (albeit unreliable) solvers BARON and Couenne produce numerically erroneous results that cannot be trusted. State-of-the-art interval-based solvers do not converge within reasonable time.

\newpage

\section*{Perspectives}

\subsection*{Inequality constrained optimization}
Using interval methods to handle the combinatorics of the complementarity equations inherent to inequality constrained optimization appears to be a promising field of research. Remember that:
\begin{enumerate}
\item the inactive inequality constraints have no influence on the solution ;
\item the inequality constraints active at the solution (provided that the active set is known) may be handled as equality constraints, whose filtering power is higher than that of inequality constraints.
\end{enumerate}
\cite{Hansen1992} mentions the possibility of solving the Fritz John conditions on the current box with an interval Newton method.
The difficulty of solving large preconditioned interval systems with respect to the primal and dual variables is probably the reason why little consideration has been given to this approach. However, an approximation of the active set (provided by solving a linear program) may help make the distinction between active and inactive inequality constraints in order to improve the filtering of inconsistent values.

\subsection*{Alternative enclosure methods}
Charibde's version of X-Newton~\cite{Araya2012Contractor} implements the recursive variant of \cite{Hansen1968Solving} to evaluates Taylor forms. Although the enclosures are usually tighter than the standard Taylor form, it computes the $n$ partial derivatives independently on different subboxes.
Computing interval slopes in adjoint mode~\cite{Krawczyk1985Slopes} generalizes Hansen's variant, while evaluating all $n$ partial derivatives simultaneously.

\begin{figure*}[htbp]
\centering
\scriptsize
\def\svgwidth{0.45\columnwidth}
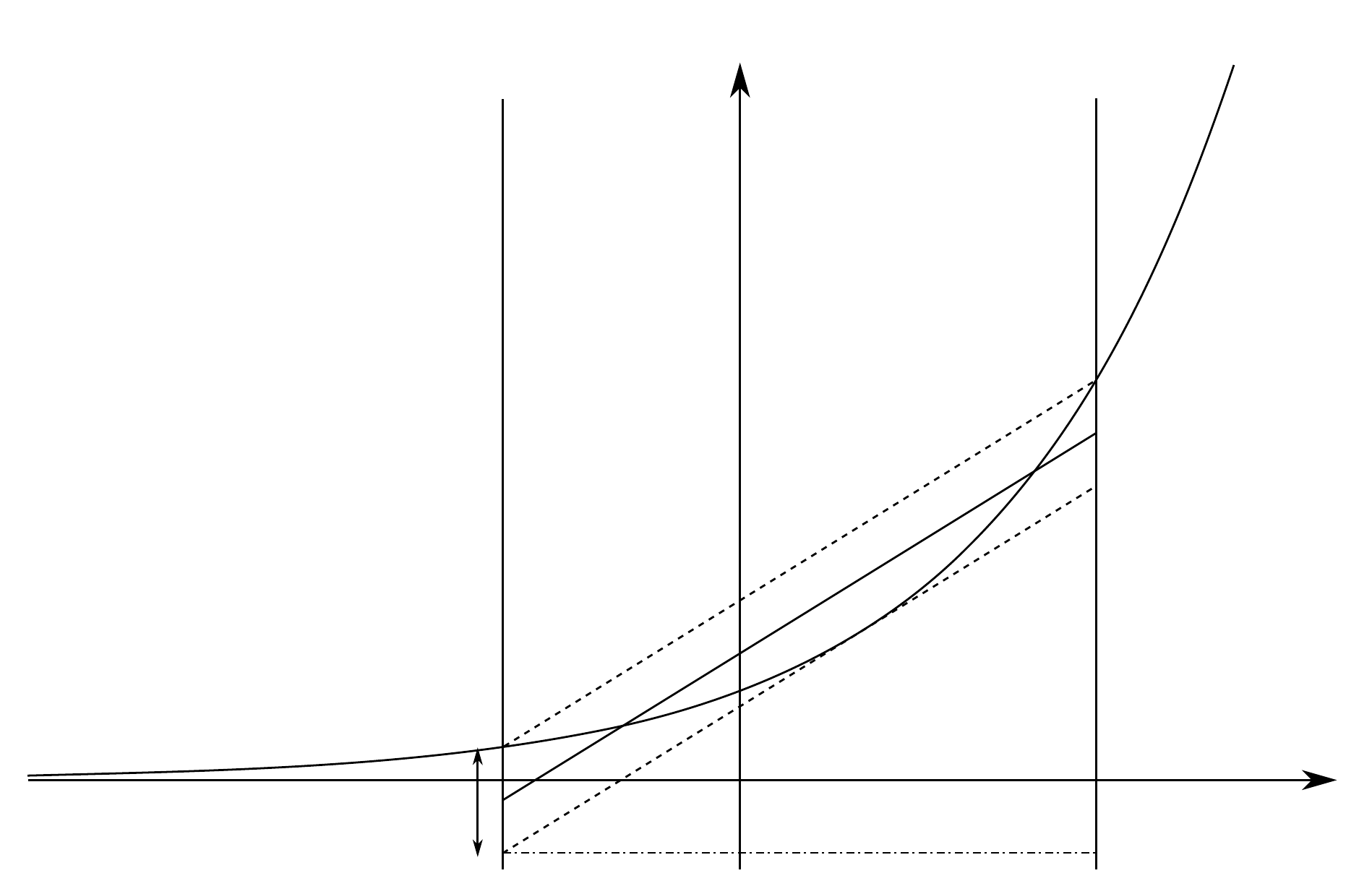
\hfill
\def\svgwidth{0.45\columnwidth}
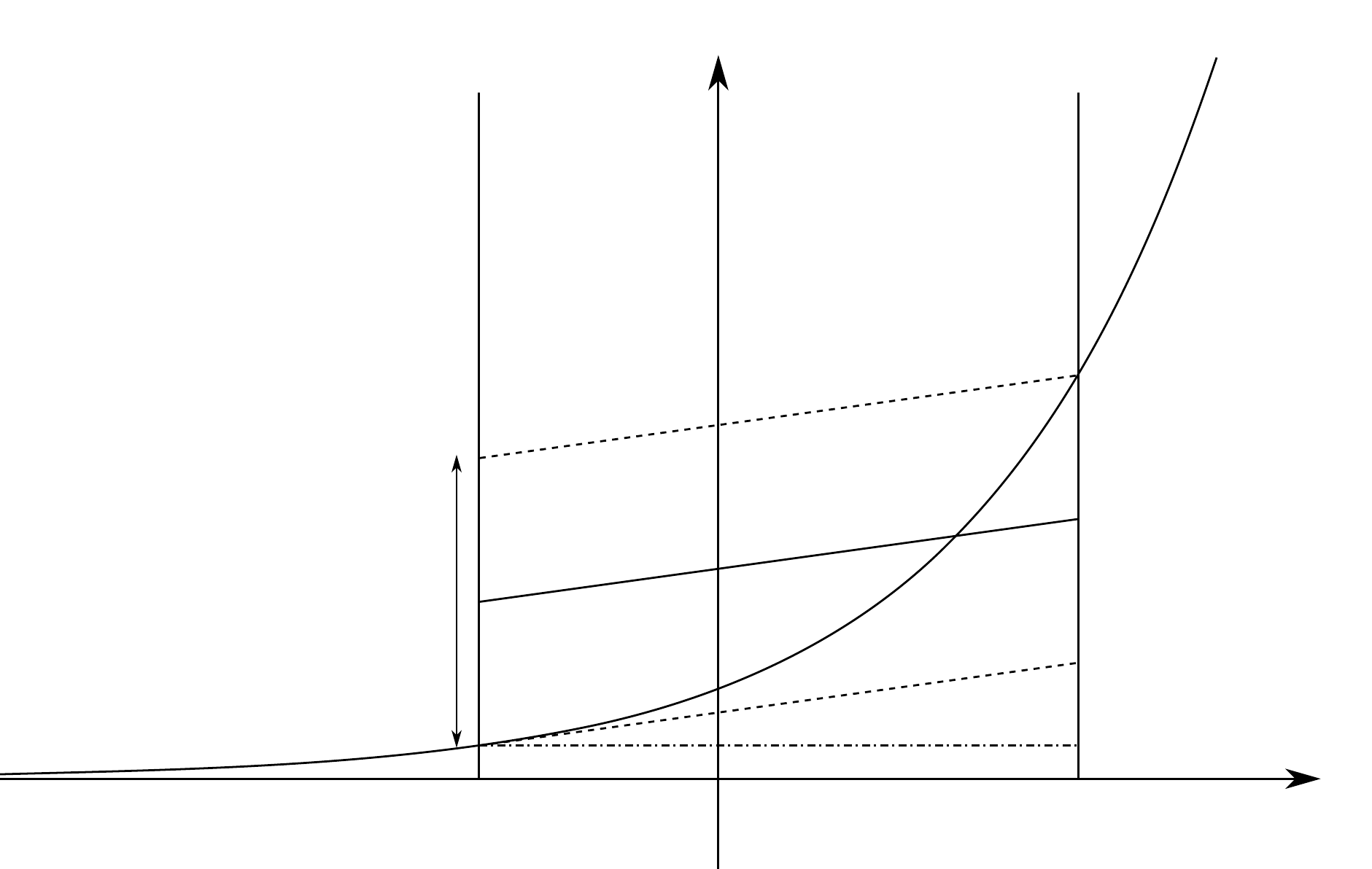
\captionsetup{labelformat=empty}
\caption[Affine approximations of $\exp$ using the Tchebychev and min-range methods]{Affine approximations of $f(x) = \exp(x)$ on $X = [-1, 1.5]$ using the Tchebychev (left) and min-range (right) methods}
\label{fig:affine-conclusion}
\end{figure*}

Affine arithmetic (Figure \ref{fig:affine-conclusion}) is an alternative enclosure method that keeps track of the linear dependencies between quantities and may reduce the dependency effect spectacularly~\cite{Comba1993Affine, Stolfi1997Self, Messine2002Extensions, Ninin2010These}.
Interval slopes and affine arithmetic, albeit complex to implement, would probably improve significantly the computations of lower bounds and the convexification-based contraction in Charibde.

\subsection*{Parallelization}
Disjoint subspaces processed by the interval branch and contract algorithms may be explored independently and in parallel. The exploration of the search space may thus be distributed over several slave processes, while the master process maintains the best known upper bound of the global minimum. Compelling speedups may be obtained if the workloads of the processes are well balanced.

\subsection*{Cutting-plane approach}

Convexification-based techniques (Sections \ref{sec:convexification} and \ref{sec:convexification-contraction}) compute a convex linear relaxation of the problem ; the objective value of the solution to the linear problem is a lower bound of the original problem.
In this section, we discuss an extension that further exploits the contracted intermediary nodes of the syntax tree when contraction and automatic differentiation are sequentially combined~\cite{Schichl2005Interval}.
In this sense, our suggestion is similar to the method of~\cite{Lebbah2005Efficient}: the constraints are decomposed into elementary constraints and a convex relaxation is computed for each intermediary node.
Example \ref{ex:hc4-ad-convexified} illustrates how adding constraints to the relaxed problem may improve the lower bound of the initial problem.
\begin{example}
\label{ex:hc4-ad-convexified}
Consider the following problem with a linear objective function and a nonlinear constraint:
\begin{equation*}
\begin{aligned}
\min_{(x, y) \in \bm{X}} 	\quad & f(x, y) = -x - 2y \\
\text{s.t.} 	\quad & g(x, y) = x + (x+y)^2 - 1 = 0
\end{aligned}
\end{equation*}
where $\bm{X} = X \times Y = [0, 5] \times [0, 5]$. Invoking HC4Revise on the constraint $g = 0$ reduces the bounds of $x$ and $y$: $X = Y = [0, 1]$ (see Example \ref{ex:hc4-ad}).

In order to determine a lower bound of the problem using convexification techniques, we must compute convex linear under- and overapproximations of $g$. Since $g = 0$ is equivalent to $\{g \le 0, -g \le 0\}$, $g$ can be enclosed between two extremal Taylor forms: one hyperplane underestimates $g$, the other underestimates $-g$:
\begin{equation*}
\begin{aligned}
g(\underline{X}, \underline{Y}) + L_1(x - \underline{X}) + L_2(y - \underline{Y}) & \le 0 \\
-g(\underline{X}, \underline{Y}) - U_1(x - \underline{X}) - U_2(y - \underline{Y}) & \le 0
\end{aligned}
\end{equation*}
where the partial derivatives of $g$ on $\bm{X}$ are evaluated after the top-down phase of HC4Revise:
\begin{equation*}
\begin{aligned}
\frac{\partial g}{\partial x}(x, y) & = 1 + 2(x+y) \in 1 + 2[0, 1] = [1, 3] =: [L_1, U_1] \\
\frac{\partial g}{\partial y}(x, y) & = 2(x+y) \in 2[0, 1] = [0, 2] =: [L_2, U_2]
\end{aligned}
\end{equation*}
The two hyperplanes reduce to:
\begin{equation*}
\begin{aligned}
\quad x & \le 1 \\
\quad -3x-2y & \le -1
\end{aligned}
\end{equation*}
The corresponding polytope is represented in Figure \ref{fig:polytope1}. The optimal solution to this convexified problem is $(1, 1)$, with objective value $f(1, 1) = -3$.

\begin{figure*}[htb]
\centering
\def\svgwidth{0.4\columnwidth}
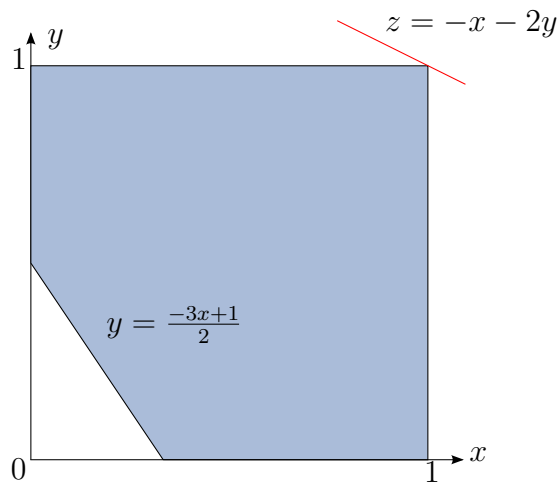
\captionsetup{labelformat=empty}
\caption{Polytope of the relaxed problem}
\label{fig:polytope1}
\end{figure*}

The information provided by the domain $[0, 1]$ of the node $x+y$ after the top-down phase of HC4Revise is not fully exploited: adding the explicit constraints $x + y \in [0, 1]$ to the initial set of constraints substantially reduces the size of the polytope (Figure \ref{fig:polytope2}) and improve the lower bound of the initial problem. The solution to the new problem is $(0, 1)$, with objective value $f(0, 1) = -2$, which is higher than $f(1, 1) = -3$.

\begin{figure*}[htb]
\centering
\def\svgwidth{0.4\columnwidth}
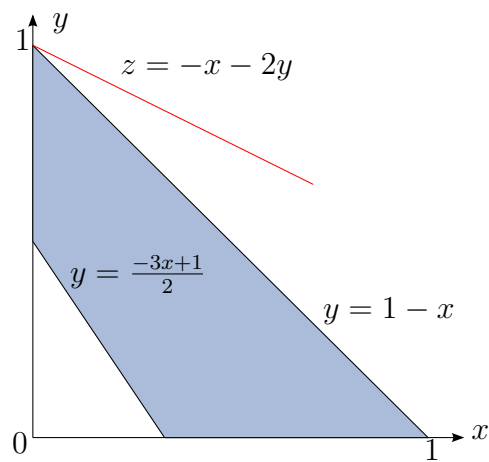
\captionsetup{labelformat=empty}
\caption{Adding constraints to the polytope of the relaxed problem}
\label{fig:polytope2}
\end{figure*}
\end{example}

\cleardoublepage
\phantomsection
\bibliographystyle{apalike}
\addcontentsline{toc}{chapter}{Bibliography}
\bibliography{Vanaret}

\end{document}